\documentclass{amsart}

\usepackage{tikz-cd}
\usepackage{booktabs,tabularx,array}
\usepackage{xcolor}
\usepackage[hyphens]{url}

\usepackage{amsthm,amsmath,amssymb,amsfonts,dsfont,amscd}
\usepackage[colorlinks, linkcolor=blue, citecolor=blue, urlcolor=blue]{hyperref}
\usepackage[indentafter]{titlesec}
\titleformat{name=\section}{}{\thetitle.}{0.8em}{\centering\scshape}
\usepackage{mathrsfs}
\usepackage{bm}
\numberwithin{equation}{section}
\titleformat{\subsection}[runin]
  {\normalfont\bfseries}{\thesubsection}{1em}{}
\makeatletter
\def\l@section{\@tocline{1}{0,2pt}{2pc}{8mm}{\ \ }}

\newtheorem{definition}{\textbf{Definition}}[section]
\newtheorem{theorem}[definition]{\textbf{Theorem}}
\newtheorem{corollary}[definition]{\textbf{Corollary}}
\newtheorem{lemma}[definition]{\textbf{Lemma}}
\newtheorem{proposition}[definition]{\textbf{Proposition}}

\theoremstyle{definition}
\newtheorem{example}[definition]{\textbf{Example}}
\theoremstyle{plain}
\newtheorem{remark}[definition]{\textbf{Remark}}
\newtheorem{letterthm}{Theorem}

\newtheorem{lettercor}[letterthm]{Corollary}

\def\wx{\infty}

\def\d{\mathrm{d}}

\def\prob{\mathrm{Prob}}
\def\BL{B(L^2(M))}

\def\CB{\mathcal{B}}
\def\CP{\mathcal{P}}
\def\CQ{\mathcal{Q}}
\def\CA{\mathcal{A}}
\def\har{\mathrm{Har}}
\def\id {\mathrm{id}}
\def\lag {\langle}
\def\rag {\rangle}

\def\N {\mathbb{N}}
\def\R {\mathbb{R}}
\def\C {\mathbb{C}}

\def\AD{\mathrm{Ad}}

\title[Poisson boundaries, bicentralizers, and Kadison's problem]{A Poisson boundary solution to Connes' bicentralizer conjecture and Kadison's problem on MASA}
\author{Shuoxing Zhou}
\address{Mathematical Institute, University of Oxford, Oxford, United Kingdom}
\email{shuoxing.zhou@maths.ox.ac.uk}
\thanks{The author was supported by the Engineering and Physical Sciences Research Council (EP/X026647/1). For the purpose of open access, the authors have applied a CC-BY license to any author accepted manuscript arising from this submission.}

\begin{document}

\begin{abstract}
We prove that, for every inclusion $M\subset P$ of countably decomposable von Neumann algebras with a normal faithful conditional expectation, $M\subset c(P)$ has the weak Dixmier property, where $c(P)$ is the continuous core of $P$. By Haagerup's characterization, we have a completely new proof of Connes' bicentralizer conjecture, independent of the recent proof of Houdayer--Marrakchi. By Marrakchi's characterization, we also prove his relative bicentralizer conjecture. We also give a new solution to Kadison's problem after Popa and Marrakchi: further assuming that $M$ and $P$ have separable preduals and $M'\cap P=\mathcal Z(M)$, then there is a maximal abelian subalgebra of $M$ that is also maximal abelian in $P$. The proofs are based on a Poisson boundary theory for general von Neumann algebras, developed here as an extension of the theory of Das--Peterson for finite von Neumann algebras. Once this theory is established, the proofs for Connes' bicentralizer conjecture and Kadison's problem are short (around 7 pages each) and self-contained, relying only on results established by 1990. In both proofs, no separate treatment is needed for the type $\mathrm{III}_1$ case.
\end{abstract}

\maketitle

\setcounter{tocdepth}{1}
\tableofcontents

\section{Introduction}\label{sec:introduction}

\subsection{Connes' bicentralizer conjecture and the weak Dixmier property.}

In the theory of von Neumann algebras, as in many areas of mathematics, averaging and averagable objects are a central theme of research. Amenability is the property that describes the possibility of averaging and is a central property in geometric and ergodic group theory. A discrete group $\Gamma$ is amenable if and only if there exists a $\Gamma$-invariant mean on $\ell^\infty(\Gamma)$ \cite{vonNeumann1929Masses}.

For von Neumann algebras, \textbf{injectivity} is the appropriate notion of \textbf{amenability}. A von Neumann algebra $M\subset B(\mathcal{H})$ is injective if and only if there exists a conditional expectation from $B(\mathcal{H})$ onto $M$. In his celebrated work, Connes \cite{Connes1976Injective} proved that, for von Neumann algebras with separable predual, injectivity is equivalent to being \textbf{hyperfinite}, which means that the algebra can be approximated by finite-dimensional subalgebras. For hyperfinite $M$, we can use these approximations to average over the compact unitary groups of finite-dimensional subalgebras. Passing to weak limits, we can average elements of $B(\mathcal{H})$ into $M'$. That is, $M\subset B(\mathcal{H})$ has the weak Dixmier property, which is an averaging property of inclusions of von Neumann algebras in the following sense: a \textbf{W$^*$-inclusion} $M\subset\CA$, i.e., an inclusion of von Neumann algebras, has the \textbf{weak Dixmier property} if every $x\in\CA$ can be averaged by unitaries of $M$ into the relative commutant $M'\cap\CA$, that is,
$$\overline{\operatorname{conv}}^{w}\{uxu^*\mid u\in\mathcal U(M)\}\cap(M'\cap\CA) \ne\emptyset\qquad(x\in\CA).$$
This connection between finite-dimensional approximation and averaging is a tool for the classification of amenable von Neumann algebras.

The classification of von Neumann algebras was a central goal of the work of Murray--von Neumann \cite{MurrayvonNeumann1936Rings}. Von Neumann algebras with separable predual decompose as direct integrals of \textbf{factors}, which are von Neumann algebras $M$ with trivial center, i.e., $\mathcal Z(M)=M'\cap M=\C1$. Murray--von Neumann divided factors into types $\mathrm{I}$, $\mathrm{II}_1$, $\mathrm{II}_\infty$, and $\mathrm{III}$. A type $\mathrm{I}$ factor is isomorphic to $B(\mathcal{H})$ for some Hilbert space $\mathcal{H}$. A type $\mathrm{II}_1$ factor admits a normal faithful finite trace, while a type $\mathrm{II}_\infty$ factor admits a normal faithful semifinite trace $\tau$ with $\tau(1)=\infty$. Type $\mathrm{III}$ factors admit no nonzero normal semifinite trace. Murray--von Neumann \cite{MurrayvonNeumann1943Rings} proved that, up to isomorphism, there is a unique hyperfinite $\mathrm{II}_1$ factor.

Connes \cite{Connes1973Classification} refined this classification by dividing type $\mathrm{III}$ into types $\mathrm{III}_\lambda$ for $0\leq\lambda\leq1$. He subsequently classified the amenable type $\mathrm{II}_1$, type $\mathrm{II}_\infty$, and type $\mathrm{III}_\lambda$ factors with separable predual for $0\leq\lambda<1$ \cite{Connes1976Injective}. For the remaining type $\mathrm{III}_1$ case, Connes \cite{Connes1985Factors} introduced the bicentralizer conjecture and showed that a positive solution for amenable factors would imply the uniqueness of the amenable type $\mathrm{III}_1$ factor. Let $M$ be a countably decomposable von Neumann algebra and let $\phi$ be a normal faithful state on $M$. The \textbf{bicentralizer} $\mathrm B(M,\phi)$ consists of the $x\in M$ such that $[x,a_n]\to0$ strongly for every bounded sequence $(a_n)\subset M$ satisfying $\|a_n\phi-\phi a_n\|\to0$, see Subsection~\ref{sec:centralizers}. Connes' bicentralizer conjecture states that $\mathrm B(M,\phi)=\C1$ for every type $\mathrm{III}_1$ factor $M$ and every normal faithful state $\phi$ on $M$.

The finite trace $\tau$ on a type $\mathrm{II}_1$ factor $M$ makes the analogous bicentralizer statement immediate: $\mathrm B(M,\tau)=\mathcal Z(M)=\C1$. More generally, much of the power of averaging in finite von Neumann algebras comes from the trace: it provides an invariant Hilbert-space geometry in which convexity produces fixed points. Type III factors have no such trace. For type $\mathrm{III}_\lambda$ factors with separable predual, $0<\lambda<1$, by the periodic modular structure we have a discrete decomposition as a crossed product of a type $\mathrm{II}_\infty$ factor by $\mathbb Z$ \cite{Connes1973Classification}. This reduces many structural questions to the semifinite setting. Such a periodic reduction is not available in general for type $\mathrm{III}_1$ factors. One must therefore find methods that work without either a trace on the original algebra or such a discrete reduction. This is one reason why the type $\mathrm{III}_1$ case often presents the greatest difficulty.

The same difficulty appears in Popa's problem on the weak Dixmier property (Corollary~\ref{lettercor expectation}) and Kadison's problem on MASA (Theorem~\ref{letterthm kadison}), discussed below. For inclusions $M\subset\CA$ with expectation and finite $M$, the weak Dixmier property follows from classical averaging arguments. The type $\mathrm{II}_1$ case of Kadison's problem was proved by Popa \cite{Popa1981Kadison} in 1981. The type $\mathrm{III}_1$ cases of these problems were resolved only much later, in 2019 and 2026, respectively.

Haagerup \cite[Theorem~2.3]{Haagerup1987Uniqueness} solved the bicentralizer problem for amenable factors, thereby proving the uniqueness of the amenable type $\mathrm{III}_1$ factor and completing the classification of amenable factors with separable predual. Haagerup \cite[Theorem~3.1]{Haagerup1987Uniqueness} proved that a type $\mathrm{III}_1$ factor $M$ with separable predual has trivial bicentralizer if and only if $M\subset c(M)$ has the weak Dixmier property, where $c(M)=M\rtimes_{\sigma^\phi}\R$ is the \textbf{continuous core} of $M$. Here $\sigma^\phi$ is the \textbf{modular automorphism group} associated with $\phi$, see Subsection~\ref{sec:core-weights} for background on modular theory and continuous cores. For amenable $M$, by hyperfiniteness we can average each element of $c(M)$ by unitaries of $M$ into $M'\cap c(M)$ in the weak topology, so $M\subset c(M)$ has the weak Dixmier property and hence $M$ has trivial bicentralizer.

Houdayer--Marrakchi \cite[Theorem~C]{HM26} recently solved Connes' bicentralizer conjecture for arbitrary type $\mathrm{III}_1$ factors. Their proof does not follow Haagerup's approach. It relies on the bicentralizer flow \cite{AHHM2020Bicentralizer}, on the ergodicity theorems for this flow proved by Marrakchi \cite[Theorem~D]{Marrakchi2020Full}, \cite[Theorem~C]{Ma26}, and on a resonance phenomenon, discovered in \cite{HM26}, for the joint action of the modular flow and the dual bicentralizer flow. In particular, their proof does not use Haagerup's theorem. It is therefore also a new proof that amenable type $\mathrm{III}_1$ factors have trivial bicentralizer, and hence, by \cite{Connes1985Factors}, a new proof of the uniqueness of the amenable type $\mathrm{III}_1$ factor.

In this paper, we follow Haagerup's approach to give a new proof of Connes' bicentralizer conjecture, independent of \cite{HM26}. We prove the weak Dixmier property for $M\subset c(M)$ directly, and more generally for $M\subset c(P)$ whenever $M\subset P$ admits a normal faithful conditional expectation, and then apply Haagerup's characterization. Our proof treats all factors uniformly, with only minor adjustments between the finite and properly infinite cases. In particular, the previously most difficult type $\mathrm{III}_1$ case is treated by the same argument as the other properly infinite cases, without using any recent advances in modular and bicentralizer theory.

\begin{letterthm}
\label{letterthm relative bicentralizer}
Let $M\subset P$ be an inclusion of countably decomposable von Neumann algebras with a normal faithful conditional expectation. Then $M\subset c(P)$ has the weak Dixmier property. In particular, $M\subset c(M)$ has the weak Dixmier property.
\end{letterthm}

By the result of Haagerup \cite[Theorem~3.1]{Haagerup1987Uniqueness}, Theorem~\ref{letterthm relative bicentralizer} implies Connes' bicentralizer conjecture.

\begin{lettercor}[Connes' bicentralizer conjecture {\cite{HM26}}]
\label{lettercor Connes bicentralizer}
Let $M$ be a type $\mathrm{III}_1$ factor with separable predual. Then $\mathrm B(M,\phi)=\C1$ for every normal faithful state $\phi$ on $M$.
\end{lettercor}

By \cite[Theorem~D]{Ma25}, Theorem~\ref{letterthm relative bicentralizer} also implies Marrakchi's relative bicentralizer conjecture \cite[Conjecture~C]{Ma25}, which generalizes Connes' bicentralizer conjecture. We refer to Subsection~\ref{sec:centralizers} for the definition of the relative bicentralizer and the statement of the relative bicentralizer conjecture.

\begin{lettercor}[The relative bicentralizer conjecture]
\label{lettercor relative bicentralizer}
Let $M\subset P$ be an inclusion of countably decomposable von Neumann algebras with a normal faithful conditional expectation. Then $M\subset P$ satisfies the relative bicentralizer conjecture.
\end{lettercor}

Our proof of Theorem~\ref{letterthm relative bicentralizer} is based on Poisson boundary techniques. To this end, we extend the Poisson boundary theory of Das--Peterson \cite{DP20} from finite von Neumann algebras to general von Neumann algebras with separable preduals. Subsection~\ref{sec:intro-weak-dixmier-pb} introduces this theory, starting from the averaging interpretation of Poisson boundaries of groups, and explains how we use it to prove the weak Dixmier property. Subsection~\ref{sec:intro-kadison} discusses another application of this theory, to Kadison's problem on MASA. Subsection~\ref{sec:intro-proof-strategy} explains the proof strategy for Theorem~\ref{letterthm relative bicentralizer}.

\subsection{The weak Dixmier property via Poisson boundaries.}\label{sec:intro-weak-dixmier-pb}

We first recall Poisson boundaries for groups, which are also objects defined through averaging. Let $\Gamma$ be a countable discrete group, and let $\mu\in\prob(\Gamma)$ be a generating probability measure, i.e., $\bigcup_{n\geq1}(\operatorname{supp}\mu)^n=\Gamma$. Define the Markov operator $P_\mu$ on $\ell^\infty(\Gamma)$ by
$$P_\mu(f)(g)=\sum_{s\in\Gamma}\mu(s)f(gs)\qquad(g\in\Gamma).$$
Hence $P_\mu$ can be viewed as a $\mu$-weighted averaging operator. Let $\har(\Gamma,\mu)=\{f\in\ell^\infty(\Gamma)\mid P_\mu(f)=f\}$ be the space of bounded $\mu$-harmonic functions.

The \textbf{Poisson boundary} of $(\Gamma,\mu)$ is the essentially unique nonsingular $\Gamma$-probability space $(B,\nu_B)$ such that the Poisson transform
$$P_{\nu_B}:L^\infty(B,\nu_B)\longrightarrow \har(\Gamma,\mu),\qquad P_{\nu_B}(F)(g)=\int_B F(gb)\,\d\nu_B(b)$$
is an isometric isomorphism \cite{Fu63a}.

For example, let $\Gamma$ be a non-elementary hyperbolic group with a finite generating set $S$, and let $\mu$ be a generating probability measure supported on $S$. Then its Poisson boundary can be identified with the Gromov boundary $(\partial_\infty\Gamma,\nu_\infty)$ of the Cayley graph with respect to $S$, where $\nu_\infty$ is the hitting measure of the $\mu$-random walk \cite[Theorems~7.6 and~7.7]{Kai00}.

Another example comes from lattices in semisimple Lie groups. Furstenberg \cite{Fu67} proved that, for every irreducible lattice $\Gamma<G$ in a connected semisimple real Lie group with finite center and no compact factors, there exists a fully supported probability measure $\mu$ on $\Gamma$ such that the $(\Gamma,\mu)$-Poisson boundary is $(G/P,\nu_P)$, where $P<G$ is a minimal parabolic subgroup. A major application of Poisson boundary theory is the proof of Margulis' superrigidity theorem for irreducible lattices in higher rank semisimple Lie groups \cite[Chapter~VII]{Mar91}.

We now turn to the corresponding construction for von Neumann algebras. Izumi \cite{Izu02,Izu04} developed a noncommutative Poisson boundary construction for normal ucp maps on von Neumann algebras, motivated by the study of discrete quantum groups. Building on Izumi's construction, Das--Peterson \cite{DP20} initiated Poisson boundary theory for finite von Neumann algebras, which was subsequently studied in \cite{Zh23,Zh23b,Zh24,Zh25}. We extend this theory to general von Neumann algebras $M$ with separable preduals. In this general setting, as analogues of generating probability measures on groups, our main objects are \textbf{normal, regular, strongly generating (nrsg) hyperstates} on $B(L^2(M))$ (Definition~\ref{def strong generation}). More concretely (see Example~\ref{exists nrsg hyperstate}), given any normal faithful state $\phi$ on $M$ with canonical cyclic vector $\xi_\phi\in L^2(M)$, we can always construct a normal state $\varphi$ on $B(L^2(M))$ with $\varphi|_M=\phi$ and sequences $(z_k),(y_k)\subset M$ such that
$$\overline{\C[(z_k)]}^{\,w}=\overline{\C[(y_k)]}^{\,w}=M,\qquad \sum_k z_k^*z_k=\sum_k y_ky_k^*=1.$$
The corresponding normal inner ucp maps
\begin{equation}\label{eq:intro-inner-ucp-maps}
\CP_\varphi(T)=\sum_k(Jz_k^*J)T(Jz_kJ),\qquad
\CQ_\varphi(T)=\sum_k y_kTy_k^*
\end{equation}
satisfy
\begin{equation}\label{eq:intro-hyperstate-representations}
\varphi(T)=\langle\CP_\varphi(T)\xi_\phi,\xi_\phi\rangle
=\langle\CQ_\varphi(T)\xi_\phi,\xi_\phi\rangle
\qquad(T\in B(L^2(M))).
\end{equation}
Here $J$ is the modular conjugation on $L^2(M)$, and $\C[(z_k)]$ and $\C[(y_k)]$ are the unital subalgebras (rather than the unital $*$-subalgebras) generated by $(z_k)$ and $(y_k)$, respectively.

Let $\har(\CP_\varphi)=\{T\in B(L^2(M))\mid \CP_\varphi(T)=T\}$. We define the \textbf{Poisson boundary} $\CB_\varphi$ to be the von Neumann algebra isomorphic to $\har(\CP_\varphi)$ as operator systems (Definition~\ref{def:poisson-boundary}). Its existence and uniqueness are proved in Proposition~\ref{amenability of PB}. Under this identification, $M$ is a von Neumann subalgebra of $\CB_\varphi$.

For example, let $\Gamma$ be a countable discrete group, let $\mu\in\prob(\Gamma)$ be a generating probability measure, and let $(B,\nu_B)$ be its Poisson boundary. For $M=L(\Gamma)$, consider the hyperstate $\varphi\in B(\ell^2(\Gamma))_*$ defined by $\varphi(T)=\sum_{g\in\Gamma}\mu(g)\langle T\delta_g,\delta_g\rangle$. Izumi \cite{Izu04} proved that $M\subset\CB_\varphi$ is naturally isomorphic to the inclusion $L(\Gamma)\subset L^\infty(B,\nu_B)\rtimes\Gamma$.

The following result is the starting point of our approach to the weak Dixmier property. We deduce this property from the existence of a suitable map into $\CB_\varphi$.

\begin{letterthm}[Theorem~\ref{Thm ucp to PB}]\label{letterthm ucp to PB}
Let $M$ be a von Neumann algebra with separable predual, and let $\varphi$ be an nrsg hyperstate on $\BL$ with Poisson boundary $\CB_\varphi$. Assume that one of the following conditions holds:
\begin{itemize}
\item[(i)] $M$ is properly infinite, i.e., $M=M_{\mathrm{I}_\wx}\oplus M_{\mathrm{II}_\wx}\oplus M_{\mathrm{III}}$;
\item[(ii)] $M$ is finite with a normal faithful trace $\tau$, and $\varphi=\int_{\mathcal U(M)}\langle\,\cdot\,u\xi_\tau,u\xi_\tau\rangle\,\d\mu(u)$ for some $\mu\in\prob(\mathcal U(M))$.
\end{itemize}
Then every W$^*$-inclusion $M\subset\CA$ equipped with a normal faithful ucp map $\Phi:\CA\to\CB_\varphi$ satisfying $\Phi|_M=\id_M$ has the weak Dixmier property.
\end{letterthm}

As a direct corollary, we give a new proof of the weak Dixmier property for W$^*$-inclusions with a normal faithful conditional expectation, a problem posed by Popa \cite[Remark~1.5]{Pop99} and solved by Marrakchi \cite{Mar19}.

\begin{lettercor}[{\cite{Pop99,Mar19}}, Corollary~\ref{cor expectation}]\label{lettercor expectation}
Let $M\subset\CA$ be an inclusion of countably decomposable von Neumann algebras with a normal faithful conditional expectation. Then $M\subset\CA$ has the weak Dixmier property.
\end{lettercor}

We also prove a semifinite version of Theorem~\ref{letterthm ucp to PB}, replacing the normal faithful ucp map into $\CB_\varphi$ by a normal faithful semifinite $M$-bimodular positive map $T:\CA_+\to\widehat{\CB_\varphi}_+$, with the averaging conclusion for positive elements $x$ such that $T(x)$ is bounded (Theorem~\ref{Thm semifinite map}). In particular, we have a new proof of the weak Dixmier result of Isono \cite{Iso25} for inclusions $M\subset\CA$ of countably decomposable von Neumann algebras with a normal faithful semifinite operator-valued weight (Corollary~\ref{cor OVW}).

The proof of Theorem~\ref{letterthm ucp to PB} is rather simple. We explain the key case where $M$ is properly infinite with separable predual, and there is a normal faithful conditional expectation $E:\CA\to M$. Fix a normal faithful state $\phi$ on $M$ and take an nrsg hyperstate $\varphi$ extending $\phi$ as above. Using the coefficients $(y_k)$ of $\CQ_\varphi$ as in (\ref{eq:intro-inner-ucp-maps}), we define the normal $M$-inner ucp map
$$\CQ:\CA\to\CA,\qquad \CQ(a)=\sum_k y_kay_k^*\qquad(a\in\CA),$$
and let $\har(\CQ,\CA)=\{a\in\CA\mid \CQ(a)=a\}$. We first prove that $\har(\CQ,\CA)=M'\cap\CA$.

Define $\psi_\CA=\phi\circ E$. For $a\in\CA$ with $\CQ(a)=a$, by Kadison's inequality we have
$$\CQ(aa^*)-aa^*=\CQ(aa^*)-\CQ(a)\CQ(a)^*\geq0.$$
By (\ref{eq:intro-hyperstate-representations}), $\phi\circ \CQ|_M=\phi$. Since $E$ is $M$-bimodular, we also have
\begin{equation}\label{eq:intro-weight-preserving}
\psi_\CA\circ \CQ=\phi\circ E\circ \CQ
=\phi\circ \CQ\circ E=\phi\circ E=\psi_\CA.
\end{equation}
Hence $\psi_\CA(\CQ(aa^*)-aa^*)=0$, and by faithfulness we have $\CQ(aa^*)=aa^*$. Since $\CQ(a^*)=a^*$, we have
\begin{equation}\label{eq:intro-commutator}
\begin{aligned}
&\sum_k(y_ka-ay_k)(y_ka-ay_k)^*\\
={}&\sum_k\bigl(y_kaa^*y_k^*-y_kay_k^*a^*
-ay_ka^*y_k^*+ay_ky_k^*a^*\bigr)\\
={}&\CQ(aa^*)-\CQ(a)a^*-a\CQ(a^*)+aa^*=0.
\end{aligned}
\end{equation}
Hence $y_ka=ay_k$ for every $k$. Since $\C[(y_k)]$ is weakly dense in $M$, we have $a\in M'\cap\CA$, proving $\har(\CQ,\CA)=M'\cap\CA$.

For every $a\in\CA$, any weak cluster point of the Ces\`aro averages $\mathcal C_n(a)=(1/n)\sum_{j=0}^{n-1}\CQ^j(a)$ lies in $\har(\CQ,\CA)=M'\cap\CA$. If $y_k=\|y_k\|u_k$ with $u_k\in\mathcal U(M)$, these averages would lie in the weakly closed convex hull of the unitary orbit of $a$, which would prove the weak Dixmier property for $M\subset\CA$ directly. In general, the coefficients $(y_k)$ need not have this form. Nevertheless, the same argument shows that $M\subset\CA$ has the inner UCP-Dixmier property defined below.

\begin{definition}\label{def:intro-inner-ucp-dixmier}
Let $M\subset\CA$ be a W$^*$-inclusion. A ucp map $\CP:\CA\to\CA$ is called \textbf{$M$-inner} if
$$\CP(a)=\sum_k z_kaz_k^*,\qquad \sum_k z_kz_k^*=1,$$
for a sequence $(z_k)\subset M$. We say that $M\subset\CA$ has the \textbf{inner UCP-Dixmier property} if, for every $a\in\CA$,
$$\overline{\{\CP(a)\mid \CP\text{ is an }M\text{-inner ucp map}\}}^{w} \cap(M'\cap\CA)\ne\emptyset.$$
\end{definition}

Indeed, each Ces\`aro averaging map $\mathcal C_n$ is $M$-inner, so $M\subset\CA$ has the inner UCP-Dixmier property.

We prove that the two properties are equivalent when $M$ is properly infinite, which completes the proof in the case considered above. The proof uses proper infiniteness to pass from inner ucp maps to unitary averages, an idea similar to that used by Haagerup \cite[Section~2]{Haagerup1985Injectivity}.

\begin{letterthm}[Theorem~\ref{Thm inner UCP eq to Dixmier}]
\label{letterthm inner UCP eq to Dixmier}
Let $M\subset\CA$ be a W$^*$-inclusion such that $M$ is properly infinite. Then $M\subset\CA$ has the weak Dixmier property if and only if it has the inner UCP-Dixmier property.
\end{letterthm}

We briefly explain the proof of Theorem~\ref{letterthm inner UCP eq to Dixmier}. We only need to prove the if direction. By proper infiniteness, we can identify $M\subset\CA$ with $M\bar\otimes B(\ell^2(\N))\subset\CA\bar\otimes B(\ell^2(\N))$. We only need to prove the weak Dixmier property for the latter inclusion. Since $B(\ell^2(\N))$ is hyperfinite, every $x\in\CA\bar\otimes B(\ell^2(\N))$ can be averaged by unitaries of $1\otimes B(\ell^2(\N))$ to some $a\otimes1$ with $a\in\CA$. By the inner UCP-Dixmier property of $M\subset\CA$, we only need to prove that, for every $M$-inner ucp map $\CP:\CA\to\CA$,
$$\CP(a)\otimes1\in\overline{\operatorname{conv}}^{w} \{u(a\otimes1)u^*\mid u\in\mathcal U(M\bar\otimes B(\ell^2(\N)))\}.$$
Indeed, we can then average $x$ into $(M'\cap\CA)\bar\otimes\C1$, the relative commutant of the amplified inclusion.

Write $\CP(a)=\sum_{k\geq1}z_kaz_k^*$ with $(z_k)\subset M$ and $\sum_kz_kz_k^*=1$. For $n\geq1$, let $Z_n=(0_{n\times n},z_1\mathrm I_n,z_2\mathrm I_n,\ldots)$. Since $M$ is properly infinite, we can complete this block-row to a unitary $U_n\in M\bar\otimes B(\ell^2(\N))$ with first block-row $Z_n$, so that
$$U_n=\begin{pmatrix}Z_n\\\ast\end{pmatrix}\qquad\text{and}\qquad U_n(a\otimes1)U_n^*=
\begin{pmatrix}
\CP(a)\mathrm I_n&\ast\\
\ast&\ast
\end{pmatrix}.$$
Then we have $U_n(a\otimes1)U_n^*\to\CP(a)\otimes1$ weakly as $n\to\infty$, which proves the desired inclusion.

To summarize, the Poisson boundary theory provides inner ucp maps with the properties needed to prove the inner UCP-Dixmier property in Definition~\ref{def:intro-inner-ucp-dixmier}; when $M$ is properly infinite, we then apply Theorem~\ref{letterthm inner UCP eq to Dixmier} to prove the weak Dixmier property.

\subsection{Kadison's problem on MASA.}\label{sec:intro-kadison}

The same Poisson boundary tools also give a direct approach to Kadison's problem on maximal abelian subalgebras (MASA). In 1967, Kadison \cite{Kadison1967Problems} asked whether, for every inclusion of factors $M\subset P$ with $M'\cap P=\C1$, there is a maximal abelian subalgebra $A\subset M$ which is also maximal abelian in $P$, that is, $A'\cap P=A$. In 1981, Popa \cite{Popa1981Kadison} answered this question affirmatively when $M$ is of type $\mathrm{II}$ and there is a normal faithful conditional expectation from $P$ onto $M$, by constructing the desired maximal abelian subalgebra as the limit of an increasing sequence of finite-dimensional abelian subalgebras. Ge--Popa \cite[Corollary~5.1]{GP98} showed that the answer can be negative without such a conditional expectation, even when $M$ is a type $\mathrm{II}_1$ factor and $P$ is a type $\mathrm{II}_\infty$ factor, and asked whether the existence of a normal faithful conditional expectation suffices for factors of arbitrary type, see the discussion after \cite[Corollary~5.1]{GP98} and \cite[Problem~7.6]{Popa2021Ergodic}. Marrakchi \cite[Theorem~A]{Ma25} solved this problem when $M$ is not of type $\mathrm{III}_1$, or when $M$ is a type $\mathrm{III}_1$ factor with trivial bicentralizer. Marrakchi \cite[Theorem~A(2)]{Ma26} completed the solution for arbitrary type $\mathrm{III}_1$ factors, by reducing to the existence of a common irreducible hyperfinite subalgebra \cite[Theorem~B]{Ma26}, which relies on the ergodicity of the relative bicentralizer flow \cite[Theorem~C]{Ma26}.

We give a new proof of the solution to Kadison's problem, different from the proofs of Marrakchi \cite{Ma25,Ma26}. As in our proof of Theorem~\ref{letterthm relative bicentralizer}, the most difficult type $\mathrm{III}_1$ case is treated by the same argument as the other continuous properly infinite cases, without using any recent advances in modular and bicentralizer theory. Moreover, although both proofs use the same Poisson boundary framework, this proof is independent of our proof of the bicentralizer conjecture (Theorem~\ref{letterthm relative bicentralizer} and Corollaries~\ref{lettercor Connes bicentralizer} and~\ref{lettercor relative bicentralizer}).

\begin{letterthm}[{\cite{Popa1981Kadison,Ma25,Ma26}}]\label{letterthm kadison}
Let $M\subset P$ be an inclusion of von Neumann algebras with separable preduals. Assume that $M'\cap P=\mathcal Z(M)$ and that there is a normal faithful conditional expectation $E:P\to M$. Then there is an abelian von Neumann subalgebra $A\subset M$ such that $A'\cap P=A$.
\end{letterthm}

In the key case where $M$ is continuous and properly infinite, we construct finite-dimensional subalgebras $F_n\subset M$ such that $\|\theta|_{F_n'\cap P}\|\to0$ as $n\to\infty$ for every $\theta\in P_*$ with $\theta|_{M'\cap P}=0$. Once we have this, the rest follows the well-established route of Popa \cite{Popa1981Kadison}, also used in \cite{Ma25,Ma26}. We construct the $F_n$ directly using $\har(\CQ,P)=M'\cap P$ for the inner ucp map $\CQ$ associated with an nrsg hyperstate, as in Subsection~\ref{sec:intro-weak-dixmier-pb}, and a construction similar to that in the proof of Theorem~\ref{letterthm inner UCP eq to Dixmier}; while Marrakchi \cite[Theorem~B]{Ma26} first constructs a hyperfinite subalgebra $N\subset M$ with $N'\cap P=M'\cap P$ and then takes $(F_n)$ to be an increasing sequence of finite-dimensional subalgebras such that $N=\bigvee_nF_n$. See also \cite[Theorem~3.2 and Remark~3.5(2)]{Popa1981Kadison} for the construction of such a hyperfinite subalgebra $N\subset M$ when $M$ is semifinite. In Marrakchi's proof, the ergodicity of the relative bicentralizer flow \cite[Theorem~C]{Ma26} was used in the type $\mathrm{III}_1$ case. In our proof, the corresponding ergodicity input is the identity $\har(\CQ,P)=M'\cap P$.

With our Poisson boundary tools, we can in fact prove the following stronger result for properly infinite $M$, replacing the conditional expectation by a normal faithful ucp map into a Poisson boundary. It generalizes the hyperfinite-subalgebra existence results of Popa and Marrakchi in this case and also gives a new proof of them.

\begin{letterthm}[Theorem~\ref{thm:harmonic-implies-afd}]\label{letterthm harmonic afd}
Let $M\subset P$ be an inclusion of von Neumann algebras with separable preduals, and assume that $M$ is properly infinite. Let $\varphi$ be an nrsg $M$-hyperstate. If there is a normal faithful ucp map $\Phi:P\to\CB_\varphi$ with $\Phi|_M=\id_M$, then there is a hyperfinite von Neumann subalgebra $N\subset M$ such that $N'\cap P=M'\cap P$.
\end{letterthm}

Note that Theorem~\ref{letterthm harmonic afd} also generalizes Theorem~\ref{letterthm ucp to PB}(i) for inclusions with separable preduals. Indeed, since $N$ is hyperfinite, every element of $P$ can be averaged by unitaries in $\mathcal U(N)\subset\mathcal U(M)$ into $N'\cap P=M'\cap P$. Hence $M\subset P$ has the weak Dixmier property.

\subsection{Proof strategy for Theorem~\ref{letterthm relative bicentralizer}.}\label{sec:intro-proof-strategy}

Let $M\subset P$ be a W$^*$-inclusion with a normal faithful conditional expectation $E:P\to M$ as in Theorem~\ref{letterthm relative bicentralizer}.

We again focus on the key case where $M$ is properly infinite with separable predual. Fix a normal faithful state $\phi$ on $M$, and let $\rho=\phi\circ E$. In Subsections~\ref{sec:intro-weak-dixmier-pb} and~\ref{sec:intro-kadison} an arbitrary nrsg hyperstate suffices. Here we take the specific nrsg $(M,\phi)$-hyperstate $\varphi$ from Example~\ref{exists nrsg hyperstate}, with the symmetry $\varphi(Ja^*J)=\varphi(a)$ $(a\in\BL)$. This is the natural construction used to prove the existence of nrsg hyperstates. It also gives two further properties: the coefficient family $(y_k)$ of $\CQ_\varphi$ is permuted by $y_k\mapsto\sigma^\phi_{i/2}(y_k^*)$, and the operator norms of the $y_k$ satisfy $\sum_k\|y_k\|_\infty^2<\infty$, see the discussion before Example~\ref{exists nrsg hyperstate}. These are exactly the two properties needed for Proposition~\ref{letterprop:intro-vanishing} below.

Let $\CA=c_\rho(P)$, and define $\CQ(x)=\sum_k y_kxy_k^*$ ($x\in\CA$) using the coefficients $(y_k)$ of $\CQ_\varphi$. As in Subsection~\ref{sec:intro-weak-dixmier-pb}, to prove the weak Dixmier property for $M\subset\CA$, we only need to prove that $\har(\CQ,\CA)=M'\cap\CA$.

Let $T:\CA_+\to\widehat P_+$ be the canonical normal faithful semifinite operator-valued weight of the continuous core. Define $\psi_\CA=\rho\circ T$, the dual weight of $\rho$, see Subsection~\ref{sec:core-weights}.

By (\ref{eq:intro-commutator}), to prove $\har(\CQ,\CA)=M'\cap\CA$, we only need to prove that $\CQ(aa^*)=aa^*$ for every $a\in\har(\CQ,\CA)$. Fix such an $a$, and let $d:=\CQ(aa^*)-aa^*\geq0$.

As in (\ref{eq:intro-weight-preserving}), we have $\psi_\CA\circ \CQ=\psi_\CA$. Intuitively, we would have $\psi_\CA(d)=\psi_\CA(\CQ(aa^*))-\psi_\CA(aa^*)=0$, but this calculation is justified only when $\psi_\CA(aa^*)<\infty$. Instead, Lemma~\ref{lem:positive-localization} shows that, for every $0\ne d'\in\CA_+$, there exists $x'\in\CA_+$ such that $0\ne x'\leq d'$ and $0<\psi_\CA(x')<\infty$. Hence we only need to prove that every $x\in\CA_+$ satisfying $x\leq d$ and $\psi_\CA(x)<\infty$ must satisfy $x=0$. Since $x\leq d=\CQ(aa^*)-aa^*$ and $\sum_{k=0}^{n-1}\CQ^k(d)=\CQ^n(aa^*)-aa^*$, we have $\sum_{k=0}^{n-1}\CQ^k(x)\leq C1$ for every $n\geq1$, with $C=\|a\|^2$.

All we need to prove is the following proposition.

\begin{proposition}[Lemmas~\ref{lem:partial-sum-estimate} and~\ref{lem:cutoff-estimates}]\label{letterprop:intro-vanishing}
Let $\CA$, $\psi_\CA$ and $\CQ$ be as above, and let $x\in\CA_+$ satisfy $\psi_\CA(x)<\infty$. Assume that there is a constant $0\leq C<\infty$ such that
$$\sum_{k=0}^{n-1}\CQ^k(x)\leq C1\qquad(n\geq1).$$
Then $x=0$.
\end{proposition}

The estimates in the proof of Proposition~\ref{letterprop:intro-vanishing} use a noncommutative analogue of the classical cutoff argument for recurrence of symmetric Markov operators \cite[Section~1.6, Theorem~1.6.3]{FOT11}. We explain this argument through the commutative model in Proposition~\ref{prop:commutative-recurrence-model}. The symmetry of the nrsg hyperstate $\varphi$ is essential here. Without it, the vanishing conclusion already fails in the commutative model, see Remark~\ref{rem:commutative-symmetry}.

\subsection*{Acknowledgments.}
The author thanks Cyril Houdayer, Amine Marrakchi, Soham Chakraborty and Stuart White for many helpful discussions and encouragement. Proposition~\ref{exists Pphi} benefited from a short note written by Soham Chakraborty and Se-Jin (Sam) Kim. Theorem~\ref{Thm inner UCP eq to Dixmier} benefited from a discussion with Amine Marrakchi during the author's visit to ENS Lyon in October 2025. Part of this work (Sections~\ref{sec:poisson-boundaries} and~\ref{sec:isono-theorem}) was completed while the author was a PhD student at ENS Paris.

\subsection*{Independent and concurrent work.}
Our solution to the relative bicentralizer conjecture (Theorem~\ref{letterthm relative bicentralizer} and Corollary~\ref{lettercor relative bicentralizer}) was obtained in parallel with and independently of the concurrent work \cite{OAI26}, released on 6 October 2026. The two proofs use completely different approaches. OpenAI's proof builds on the work of Houdayer--Marrakchi \cite{HM26}, while ours uses the Poisson boundary techniques developed here. See the subsection on AI-assisted research methodology below for the timeline.

\subsection*{AI-assisted research methodology.}
It is well known that AI can help solve a given problem by identifying techniques whose relevance had previously gone unnoticed. Our solution to the relative bicentralizer conjecture (Theorem~\ref{letterthm relative bicentralizer}), however, illustrates the reverse process: starting from a given body of techniques developed by the human author, AI helped identify a solvable problem whose connection with these techniques the author had not previously noticed.

The underlying Poisson boundary theory in Section~\ref{sec:poisson-boundaries} was developed entirely by the human author between February 2025 and June 2025, benefiting from helpful discussions with Cyril Houdayer and Soham Chakraborty.

The new proof for Popa's problem on the weak Dixmier property (Corollary~\ref{lettercor expectation}) was also purely human work. During the author's visit to ENS Lyon in October 2025, helpful discussions with Amine Marrakchi led to the observation that the framework in Section~\ref{sec:poisson-boundaries} could be used to simplify the proof of this result. Between October 2025 and May 2026, the author further generalized this result to inclusions with maps into Poisson boundaries (Theorem~\ref{letterthm ucp to PB}) and to the semifinite setting of Isono's result (Corollary~\ref{cor OVW}), and prepared a manuscript containing the material of the present Sections~\ref{sec:poisson-boundaries} and~\ref{sec:isono-theorem}. A \$20 version of Gemini 3.1 Pro served only as a writing aid during April 2026 (for parts that remain in the present text). During subsequent checking of this manuscript, GPT-6 Astra helped identify and repair minor gaps.

In September 2026, the author began using GPT-6 Astra to investigate whether every type $\mathrm{II}_1$ factor admits a type $\mathrm{III}_1$ Poisson boundary, a question that had interested the author since the early stages of his doctoral studies. This investigation did not resolve the original question. The author then asked the model to consider the weaker conclusion that every type $\mathrm{II}_1$ factor admits an irreducible embedding into the hyperfinite type $\mathrm{III}_1$ factor. This led to a construction using Poisson boundaries and crossed products. To extend this construction to arbitrary continuous factors, the author gave the model his earlier Poisson boundary manuscript on general von Neumann algebras. The model then helped develop the corresponding embedding arguments. To keep this paper concise and focused, these irreducible embedding results are not included in the present version. We refer interested readers to a preliminary version~\cite{hexagon:2610.00134v2} of this paper.

Having recognized Poisson boundaries as a very powerful tool, the author then asked GPT-6 Astra on 26 September 2026 which other open problems might be solved using the existing Poisson boundary theory. The model suggested three problems that were still open at the time, before OpenAI's public release of 722 math manuscripts on 6 October 2026~\cite{OAI26Announcement}:
\begin{enumerate}
\item the relative bicentralizer conjecture;
\item the free group factor problem;
\item the vanishing cohomology conjecture of Kadison--Ringrose.
\end{enumerate}
The author chose the relative bicentralizer conjecture for further study. The main tools needed for the proof, including the Poisson boundary construction, symmetric nrsg hyperstates, and Theorem~\ref{letterthm inner UCP eq to Dixmier}, had already been developed in the earlier Poisson boundary manuscript. GPT-6 Astra helped connect these tools to the conjecture and work out their application. With this assistance, the author obtained a preliminary proof on 27 September 2026.

The human author identified Kadison's problem (Theorem~\ref{letterthm kadison}) as another application of the Poisson boundary tools: Following the successful treatment of Popa's problem on the weak Dixmier property and the relative bicentralizer conjecture, the author recognized the parallels between these two problems: both are formulated for general von Neumann algebras, with classical type $\mathrm{II}_1$ cases but extremely difficult type $\mathrm{III}_1$ cases. The author saw the same contrast in Kadison's problem on MASA and therefore proposed applying these tools to it. The proof uses the generating inner ucp maps constructed in the author's earlier Poisson boundary manuscript and the unitary dilation argument used in Theorem~\ref{letterthm inner UCP eq to Dixmier}. With the assistance of GPT-6 Astra in applying these tools to Kadison's problem, the author obtained a substantially simplified proof on 9 October 2026.

Lean was also used to check a formalization of a short note consisting of the preliminaries in Subsection~\ref{sec:core-weights}, Example~\ref{exists nrsg hyperstate}, Theorem~\ref{Thm inner UCP eq to Dixmier}, and the whole of Section~\ref{sec relative bicentralizer}. These components form the main line of the proof of Theorem~\ref{letterthm relative bicentralizer}. Results from the published literature were taken as explicit assumptions.

The author is responsible for the mathematical claims, citations, and final presentation.

\section{Preliminaries}\label{sec:preliminaries}
\subsection{The standard form of von Neumann algebras.}

Let $M$ be a countably decomposable von Neumann algebra. Following \cite{Ha75}, the standard form of $M$ is a quadruple $(M,\mathcal{H},J,\mathcal{H}_+)$ such that $\mathcal{H}$ is a Hilbert space with $M\subset B(\mathcal{H})$; $J:\mathcal{H}\to \mathcal{H}$ is an anti-unitary with $J^2=1$; $\mathcal{H}_+\subset \mathcal{H}$ is a self-dual cone, i.e., its dual cone $\mathcal{H}_+^{\mathrm{o}}:=\{\xi\in \mathcal{H}\mid\langle\xi,\eta\rangle\geq0,\ \forall\eta\in \mathcal{H}_+\}$ satisfies $\mathcal{H}_+^{\mathrm{o}}=\mathcal{H}_+$; and the following conditions hold:
\begin{itemize}
    \item[(1)] $JMJ=M'\cap B(\mathcal{H})$;
    \item[(2)] $JcJ=c^*$ for any $c\in\mathcal Z(M)$, the center of $M$;
    \item[(3)] $J\xi=\xi$ for any $\xi\in \mathcal{H}_+$;
    \item[(4)] $a(JaJ)\mathcal{H}_+\subset \mathcal{H}_+$ for any $a\in M$.
\end{itemize}
The standard form is unique up to unitary equivalence. Assume that $(M,\tilde{\mathcal{H}},\tilde J,\tilde{\mathcal{H}}_+)$ is another quadruple satisfying the above conditions. Then there exists a unique unitary $u:\mathcal{H}\to\tilde{\mathcal{H}}$ such that
\begin{itemize}
    \item[(a)] $\pi_{\tilde{\mathcal{H}}}(a)=u\pi_{\mathcal{H}}(a)u^*$ ($a\in M$), where $\pi_{\mathcal{H}}:M\to B(\mathcal{H})$ and $\pi_{\tilde{\mathcal{H}}}:M\to B(\tilde{\mathcal{H}})$ are the representations of $M$ on $\mathcal{H}$ and $\tilde{\mathcal{H}}$, respectively;
    \item[(b)] $\tilde J=uJu^*$;
    \item[(c)] $\tilde{\mathcal{H}}_+=u\mathcal{H}_+$.
\end{itemize}
Therefore, the standard form of $M$ depends only on $M$ itself, and we denote it by $(M,L^2(M),J_M,L^2(M)_+)$. In particular, for any normal faithful state $\phi$ on $M$, we can take
$$(M,L^2(M),J_M,L^2(M)_+)=(M,L^2(M,\phi),J_\phi,L^2(M,\phi)_+),$$
where $J_\phi$ is the modular conjugation operator determined by the polar decomposition $S_\phi=J_\phi\Delta_\phi^{1/2}$ of the Tomita operator, and $L^2(M,\phi)_+$ is the closure of $\{a(J_\phi aJ_\phi)\xi_\phi\mid a\in M\}$ \cite[Lemma~2.9]{Ha75}.

By \cite[Lemma~2.10]{Ha75}, any normal positive functional $\omega\in M_*^+$ admits a unique representing vector $\xi_\omega\in L^2(M)_+$, i.e., $\omega=\langle\,\cdot\,\xi_\omega,\xi_\omega\rangle$. Moreover, $\xi\mapsto\omega_\xi=\langle\,\cdot\,\xi,\xi\rangle$ is a homeomorphism between $L^2(M)_+$ and $M_*^+$.

\subsection{Continuous cores and dual weights.}\label{sec:core-weights}

For background on modular theory and continuous cores, we refer to \cite{TakII}.

Let $P$ be a countably decomposable von Neumann algebra, and let $\rho$ be a normal faithful state on $P$. Let $\xi_\rho\in L^2(P)_+$ be the vector representing $\rho$. The Tomita operator $S_\rho$ is the closure on $L^2(P)$ of the densely defined antilinear map $x\xi_\rho\mapsto x^*\xi_\rho$ ($x\in P$). Define the modular operator by $\Delta_\rho=S_\rho^*S_\rho$. The modular automorphism group $\sigma^\rho:\R\to\operatorname{Aut}(P)$ is given by
$$\sigma_t^\rho(x)=\Delta_\rho^{it}x\Delta_\rho^{-it} \qquad(x\in P,\ t\in\R).$$
We say that $x\in P$ is \textbf{entire analytic} for $\sigma^\rho$ if $t\mapsto\sigma_t^\rho(x)$ extends to a map $z\mapsto\sigma_z^\rho(x)$ from $\C$ to $P$ such that $z\mapsto\omega(\sigma_z^\rho(x))$ is entire for every $\omega\in P_*$.

Define the continuous core of $P$ by
$$\CA=c_\rho(P)=P\rtimes_{\sigma^\rho}\R.$$
We identify $P$ with its canonical image in $\CA$ and denote by $(\lambda_t)_{t\in\R}$ the canonical implementing unitaries, satisfying $\lambda_t x\lambda_t^*=\sigma_t^\rho(x)$ ($x\in P$, $t\in\R$). The dual action $\widehat\sigma:\R\curvearrowright\CA$ is determined by
$$\widehat\sigma_r(x)=x,\qquad \widehat\sigma_r(\lambda_t)=e^{-irt}\lambda_t \qquad(x\in P,\ r,t\in\R).$$
The canonical operator-valued weight $T:\CA_+\to\widehat P_+$ is defined by
$$T(x)=\int_\R\widehat\sigma_r(x)\,\d r\qquad(x\in\CA_+).$$
The dual weight of $\rho$ is $\psi_\CA=\rho\circ T$. Write $\tau_\CA$ for the canonical trace on $\CA$, and let
$$h=\frac{\d\psi_\CA}{\d\tau_\CA},\qquad H=\log h,\qquad \lambda_t=h^{it}=e^{itH}.$$
Products of measurable operators are understood as closed products.

We next prove two lemmas that will be used in Section~\ref{sec relative bicentralizer}.

\begin{lemma}\label{prop:core-integration}
For every positive integrable function $g$ on $\R$,
\begin{equation}\label{eq:core-normalization}
\psi_\CA(g(H))=\rho(T(g(H)))=\int_\R g(s)\,\d s.
\end{equation}
In particular, we have
\begin{equation}\label{eq:tail}
\tau_\CA(\mathds{1}_{(s,\infty)}(h))=\frac1s\qquad(s>0).
\end{equation}
\end{lemma}
\begin{proof}
For $r,t\in\R$, by functional calculus we have
$$e^{it\widehat\sigma_r(H)} =\widehat\sigma_r(e^{itH}) =\widehat\sigma_r(\lambda_t) =e^{-irt}e^{itH} =e^{it(H-r)}.$$
By uniqueness of the self-adjoint generator of a strongly continuous unitary group,
\begin{equation}\label{eq:dual-action-H}
\widehat\sigma_r(H)=H-r.
\end{equation}
Let $\mathcal E_H$ be the spectral measure of $H$. For a positive integrable function $g$, by (\ref{eq:dual-action-H}) we have
\begin{equation}\label{eq:core-operator-valued-normalization}
T(g(H)) =\int_\R g(H-r)\,\d r =\int_\R\left(\int_\R g(s-r)\,\d r\right)\,\d\mathcal E_H(s) =\left(\int_\R g(s)\,\d s\right)1.
\end{equation}
Since $\psi_\CA=\rho\circ T$ and $\rho(1)=1$, by (\ref{eq:core-operator-valued-normalization}) we have (\ref{eq:core-normalization}).

For $s>0$, by the definition of $h$ and (\ref{eq:core-normalization}), we have
$$
\tau_\CA(\mathds{1}_{(s,\infty)}(h)) =\psi_\CA\bigl(h^{-1}\mathds{1}_{(s,\infty)}(h)\bigr) =\psi_\CA\bigl(e^{-H}\mathds{1}_{(\log s,\infty)}(H)\bigr) =\int_{\log s}^{\infty}e^{-t}\,\d t=\frac1s.
$$
See also \cite[Theorem~1.2(2)]{Ha79Lp}. In particular, $h$ is $\tau_\CA$-measurable.
\end{proof}

\begin{lemma}\label{lem:positive-localization}
For every $0\ne d\in\CA_+$, there exists $x\in\CA_+$ such that
$$0\ne x\leq d,\qquad\psi_\CA(x)<\infty.$$
\end{lemma}
\begin{proof}
Take $\delta>0$ such that $\mathds{1}_{[\delta,\infty)}(d)\ne0$. By semifiniteness of $\tau_\CA$, there is a nonzero projection
$$p\leq\mathds{1}_{[\delta,\infty)}(d),\qquad 0<\tau_\CA(p)<\infty.$$
By (\ref{eq:tail}), take $s>0$ such that $\tau_\CA(\mathds{1}_{(s,\infty)}(h))<\tau_\CA(p)$, and let $e=\mathds{1}_{[0,s]}(h)$ and $q=p\wedge e$. Let $(1-e)p=v|(1-e)p|$ be the polar decomposition. The kernel of $(1-e)p$ inside the range of $p$ is the range of $q=p\wedge e$, and the range of $(1-e)p$ is contained in the range of $1-e$. Hence
$$v^*v=p-q,\qquad vv^*\leq1-e.$$
Hence $p-q\sim vv^*\leq1-e$, so $p-q\precsim1-e$ and
$$\tau_\CA(q)\geq\tau_\CA(p)-\tau_\CA(1-e)>0.$$
Since $q\leq e$, we have
$$\psi_\CA(q)=\tau_\CA(qhq)=\tau_\CA(q(ehe)q)\leq s\tau_\CA(q)<\infty.$$
Hence $x=\delta q$ satisfies
$$0\ne x=\delta q\leq\delta p\leq d,\qquad \psi_\CA(x)=\delta\psi_\CA(q)<\infty.$$
\end{proof}

\subsection{Centralizers and bicentralizers.}\label{sec:centralizers}

Let $M$ be a countably decomposable von Neumann algebra, and let $\phi$ be a normal faithful state on $M$. The \textbf{centralizer} of $\phi$ is
$$M_\phi=\{a\in M\mid \phi(ax)=\phi(xa),\ \forall x\in M\} =\{a\in M\mid \sigma_t^\phi(a)=a,\ \forall t\in\R\}.$$
The \textbf{asymptotic centralizer} of $\phi$ is
$$\operatorname{AC}(M,\phi) =\{(a_n)\in\ell^\infty(\N,M)\mid \|a_n\phi-\phi a_n\|_{M_*}\to0\},$$
where $(a\phi)(x)=\phi(xa)$ and $(\phi a)(x)=\phi(ax)$ for $a,x\in M$. For a W$^*$-inclusion $M\subset P$, the \textbf{relative bicentralizer} is defined by
$$\mathrm B(M\subset P,\phi) =\{x\in P\mid [x,a_n]\xrightarrow{\mathrm{s}}0,\quad \forall(a_n)\in\operatorname{AC}(M,\phi)\},$$
where $\mathrm{s}$ denotes strong operator convergence, see \cite[Definition~3.1]{AHHM2020Bicentralizer}. For $P=M$, we have the \textbf{bicentralizer} of $\phi$
$$\mathrm B(M,\phi)=\mathrm B(M\subset M,\phi).$$

Marrakchi's \textbf{relative bicentralizer conjecture} \cite[Conjecture~C]{Ma25} asserts that, for every inclusion $M\subset P$ of countably decomposable von Neumann algebras with a normal faithful conditional expectation $E:P\to M$ and every normal faithful state $\phi$ on $M$,
\begin{equation}\label{eq:target}
\mathrm B(M\subset c(P),\phi)
=\{\phi^{it}\mid t\in\R\}''\vee(M'\cap c(P)).
\end{equation}
Here we use the compatible realizations $c_\phi(M)\subset c_{\phi\circ E}(P)=c(P)$, and $\phi^{it}$ denotes the canonical implementing unitary in $c_\phi(M)$.

\section{Poisson boundaries of general von Neumann algebras}\label{sec:poisson-boundaries}
In this section, we develop the Poisson boundary theory for general von Neumann algebras, following the framework of Das--Peterson \cite{DP20} for tracial von Neumann algebras. The proofs of ergodicity, double ergodicity, and rigidity are in the same spirit as those of \cite{DP20}. However, the coefficient families $(z_k)$ and $(y_k)$ for $\CP_\varphi$ and $\CQ_\varphi$ need not coincide in the nontracial case, whereas we can take $y_k=z_k$ if $\varphi|_M$ is tracial (Remark~\ref{rem:tracial-hyperstate-coefficients}). Hence, in the general setting, we cannot recover all the results proved in \cite{DP20} for finite von Neumann algebras.

Throughout this section, we fix a von Neumann algebra $M$ with separable predual.

\begin{definition}
Let $M \subset \mathcal{A}$ be a W$^*$-inclusion. A normal state $\varphi$ on $\mathcal{A}$ is called an \textbf{$M$-hyperstate} if its restriction $\varphi|_M$ to $M$ is faithful. Given a fixed normal faithful state $\phi$ on $M$, we say that $\varphi$ is an \textbf{$(M,\phi)$-hyperstate} if $\varphi|_M = \phi$.
\end{definition}

The framework of Das--Peterson \cite{DP20} concerns $(M,\tau)$-hyperstates for tracial von Neumann algebras $(M,\tau)$, whereas we work with general $M$-hyperstates without fixing a normal faithful state $\phi$ on $M$ in advance.

The following proposition associates with each normal $M$-hyperstate on $\BL$ an $M$-bimodular ucp map, which we use to define the Poisson boundary.

\begin{proposition}\label{exists Pphi}
Let $\varphi \in B(L^2(M))_*$ be a normal $M$-hyperstate. Set $\phi = \varphi|_M$, and let $\xi_\phi \in L^2(M)_+$ be the canonical cyclic and separating vector implementing $\phi$. Then there exists a unique normal $M$-bimodular ucp map $\mathcal{P}_\varphi : B(L^2(M)) \to B(L^2(M))$ such that for all $T \in B(L^2(M))$,
$$ \varphi(T) = \langle \mathcal{P}_\varphi(T)\xi_\phi, \xi_\phi \rangle. $$
Moreover, there exists a (not necessarily unique) sequence $(z_k)_{k=1}^\infty \subset M$ such that
$$\mathcal P_\varphi(T)=\sum_{k=1}^\infty(Jz_k^*J)T(Jz_kJ) \qquad(T\in B(L^2(M))),$$
where $J$ is the modular conjugation operator on $L^2(M)$.
\end{proposition}

\begin{proof}
We first prove the existence and construct the explicit form of $\mathcal{P}_\varphi$.

Since $\varphi$ is normal, there exists a sequence $(\xi_k)\subset L^2(M)$ such that
$$\varphi(T)=\sum_k\langle T\xi_k,\xi_k\rangle\qquad(T\in B(L^2(M))).$$

For any $x \in M$, since $\varphi|_M = \phi$, we have
$$ \|x\xi_\phi\|^2 = \phi(x^*x) = \varphi(x^*x) = \sum_k \|x\xi_k\|^2. $$
Hence
\begin{equation}\label{eq:vector-component-bound}
\|x\xi_n\|^2\leq\|x\xi_\phi\|^2
\qquad(x\in M,\ n\geq1).
\end{equation}

Define the operator $R_n : M\xi_\phi \to L^2(M)$ by $R_n(x\xi_\phi) = x\xi_n$. By (\ref{eq:vector-component-bound}), $\|R_n(x\xi_\phi)\| \le \|x\xi_\phi\|$, so $R_n$ extends to a contraction on $L^2(M)$. Moreover, for any $x, y \in M$, we have $R_n(yx\xi_\phi) = yx\xi_n = yR_n(x\xi_\phi)$. Hence $R_n\in M'$. By standard Tomita--Takesaki theory, $M' = J M J$, so there exists $z_n \in M$ such that $R_n = J z_n J$.

Now, for any $x, y \in M$ and $T \in B(L^2(M))$, we calculate
\begin{equation}\label{eq:hyperstate-bimodule-calculation}
\begin{aligned}
\varphi(y^*Tx) &= \sum_k \langle Tx\xi_k,y\xi_k\rangle
= \sum_k \langle T R_k(x\xi_\phi), R_k(y\xi_\phi) \rangle \\
&= \sum_k \langle T J z_k J x\xi_\phi, J z_k J y\xi_\phi \rangle \\
&= \left\langle \left( \sum_k J z_k^* J \, T \, J z_k J \right) x\xi_\phi, y\xi_\phi \right\rangle.
\end{aligned}
\end{equation}

Let $\mathcal{P}_\varphi(T) = \sum_k (J z_k^* J) T ( J z_k J)$. By (\ref{eq:hyperstate-bimodule-calculation}), $\varphi(y^*Tx)=\langle\mathcal P_\varphi(T)x\xi_\phi,y\xi_\phi\rangle$. Taking $x=y=1$, we immediately obtain $\varphi(T) = \langle \mathcal{P}_\varphi(T) \xi_\phi, \xi_\phi \rangle$.

Since $M\xi_\phi$ is dense in $L^2(M)$, the identity
$$\langle\mathcal P_\varphi(T)x\xi_\phi,y\xi_\phi\rangle=\varphi(y^*Tx)$$
uniquely determines $\mathcal P_\varphi(T)$ for every $T$. Hence $\mathcal P_\varphi$ is unique.

Taking $T=1$ in (\ref{eq:hyperstate-bimodule-calculation}) and using $\varphi|_M=\phi$, we have $\sum_k z_k^*z_k=1$. Since $Jz_kJ\in M'$, $\mathcal P_\varphi$ is a normal $M$-bimodular ucp map.
\end{proof}

\begin{definition}\label{def:poisson-boundary}
Let $\varphi$ be a normal $M$-hyperstate on $\BL$. Let $\har(\CP_\varphi) = \{ T \in \BL \mid \CP_\varphi(T) = T \}$ denote the space of \textbf{$\varphi$-harmonic operators}. The \textbf{Poisson boundary} $(\CB_\varphi,\zeta)$ is defined as a von Neumann algebra that is isomorphic to $\har(\CP_\varphi)$ as an operator system, together with the canonical hyperstate $\zeta:=\varphi|_{\har(\CP_\varphi)}$.
\end{definition}

The next proposition shows that the Poisson boundary is well defined and amenable.

\begin{proposition}\label{amenability of PB}
Let $\varphi$ be a normal $M$-hyperstate on $\BL$. The Poisson boundary $\CB_\varphi$ always exists, is unique, and is amenable.
\end{proposition}

\begin{proof}
Let $\omega \in \beta\N \setminus \N$ be a free ultrafilter, and let $\mathcal C_n=(1/n)\sum_{k=0}^{n-1}\CP_\varphi^k$. We define a conditional expectation $E : \BL \to \har(\CP_\varphi)$ by
$$E(T) = \operatorname*{wo-lim}_{n \to \omega} \mathcal C_n(T) \quad(T\in\BL).$$

By the Choi--Effros theorem \cite{CE77}, equipping the operator system $\har(\CP_\varphi)$ with the \textbf{Choi--Effros product} defined by
$$x \circ y = E(xy) \quad(x,y\in\har(\CP_\varphi))$$
makes $(\har(\CP_\varphi),\circ)$ a C$^*$-algebra with its original norm. The harmonic operator system is ultraweakly closed, hence this C$^*$-algebra is a von Neumann algebra, with its canonical ultraweak topology agreeing with the inherited one, see also \cite{Izu02}. The inclusion $M\subset\har(\CP_\varphi)$ is multiplicative for this product, since $E$ is $M$-bimodular.

The uniqueness of $\CB_\varphi$ follows from the fact that any two von Neumann algebras that are completely isometrically isomorphic as operator systems are necessarily *-isomorphic.

Finally, since $\CB_\varphi\cong \har(\CP_\varphi)$ as injective operator systems, $\CB_\varphi$ is amenable as a von Neumann algebra.
\end{proof}

Under the following symmetry condition, we can also represent $\varphi$ by an $M$-inner ucp map.

\begin{definition}
We say that a normal $M$-hyperstate $\varphi \in B(L^2(M))_*$ is \textbf{regular} if $\varphi(Jz^*J) = \varphi(z)$ for all $z \in M$.
\end{definition}

For regular hyperstates, Proposition~\ref{exists Pphi} has the following counterpart, with $M$ and $JMJ$ interchanged.

\begin{proposition}\label{prop:regular_P_phi_o}
Let $\varphi\in\BL_*$ be a normal regular $M$-hyperstate, and let $\phi=\varphi|_M$. Then there exists a sequence $(y_k)_{k=1}^\infty\subset M$ with $\sum_k y_ky_k^*=1$ such that
$$\CQ_\varphi(T)=\sum_k y_kTy_k^*\qquad(T\in\BL)$$
defines a normal $JMJ$-bimodular ucp map satisfying
$$\varphi(T)=\langle\CQ_\varphi(T)\xi_\phi,\xi_\phi\rangle \qquad(T\in\BL).$$
\end{proposition}

\begin{proof}
Let $\phi^o=\varphi|_{JMJ}$. Since $\varphi$ is regular and $J\xi_\phi=\xi_\phi$, we have
$$\langle Jz^*J\xi_\phi,\xi_\phi\rangle=\phi(z)=\phi^o(Jz^*J) \qquad(z\in M).$$
Hence $\xi_\phi$ is the canonical cyclic and separating vector implementing $\phi^o$. By applying Proposition~\ref{exists Pphi} to $JMJ$, we have the required map and coefficients.
\end{proof}

\begin{remark}\label{rem:tracial-hyperstate-coefficients}
In general, $(y_k)$ need not coincide with $(z_k)$ from Proposition~\ref{exists Pphi}. If $\phi$ is tracial, we can take $y_k=z_k$. This is the main case studied in the Das--Peterson framework \cite{DP20}.
\end{remark}

\begin{definition}\label{def strong generation}
Let $\varphi$ be a normal regular $M$-hyperstate. We say that $\varphi$ is \textbf{strongly generating} if there exist sequences $(z_k),(y_k)\subset M$ with
$$\overline{\C[(z_k)]}^{\,w}=\overline{\C[(y_k)]}^{\,w}=M$$
such that $\CP_\varphi$ and $\CQ_\varphi$ have the standard forms
$$\CP_\varphi(T)=\sum_k(Jz_k^*J)T(Jz_kJ),\qquad \CQ_\varphi(T)=\sum_k y_kTy_k^*\qquad(T\in\BL).$$
Here $\C[(z_k)]$ and $\C[(y_k)]$ are the unital subalgebras (rather than $*$-algebras) generated by $(z_k)$ and $(y_k)$, respectively. We abbreviate \textbf{normal regular strongly generating} hyperstate as \textbf{nrsg} hyperstate. In the nontracial setting we assume both density conditions. We do not assert that one follows from the other.
\end{definition}

We now prove the existence of nrsg hyperstates. Although regularity requires $\varphi(JT^*J)=\varphi(T)$ only for $T\in M$, it is easier to impose this \textbf{symmetry} for all $T\in\BL$. Indeed, for any normal faithful state $\varphi_0$ on $\BL$, the state $\varphi(T)=\frac{1}{2}\bigl(\varphi_0(T)+\varphi_0(JT^*J)\bigr)$ $(T\in\BL)$ is a normal faithful symmetric $M$-hyperstate, and hence regular.

Fix a normal faithful state $\phi$ on $M$, and let $\xi=\xi_\phi\in L^2(M)_+$. To construct a strongly generating symmetric $(M,\phi)$-hyperstate, we want the coefficient families $(z_k)$ and $(y_k)$ in Definition~\ref{def strong generation} to generate $M$ and satisfy
$$\sum_k\langle Ty_k^*\xi,y_k^*\xi\rangle =\sum_k\langle TJy_k^*\xi,Jy_k^*\xi\rangle\qquad(T\in\BL).$$
This equality holds if $(Jy_k^*\xi)$ is a permutation of $(y_k^*\xi)$. For analytic coefficients, we have
$$Jy_k^*\xi=\sigma_{-i/2}^\phi(y_k)\xi =\bigl(\sigma_{i/2}^\phi(y_k^*)\bigr)^*\xi.$$
Hence we naturally ask that $y_k\mapsto\sigma_{i/2}^\phi(y_k^*)$ permute the coefficient family. We therefore take an analytic family $(a_n)$ with weakly dense linear span and let $b_n=\sigma_{-i/2}^\phi(a_n^*)$. We pair the coefficients $\sqrt{c_n}a_n^*$ and $\sqrt{c_n}b_n^*$, using the same weight $c_n>0$.

It is not easy to arrange $\phi\circ \CQ_0=\phi$ directly for $\CQ_0(z)=\sum_n c_n(a_n^*za_n+b_n^*zb_n)$ $(z\in M)$ with these paired coefficients. However, we can take the weights $c_n$ small enough that $\psi:=\phi\circ \CQ_0\leq\phi$ and $\sum_n c_n(\|a_n\|_\infty^2+\|b_n\|_\infty^2)<\infty$. By \cite[Lemma~2.10(1)]{Ha75}, let $\eta\in L^2(M)_+$ represent $\phi-\psi$. Since $J\eta=\eta$, adding the positive functional $\langle\,\cdot\,\eta,\eta\rangle$ to the paired sum preserves this symmetry and makes its restriction to $M$ equal to $\phi$.

The following example gives this construction. The properties in (ii) follow naturally from it and will be used in Section~\ref{sec relative bicentralizer}.

\begin{example}\label{exists nrsg hyperstate}
Let $\phi$ be a normal faithful state on $M$, and write $\xi=\xi_\phi\in L^2(M)_+$. There exists a faithful nrsg $(M,\phi)$-hyperstate $\varphi$ on $\BL$ with inner ucp coefficients $(y_k)\subset M$ for $\CQ_\varphi$ having the following properties:
\begin{itemize}
\item[(i)] The hyperstate is given by $\varphi(T)=\sum_k\langle Ty_k^*\xi,y_k^*\xi\rangle$ and satisfies
\begin{equation}\label{eq:nrsg-J-symmetry}
\varphi(JT^*J)=\varphi(T)\qquad(T\in\BL).
\end{equation}
In particular, the $M$-inner ucp map $\CQ(a)=\sum_k y_kay_k^*$ $(a\in M)$ satisfies $\phi\circ \CQ=\phi$.
\item[(ii)] Each $y_k^*$ admits a bounded modular continuation to $i/2$, and the map
\begin{equation}\label{eq symmetric coefficient involution}
y_k\longmapsto\sigma_{i/2}^\phi(y_k^*)
\end{equation}
permutes the family $(y_k)$. Moreover, $\sum_k\|y_k\|_\infty^2<\infty$.
\item[(iii)] If $\phi$ is tracial, we may further take $y_k=\lambda_k u_k$ with $u_k\in\mathcal U(M)$ and $\lambda_k>0$ such that $\sum_k\lambda_k^2=1$. In particular, $\CQ(a)=\sum_k y_kay_k^*$ is a unitary averaging map.
\end{itemize}
\end{example}

\begin{proof}
Take elements $a_n\in M$, entire analytic for $\sigma^\phi$, whose linear span is weakly dense in $M$ and such that $\{a_n\xi\mid n\geq1\}$ has dense linear span in $L^2(M)$. We may construct such a sequence by applying the Gaussian regularizations
$$a\longmapsto(r/\pi)^{1/2}\int_\R e^{-rt^2}\sigma_t^\phi(a)\,\d t \qquad(r\in\N),$$
to a countable strongly dense subset of the unit ball, see \cite[Lemma~VIII.2.3 and its proof, pp.~98--99]{TakII}.

Let $b_n=\sigma_{-i/2}^\phi(a_n^*)$. By the standard-form identities, we have
\begin{equation}\label{eq:nrsg-paired-vectors}
Ja_nJ\xi=b_n\xi,\qquad a_n\xi=Jb_nJ\xi.
\end{equation}
Take $c_n>0$ such that
\begin{equation}\label{eq:nrsg-weight-summability}
\sum_n c_n(\|a_n\|_\infty^2+\|b_n\|_\infty^2)<\infty.
\end{equation}
Rescaling the coefficients, we may assume that
\begin{equation}\label{eq nrsg normalization}
C=\sum_n c_n(a_n^*a_n+b_n^*b_n),\qquad \|C\|_\infty=1,
\end{equation}
where the series converges in norm. In particular, $0\leq C\leq1$. For $\eta\in L^2(M)$, write $\omega_\eta(T)=\langle T\eta,\eta\rangle$, and set
$$\psi=\left.\sum_n c_n(\omega_{a_n\xi}+\omega_{b_n\xi})\right|_M.$$
Since $\omega_{a_n\xi}+\omega_{b_n\xi} =\omega_{Jb_nJ\xi}+\omega_{Ja_nJ\xi}$ by (\ref{eq:nrsg-paired-vectors}) and $Ja_nJ,Jb_nJ\in M'$, for $x\in M_+$ we have
$$\psi(x)=\langle JCJx^{1/2}\xi,x^{1/2}\xi\rangle \leq\phi(x).$$
By \cite[Lemma~2.10(1)]{Ha75}, let $\eta\in L^2(M)_+$ with $\phi-\psi=\langle\,\cdot\,\eta,\eta\rangle$, and define
$$\varphi=\sum_n c_n(\omega_{a_n\xi}+\omega_{b_n\xi})+\omega_\eta.$$
Then $\varphi|_M=\phi$. Since $J\eta=\eta$ and $\omega_v(JT^*J)=\omega_{Jv}(T)$, the state $\varphi$ is symmetric and hence regular. The totality of $(a_n\xi)$ makes $\varphi$ faithful on $\BL$.

By the same argument as in the proof of Proposition~\ref{exists Pphi}, the inequality $\phi-\psi\leq\phi$ implies that the map $x\xi\mapsto x\eta$ extends to a contraction $JwJ\in M'$ with $w\in M$. Then
$$w\xi=\eta=JwJ\xi,\qquad \|w\|_\infty\leq1.$$
The two coefficient families associated with $\varphi$ are
$$(z_k)=\big(\sqrt{c_n}b_n,\sqrt{c_n}a_n\big)_{n\geq1}\cup\{w\},\qquad (y_k)=\big(\sqrt{c_n}a_n^*,\sqrt{c_n}b_n^*\big)_{n\geq1}\cup\{w^*\}.$$
The two representations are
\begin{equation}\label{eq:nrsg-vector-representations}
\varphi(T) =\sum_k\langle TJz_kJ\xi,Jz_kJ\xi\rangle =\sum_k\langle Ty_k^*\xi,y_k^*\xi\rangle \qquad(T\in\BL).
\end{equation}
Since $\varphi|_M=\phi$ and $\varphi(Ja^*J)=\phi(a)$ for $a\in M$,
$$\sum_kz_k^*z_k=1,\qquad\sum_ky_ky_k^*=1.$$
Moreover, since $a_n\in\C[(z_k)]$ and $a_n^*\in\C[(y_k)]$, we have $\overline{\C[(z_k)]}^{\,w}=\overline{\C[(y_k)]}^{\,w}=M$, so $\varphi$ is nrsg. By (\ref{eq:nrsg-vector-representations}),
$$\phi(\CQ(a))=\sum_k\langle ay_k^*\xi,y_k^*\xi\rangle =\varphi(a)=\phi(a)\qquad(a\in M).$$
This proves (i).

To prove (ii), note that
$$\sigma_{i/2}^\phi(a_n)=b_n^*,\qquad \sigma_{i/2}^\phi(b_n)=a_n^*.$$
It remains to treat $w^*$. For $r>0$, define
$$w_r=(r/\pi)^{1/2}\int_\R e^{-rt^2}\sigma_t^\phi(w)\,\d t.$$
By the identity $Jw\xi=w\xi$ and the commutation relation $J\Delta_\phi^{it}=\Delta_\phi^{it}J$, we have $Jw_r\xi=w_r\xi$. The Tomita identity and the separating property of $\xi$ imply
$$\sigma_{-i/2}^\phi(w_r^*)=w_r,\qquad \sigma_{i/2}^\phi(w_r)=w_r^*.$$
The boundary values on the half-strip are bounded by $\|w\|_\infty$. By the three-lines theorem and the ultraweak closedness of the analytic generator, we have $\sigma_{i/2}^\phi(w)=w^*$, see \cite[p.~192]{Ha79OVW}. Hence (\ref{eq symmetric coefficient involution}) exchanges the paired coefficients and fixes $w^*$.

Moreover, by the choice of $c_n$ in (\ref{eq:nrsg-weight-summability}) and the definition of $(y_k)$, we have
$$\sum_k\|y_k\|_\infty^2 =\sum_n c_n(\|a_n\|_\infty^2+\|b_n\|_\infty^2)+\|w\|_\infty^2<\infty.$$
This proves (ii).

(iii): Assume that $\phi$ is tracial. The modular action $\sigma^\phi$ is trivial, so every element of $M$ is entire analytic. Since $M$ is linearly spanned by its unitaries, we may take $a_n\in\mathcal U(M)$ in the preceding construction with the same density properties. For this choice, $b_n=a_n^*$, and hence $C=2\left(\sum_n c_n\right)1$. By (\ref{eq nrsg normalization}), $2\sum_n c_n=1$. Since each $\omega_{a_n\xi}$ restricts to $\phi$ on $M$, we have $\psi=\phi$. Hence $\eta=0$ and $w=0$. Omitting the zero coefficient, the inner ucp coefficients of $\CQ$ are
$$y_{2n-1}=\sqrt{c_n}a_n^*,\qquad y_{2n}=\sqrt{c_n}a_n\qquad(n\geq1).$$
This has the required form $y_k=\lambda_ku_k$, with $\sum_k\lambda_k^2=2\sum_n c_n=1$, which finishes the proof.
\end{proof}

For the rest of this section, we fix a normal faithful state $\phi$ on $M$ and an nrsg $(M,\phi)$-hyperstate $\varphi$ (Definition~\ref{def strong generation}) with Poisson boundary $(\CB_\varphi,\zeta)$.

We prove ergodicity, double ergodicity, and rigidity for Poisson boundaries. The proofs are in the same spirit as Das--Peterson \cite{DP20}.

\begin{proposition}[cf.\ {\cite[Propositions~2.5(i) and~2.9]{DP20}}]\label{prop generation consequences}
For $k\geq0$, define the \textbf{$k$-th convolution} $\varphi^{*k}$ by $\varphi^{*k}(T)=\langle\mathcal P_\varphi^k(T)\xi_\phi,\xi_\phi\rangle$. Then the normal state $\overline\varphi=\sum_{n=0}^{\infty}2^{-n-1}\varphi^{*n}$ is faithful on $\BL$. In particular, $\zeta=\varphi|_{\har(\CP_\varphi)}$ is a normal faithful state on the Poisson boundary $\CB_\varphi$. Moreover, $M$ is the largest von Neumann subalgebra of $\har(\CP_\varphi)$ for the original operator product on $\BL$.
\end{proposition}

\begin{proof}
Let $T\in\BL_+$ satisfy $\overline\varphi(T)=0$. Since all the coefficients are positive, $\varphi^{*n}(T)=0$ for every $n\geq0$. Expanding the iterates of $\CP_\varphi$ and using positivity, we have
$$\langle T(JwJ)\xi_\phi,(JwJ)\xi_\phi\rangle=0$$
for every word $w$ in $(z_k)$, including the empty word. The weak density of $\C[(z_k)]$ implies that the vectors $(JwJ)\xi_\phi$ have dense linear span in $L^2(M)$. Hence $T=0$, proving that $\overline\varphi$ is faithful. Since $\varphi^{*n}(T)=\varphi(T)$ for $T\in\har(\CP_\varphi)$ and $n\geq0$, we have $\zeta=\overline\varphi|_{\har(\CP_\varphi)}$. Hence $\zeta$ is normal and faithful on $\CB_\varphi$.

Let $\mathcal N\subset\har(\CP_\varphi)$ be a von Neumann subalgebra for the original operator product, and let $v_k=Jz_k^*J$. For $T\in\mathcal N$, the operators $T,T^*,TT^*$ are fixed by $\CP_\varphi$, so
\begin{equation}\label{vT-Tv}
\sum_k(v_kT-Tv_k)(v_kT-Tv_k)^*
=\CP_\varphi(TT^*)-\CP_\varphi(T)T^*-T\CP_\varphi(T^*)+TT^*=0.
\end{equation}
Hence $[T,v_k]=0$ for every $k$. By applying (\ref{vT-Tv}) to $T^*$, we also have $[T,v_k^*]=0$ for every $k$. By the weak density of $\C[(z_k)]$, we have $T\in(JMJ)'=M$. Since $M$ itself is contained in the harmonic space, it is the largest such subalgebra.
\end{proof}

The following is a well-known result \cite[Lemma~2.6]{DP20}, which also appears in earlier literature such as \cite{FNW94}, \cite[Lemma~3.4]{BJKW00}, and \cite[Lemma~3.1]{CD20}. For completeness, we include the proof.

\begin{lemma}\label{lem:faithful-state-fixed-points}
Let $\CA$ be a von Neumann algebra, and let $\CP:\CA\to\CA$ be a normal ucp map preserving a normal faithful state $\psi_\CA$ on $\CA$. Then $\har(\CP,\CA)=\{x\in\CA\mid \CP(x)=x\}$ is a von Neumann subalgebra of $\CA$.
\end{lemma}

\begin{proof}
Since $\har(\CP,\CA)$ is a weakly closed operator system, by polarization we only need to prove that $aa^*\in\har(\CP,\CA)$ for every $a\in\har(\CP,\CA)$. For such an $a$, Kadison's inequality gives $\CP(aa^*)-aa^*=\CP(aa^*)-\CP(a)\CP(a)^*\geq0$. Since $\psi_\CA\circ\CP=\psi_\CA$, we have $\psi_\CA(\CP(aa^*)-aa^*)=0$. By faithfulness, $\CP(aa^*)=aa^*$, which finishes the proof.
\end{proof}

\begin{proposition}[cf.\ {\cite[Proposition~2.7]{DP20}}]\label{ergodicity of PB}
We have the \textbf{ergodicity} of the Poisson boundary, that is, $M'\cap \CB_\varphi=\mathcal{Z}(M)$.
\end{proposition}

\begin{proof}
By identifying $\CB_\varphi$ with $\har(\CP_\varphi)$, we consider elements in $M' \cap \har(\CP_\varphi)$. Note that $M' \cap B(L^2(M)) = JMJ$. Hence, we only need to compute $JMJ \cap \har(\CP_\varphi)$.

Write $\phi^o=\varphi|_{JMJ}$. Since $\CP_\varphi$ is $JMJ$-inner, we have $\CP_\varphi(JMJ)\subset JMJ$. Since $\varphi$ is regular, for any $T\in JMJ$, we have
$$\phi^o(\CP_\varphi(T))=\langle\CP_\varphi(T)\xi_\phi,\xi_\phi\rangle =\varphi(T)=\phi^o(T).$$
Here the first equality uses $\CP_\varphi(T)\in JMJ$ and the fact that $\xi_\phi$ implements $\phi^o$ on $JMJ$. Hence, the restriction $\CP_\varphi|_{JMJ}$ preserves the normal faithful state $\phi^o$ on $JMJ$. By Lemma~\ref{lem:faithful-state-fixed-points}, the fixed-point space $\har(\CP_\varphi|_{JMJ},JMJ)$ is a von Neumann subalgebra of $JMJ$.

Because $\varphi$ is strongly generating, by Proposition~\ref{prop generation consequences}, $M$ is the largest von Neumann subalgebra contained in $\har(\CP_\varphi)$. Therefore, $\har(\CP_\varphi) \cap JMJ \subset M$, which implies $\har(\CP_\varphi) \cap JMJ \subset M \cap JMJ = \mathcal{Z}(M)$. Since $\mathcal{Z}(M) \subset \har(\CP_\varphi)$ is trivial, we conclude $M'\cap \CB_\varphi=\mathcal{Z}(M)$.
\end{proof}

We next prove the double ergodicity of the Poisson boundary.

\begin{proposition}\label{double ergodicity}
Let $M \subset \CA$ be a W$^*$-inclusion. Using the coefficients $(y_k)$ of $\CQ_\varphi$ as in (\ref{eq:intro-inner-ucp-maps}), we define the normal $M$-inner ucp map $\CQ:\CA\to\CA$ by $\CQ(T)=\sum_k y_kTy_k^*$. Let
$$\har(\CQ, \CA) =\{T\in \CA\mid \CQ(T)=T\}.$$
Assume there exists a normal faithful ucp map $\Phi: \CA \to \CB_\varphi$ such that $\Phi|_M = \id_M$. Then we have
$$\har(\CQ, \CA) = M' \cap \CA.$$

In particular, we have the \textbf{double ergodicity} of the Poisson boundary (cf.\ \cite[Theorem~3.1]{DP20})
$$ \har(\CP_\varphi) \cap \har(\CQ_\varphi) = \mathcal{Z}(M). $$
\end{proposition}

\begin{proof}
Let $\zeta$ be the normal faithful stationary state on the Poisson boundary $\CB_\varphi$. We define a normal faithful state $\psi_\CA$ on $\CA$ by $\psi_\CA = \zeta \circ \Phi$.

Because $\Phi|_M = \id_M$, $\Phi$ is an $M$-bimodular map. Hence $\Phi\circ \CQ=\CQ_\varphi\circ\Phi$. Indeed, for every $T\in\CA$,
\begin{equation}\label{eq:PB-map-intertwining}
\Phi(\CQ(T))
=\sum_{k=1}^\infty\Phi(y_kTy_k^*)
=\sum_{k=1}^\infty y_k\Phi(T)y_k^*
=\CQ_\varphi(\Phi(T)).
\end{equation}
We now prove that $\zeta$ is $\CQ_\varphi$-invariant on $\CB_\varphi$. Identifying $\CB_\varphi$ with $\har(\CP_\varphi)$, the state is implemented by the vector state $\zeta(\,\cdot\,) =\lag \CP_\varphi(\,\cdot\,) \xi_\phi, \xi_\phi \rag= \lag\,\cdot\,\xi_\phi, \xi_\phi \rag$, which coincides with the original hyperstate $\varphi$ on this space. By the defining construction of $\CQ_\varphi$, we have $\lag \CQ_\varphi(T)\xi_\phi, \xi_\phi \rag = \varphi(T)$ for any $T\in\BL$. Hence, for any $S \in \CB_\varphi$, we have
\begin{equation}\label{eq:boundary-stationarity}
\zeta(\CQ_\varphi(S))
=\lag\CQ_\varphi(S)\xi_\phi,\xi_\phi\rag
=\varphi(S)=\zeta(S).
\end{equation}

By (\ref{eq:PB-map-intertwining}) and (\ref{eq:boundary-stationarity}), for $T\in\CA$,
$$ \psi_\CA(\CQ(T)) = \zeta(\Phi(\CQ(T))) = \zeta(\CQ_\varphi(\Phi(T))) = \zeta(\Phi(T)) = \psi_\CA(T). $$
This shows that $\CQ$ preserves the normal faithful state $\psi_\CA$. By Lemma~\ref{lem:faithful-state-fixed-points}, the fixed-point space $\CA_0 = \har(\CQ, \CA)$ is a von Neumann subalgebra of $\CA$.

Clearly, if $T \in M' \cap \CA$, then $\CQ(T) = \sum_{k=1}^\infty y_k T y_k^* = \sum_{k=1}^\infty y_k y_k^* T = T$, so $M' \cap \CA \subset \CA_0$.

Conversely, suppose $T \in \CA_0$. Because $\CA_0$ is a von Neumann subalgebra, we also have $TT^* \in \CA_0$, meaning $\CQ(T) = T$ and $\CQ(TT^*) = TT^*$. We compute the following sum of positive elements
\begin{equation}\label{yT-Ty}
\begin{aligned}
&\sum_{k=1}^\infty (y_k T - T y_k)(y_k T - T y_k)^* \\
=& \sum_{k=1}^\infty (y_k T T^* y_k^* - y_k T y_k^* T^* - T y_k T^* y_k^* + T y_k y_k^* T^*) \\
=& \CQ(TT^*) - \CQ(T)T^* - T\CQ(T^*) + T \left(\sum_{k=1}^\infty y_k y_k^*\right) T^* \\
=& TT^* - TT^* - TT^* + TT^* = 0.
\end{aligned}
\end{equation}
By (\ref{yT-Ty}), $[y_k,T]=0$ for every $k$. Since $\overline{\C[(y_k)]}^{\,w}=M$, $T\in M'\cap\CA$. Hence $\har(\CQ,\CA)=M'\cap\CA$.

For the double ergodicity of $\CB_\varphi$, we set $\CA = \CB_\varphi \cong \har(\CP_\varphi)$ and let $\Phi = \id_{\CB_\varphi}$. The hypothesis is satisfied, so we have
$$ \har(\CP_\varphi) \cap \har(\CQ_\varphi) = \har(\CQ, \CB_\varphi) = M' \cap \CB_\varphi. $$
By Proposition~\ref{ergodicity of PB}, $M'\cap\CB_\varphi=\mathcal Z(M)$, which concludes the proof.
\end{proof}

By taking the Ces\`aro means of $\CQ$, we have the inner UCP-Dixmier property.

\begin{corollary}\label{Cor inner UCP}
Let $M \subset \CA$ be a W$^*$-inclusion, and assume there exists a normal faithful ucp map $\Phi: \CA \to \CB_\varphi$ such that $\Phi|_M = \id_M$. Then $M \subset \CA$ has the inner UCP-Dixmier property.
\end{corollary}
\begin{proof}
Let $\CQ$ be as in Proposition~\ref{double ergodicity}, and let $\mathcal C_n=(1/n)\sum_{k=0}^{n-1}\CQ^k$. Let $\omega\in\beta\mathbb N\setminus\mathbb N$ be a free ultrafilter. By Proposition~\ref{double ergodicity}, we can define a conditional expectation $E:\CA\to M'\cap\CA$ by
$$ E(x) := \operatorname*{wo-lim}_{n \to \omega} \mathcal C_n(x). $$

Recall that $\CQ$ is an $M$-inner ucp map. Therefore, $E(x)$ lies in the weak closure of the inner UCP orbit of $x$:
$$ E(x) \in \overline{\{ \CP(x) \mid \CP \text{ is an } M\text{-inner ucp map} \}}^{w} \cap(M'\cap\CA). $$
Hence $M\subset\CA$ has the inner UCP-Dixmier property.
\end{proof}

We now prove the rigidity of the Poisson boundary, one of the most important properties of Poisson boundaries in both the group setting and the Das--Peterson framework. Although this result is not used in the subsequent proofs, we include it for completeness of our treatment of Poisson boundary theory.

\begin{proposition}[cf.\ {\cite[Theorem~4.1]{DP20}}]\label{rigidity of PB}
Let $\mathcal{C}$ be a weakly closed operator system such that $M\subset\mathcal{C}\subset \CB_\varphi$ and $M\mathcal{C}+\mathcal{C}M\subset \mathcal{C}$. If $\Psi:\mathcal{C}\to\CB_\varphi$ is a normal ucp map such that $\Psi|_M=\id_M$, then $\Psi=\id_\mathcal{C}$.
\end{proposition}

\begin{proof}
We identify $\mathcal{C}$ with its image under the isomorphism $\CB_\varphi\cong\har(\CP_\varphi)$. Then $\Psi$ is a normal ucp map from $\mathcal{C}$ to $\har(\CP_\varphi)$. Since $\Psi|_M = \id_M$, the map $\Psi$ is $M$-bimodular. Using the explicit formula $\CQ_\varphi(T)=\sum_{k=1}^{\infty}y_kTy_k^*$, for any $T\in\mathcal{C}$ we have
$$\lag\Psi(\CQ_\varphi(T))\xi_\phi,\xi_\phi\rag = \sum_{k=1}^\infty \lag \Psi(y_kTy_k^*)\xi_\phi,\xi_\phi\rag = \sum_{k=1}^\infty \lag y_k\Psi(T)y_k^*\xi_\phi,\xi_\phi\rag = \lag \CQ_\varphi(\Psi(T))\xi_\phi,\xi_\phi\rag.$$
Since $\Psi(T) \in \har(\CP_\varphi)$ and $\varphi$ is regular, the term on the right evaluates to
$$ \lag \CP_\varphi(\Psi(T))\xi_\phi,\xi_\phi\rag = \lag \Psi(T)\xi_\phi,\xi_\phi\rag. $$
Hence, $\lag\Psi(\CQ_\varphi(T))\xi_\phi,\xi_\phi\rag = \lag\Psi(T)\xi_\phi,\xi_\phi\rag$ for all $T\in\mathcal{C}$. For the Ces\`aro means $\mathcal C_n=(1/n)\sum_{k=0}^{n-1}\CQ_\varphi^k$, linearity gives
$$ \left\lag \Psi(\mathcal C_n(T))\xi_\phi, \xi_\phi \right\rag = \lag\Psi(T)\xi_\phi,\xi_\phi\rag. $$
Let $\omega\in\beta\N\setminus\N$ be a free ultrafilter. By normality and Proposition~\ref{double ergodicity}, we have
$$z:=\operatorname*{wo-lim}_{n\to\omega}\mathcal C_n(T)\in\har(\CP_\varphi)\cap\har(\CQ_\varphi)=\mathcal Z(M).$$

Since $z\in M$, we have $\Psi(z)=z$. By normality of $\Psi$, we have the first equality below. The last equality follows from (\ref{eq:boundary-stationarity}) and normality of $\zeta$:
$$ \lag\Psi(T)\xi_\phi,\xi_\phi\rag = \lag\Psi(z)\xi_\phi,\xi_\phi\rag = \lag z\xi_\phi,\xi_\phi\rag = \lag T\xi_\phi,\xi_\phi\rag. $$
Now, for any $a,b\in M$ and $T\in\mathcal{C}$, the element $b^*Ta$ is still in $\mathcal{C}$. Using the $M$-bimodularity of $\Psi$, we have
$$ \lag\Psi(T)a\xi_\phi, b\xi_\phi\rag = \lag\Psi(b^*Ta)\xi_\phi,\xi_\phi\rag = \lag b^*Ta\xi_\phi,\xi_\phi\rag = \lag Ta\xi_\phi, b\xi_\phi\rag. $$
Since $M\xi_\phi$ is dense in $L^2(M)$, we have $\Psi(T)=T$ for all $T\in\mathcal{C}$. Hence, $\Psi=\id_\mathcal{C}$.
\end{proof}

\section{Weak Dixmier properties via Poisson boundaries}\label{sec:isono-theorem}

In this section, we use Poisson boundaries to give a sufficient condition for a W$^*$-inclusion $M\subset\CA$ to have the weak Dixmier property. This relates the weak Dixmier property to Poisson boundaries and the associated inner ucp maps. We begin by proving Theorem~\ref{letterthm inner UCP eq to Dixmier}.

\begin{theorem}[Theorem~\ref{letterthm inner UCP eq to Dixmier}]\label{Thm inner UCP eq to Dixmier}
Let $M \subset \CA$ be a W$^*$-inclusion such that $M$ is properly infinite. Then $M \subset \CA$ has the weak Dixmier property if and only if it has the inner UCP-Dixmier property. Here $\CA$ need not have separable predual.
\end{theorem}
\begin{proof}
We use a dilation argument similar to that of Haagerup \cite[Proposition~2.1 and the proof of Theorem~2.2]{Haagerup1985Injectivity}, passing from completely positive maps to unitary conjugations by using proper infiniteness. We only need to prove the if direction, since every finite convex combination of unitary conjugations is an $M$-inner ucp map. Assume that $M \subset \CA$ has the inner UCP-Dixmier property.

Because $M$ is properly infinite, it contains an infinite sequence of mutually orthogonal projections $(p_k)_{k=1}^\infty$ such that $p_k \sim 1_M$ for all $k$ and $\sum_{k=1}^\infty p_k = 1_M$. Take isometries $v_k\in M$ with $v_kv_k^*=p_k$. By this partition of unity, we have a spatial isomorphism
$$\pi:\CA\xrightarrow{\sim}\CA\bar\otimes B(\mathcal{H}),\qquad \pi(a)=(v_i^*av_j)_{i,j\geq1},$$
where $\mathcal{H}=\ell^2(\N)$. Its restriction maps $M$ onto $M\bar\otimes B(\mathcal{H})$ and $p_k$ to $1_M\otimes e_{kk}$, where $(e_{ij})_{i,j\geq1}$ are the standard matrix units of $B(\mathcal{H})$. Indeed, on a Hilbert space on which $\CA$ acts, the unitary $\xi\otimes e_k\mapsto v_k\xi$ implements this isomorphism. Hence we identify the inclusion $M\subset\CA$ with
$$M\bar\otimes B(\mathcal{H})\subset\CA\bar\otimes B(\mathcal{H}).$$
We only need to prove the weak Dixmier property for an arbitrary element $x\in\CA\bar\otimes B(\mathcal{H})$. Write
$$K_x=\overline{\operatorname{conv}}^{w} \{UxU^*\mid U\in\mathcal U(M\bar\otimes B(\mathcal{H}))\}.$$

Since $B(\mathcal{H})$ is hyperfinite, take an increasing sequence of finite-dimensional subalgebras $F_n\subset B(\mathcal{H})$ with $\bigvee_nF_n=B(\mathcal{H})$. For a free ultrafilter $\omega$ on $\N$, by normalized Haar averaging we have
$$a\otimes1_{B(\mathcal{H})}:=\operatorname*{wo-lim}_{n\to\omega}\int_{\mathcal U(F_n)}
(1_M\otimes u)x(1_M\otimes u^*)\,\d u\in K_x\qquad(a\in\CA).$$
Since $K_x$ is weakly closed, convex and invariant under these unitary conjugations, it also contains the weakly closed convex hull of the unitary orbit of $a\otimes1_{B(\mathcal{H})}$.

We only need to prove that for any $M$-inner ucp map
$$\CP(a)=\sum_{k=1}^\infty z_kaz_k^*,\qquad z_k\in M,\qquad\sum_{k=1}^\infty z_kz_k^*=1_M,$$
the element $\CP(a)\otimes1_{B(\mathcal{H})}$ is contained in the weak closure of the convex hull of the unitary orbit of $a\otimes1_{B(\mathcal{H})}$.

Fix an integer $n\geq1$. We construct a partial isometry $Z_n\in M\bar\otimes B(\mathcal{H})$ structured as a block matrix where each block is of size $n\times n$. The first block-row consists of a zero block followed by scaled identity matrices $z_k\mathrm I_n$, and all subsequent block-rows are identically zero:
$$Z_n=
\begin{pmatrix}
0_{n\times n}&z_1\mathrm I_n&z_2\mathrm I_n&z_3\mathrm I_n&\cdots\\
0_{n\times n}&0_{n\times n}&0_{n\times n}&0_{n\times n}&\cdots\\
0_{n\times n}&0_{n\times n}&0_{n\times n}&0_{n\times n}&\cdots\\
\vdots&\vdots&\vdots&\vdots&\ddots
\end{pmatrix}.$$

The strong convergence of $\sum_{k=1}^\infty z_kz_k^*=1_M$ makes this row a bounded operator. Let $E_n=\sum_{j=1}^n e_{jj}$. Then
$$Z_nZ_n^*=1_M\otimes E_n.$$
Hence, $Z_n^*$ is a partial isometry with initial projection $1_M\otimes E_n$. Because the first block-column of $Z_n$ is identically zero, its final projection $Z_n^*Z_n$ is orthogonal to $1_M\otimes E_n$. Since $M$ is properly infinite, each $1_M\otimes e_{jj}$ is equivalent to $1_{M\bar\otimes B(\mathcal{H})}$. Hence
$$1-Z_n^*Z_n\geq1_M\otimes E_n \geq1_M\otimes e_{11}\sim1_{M\bar\otimes B(\mathcal{H})}.$$
We also have
$$1-(1_M\otimes E_n)\geq 1_M\otimes e_{n+1,n+1}\sim1_{M\bar\otimes B(\mathcal{H})}.$$
Therefore, the orthogonal complements of both projections are Murray--von Neumann equivalent to the identity:
$$1-(1_M\otimes E_n)\sim1-Z_n^*Z_n \sim1_{M\bar\otimes B(\mathcal{H})}.$$
Hence there exists a partial isometry $W_n\in M\bar\otimes B(\mathcal{H})$ with
$$W_n^*W_n=1-Z_n^*Z_n,\qquad W_nW_n^*=1-(1_M\otimes E_n).$$
Define $U_n=Z_n+W_n$. Since the initial and final projections of $Z_n$ and $W_n$ are respectively orthogonal, all mixed terms vanish. Therefore,
$$U_n^*U_n=Z_n^*Z_n+W_n^*W_n=1,\qquad U_nU_n^*=Z_nZ_n^*+W_nW_n^*=1.$$
Hence we complete the partial isometry into a unitary $U_n\in\mathcal U(M\bar\otimes B(\mathcal{H}))$. Since $(1_M\otimes E_n)W_n=0$, the first block-row of $U_n$ remains the first block-row of $Z_n$. In block-row notation, we write
$$U_n=\begin{pmatrix}Z_n\\W_n\end{pmatrix}.$$

We conjugate $a\otimes1_{B(\mathcal{H})}$ by $U_n$. Since
$$Z_n(a\otimes1_{B(\mathcal{H})})Z_n^*=\CP(a)\otimes E_n,$$
we have, with respect to the first $n$ coordinates and their orthogonal complement,
\begin{equation}\label{eq:unitary-block-conjugate}
U_n(a\otimes1_{B(\mathcal{H})})U_n^*=
\begin{pmatrix}
\CP(a)\mathrm I_n&\ast\\
\ast&\ast
\end{pmatrix}.
\end{equation}
Let $b_n=U_n(a\otimes1_{B(\mathcal{H})})U_n^*-\CP(a)\otimes E_n$. By (\ref{eq:unitary-block-conjugate}), we have $(1_\CA\otimes E_n)b_n(1_\CA\otimes E_n)=0$. Since $\|b_n\|\leq2\|a\|$ and $E_n\nearrow1_{B(\mathcal{H})}$ strongly, we have
$$b_n=(1-1_\CA\otimes E_n)b_n+(1_\CA\otimes E_n)b_n(1-1_\CA\otimes E_n) \longrightarrow0\quad\text{weakly}\qquad(n\to\infty).$$
Therefore,
$$\operatorname*{wo-lim}_{n\to\infty} U_n(a\otimes1_{B(\mathcal{H})})U_n^*=\CP(a)\otimes1_{B(\mathcal{H})}.$$

Because this holds for any $M$-inner ucp map $\CP$, the entire set
$$\{\CP(a)\otimes1_{B(\mathcal{H})}\mid \CP\text{ is an }M\text{-inner ucp map}\}$$
is contained in $K_x$.

By the assumption that the original inclusion $M\subset\CA$ has the inner UCP-Dixmier property, there exists an element
$$y\in\overline{\{\CP(a)\mid \CP\text{ is an }M\text{-inner ucp map}\}}^{w}\cap(M'\cap\CA).$$
Because the mapping $a\mapsto a\otimes1_{B(\mathcal{H})}$ is normal, the element $y\otimes1_{B(\mathcal{H})}$ lies in the weak closure of $\{\CP(a)\otimes1_{B(\mathcal{H})}\}$. Since $K_x$ is already weakly closed, it contains $y\otimes1_{B(\mathcal{H})}$.

Hence, $y\otimes1_{B(\mathcal{H})}$ belongs to
$$(M\bar\otimes B(\mathcal{H}))'\cap(\CA\bar\otimes B(\mathcal{H})) =(M'\cap\CA)\bar\otimes\C1_{B(\mathcal{H})}$$
and lies in the weakly closed convex hull of the unitary orbit of $x$, which finishes the proof.
\end{proof}

Combining Theorem~\ref{Thm inner UCP eq to Dixmier} with Corollary~\ref{Cor inner UCP}, we have the weak Dixmier property for inclusions that map into a Poisson boundary.

\begin{theorem}[Theorem~\ref{letterthm ucp to PB}]\label{Thm ucp to PB}
Let $M$ be a von Neumann algebra with separable predual and $\varphi$ be an nrsg hyperstate on $\BL$ with the Poisson boundary $\CB_\varphi$. Assume that one of the following conditions holds:
\begin{itemize}
    \item [(i)] $M$ is properly infinite, i.e., $M=M_{\mathrm{I}_\wx}\oplus M_{\mathrm{II}_\wx}\oplus M_{\mathrm{III}}$;
    \item [(ii)] $M$ is finite with normal faithful trace $\tau$ and $\varphi=\int_{\mathcal{U}(M)}\lag\,\cdot\,u\xi_\tau,u\xi_\tau\rag\d\mu(u)$ for some $\mu\in\prob(\mathcal{U}(M))$.
\end{itemize}
Then any W$^*$-inclusion $M\subset \CA$ equipped with a normal faithful ucp map $\Phi:\CA\to\CB_\varphi$ with $\Phi|_M=\id_M$ has the weak Dixmier property.
\end{theorem}

\begin{proof}
Let $M \subset \CA$ be a W$^*$-inclusion equipped with a normal faithful ucp map $\Phi: \CA \to \CB_\varphi$ satisfying $\Phi|_M = \id_M$. By Corollary \ref{Cor inner UCP}, $M \subset \CA$ has the inner UCP-Dixmier property.

We consider the two cases separately.

\textbf{Case (i):} Assume $M$ is properly infinite. By Theorem \ref{Thm inner UCP eq to Dixmier}, for properly infinite von Neumann algebras, the weak Dixmier property is equivalent to the inner UCP-Dixmier property. Since $M \subset \CA$ already has the inner UCP-Dixmier property by Corollary \ref{Cor inner UCP}, it has the weak Dixmier property.

\textbf{Case (ii):} Assume $M$ is finite with trace $\tau$, and $\varphi = \int_{\mathcal{U}(M)} \lag\,\cdot\,u\xi_\tau, u\xi_\tau \rag \d\mu(u)$ for some probability measure $\mu \in \prob(\mathcal{U}(M))$. In this setting, the associated ucp map $\CQ : \CA \to \CA$ is given by $ \CQ(x) = \int_{\mathcal{U}(M)} u^* x u \,\d\mu(u) \quad \text{for all } x \in \CA. $ The adjoint direction follows from the vector-state convention $\varphi=\int\omega_{u\xi_\tau}\,\d\mu(u)$. Because $\CQ$ is an integral of unitary conjugations, the elements $\CQ^k(x)$ and their Ces\`aro means $\mathcal C_n(x)=(1/n)\sum_{k=0}^{n-1}\CQ^k(x)$ lie in the weakly closed convex hull of the unitary orbit of $x$:
$$ \mathcal C_n(x) \in \overline{\operatorname{conv}}^{w} \{ u x u^* \mid u \in \mathcal{U}(M) \}. $$
By the proof of Corollary \ref{Cor inner UCP}, any weak limit point $z$ of the sequence $(\mathcal C_n(x))_{n=1}^\infty$ is an element of $M' \cap \CA$. Since the weakly closed convex hull of the unitary orbit is weakly closed, the limit point $z$ remains within it. Hence
$$ z \in \overline{\operatorname{conv}}^{w} \{ u x u^* \mid u \in \mathcal{U}(M) \} \cap (M' \cap \CA) \neq \emptyset. $$
This proves that $M \subset \CA$ has the weak Dixmier property.

\end{proof}

As a first application, we recover the weak Dixmier property for inclusions with expectation.

\begin{corollary}[Corollary~\ref{lettercor expectation}, \cite{Pop99,Mar19}]\label{cor expectation}
Let $M\subset\CA$ be an inclusion of countably decomposable von Neumann algebras with a normal faithful conditional expectation. Then $M\subset\CA$ has the weak Dixmier property.
\end{corollary}
\begin{proof}
Let $E:\CA\to M$ be the given normal faithful conditional expectation. As in \cite[p.~122, proof of the main theorem]{Mar19}, take an increasing net $(M_i)_{i\in I}$ of von Neumann subalgebras with separable preduals and normal faithful conditional expectations $E_{M_i}:M\to M_i$, such that $\bigvee_iM_i=M$. For $x\in\CA$, let
$$K_x^i=\overline{\operatorname{conv}}^{w}\{uxu^*\mid u\in\mathcal U(M_i)\}.$$
We only need to prove that $K_x^i\cap M_i'\ne\varnothing$ for every $i$: by weak compactness, any choice $y_i\in K_x^i\cap M_i'$ has a weak cluster point in
$$\overline{\operatorname{conv}}^{w}\{uxu^*\mid u\in\mathcal U(M)\}\cap(M'\cap\CA).$$
Since $E_{M_i}\circ E:\CA\to M_i$ is a normal faithful conditional expectation, we may therefore assume that $M$ has separable predual. No separability assumption on $\CA$ is needed in the argument below or in the auxiliary averaging results used here. Decompose $M=M_{\mathrm{fin}}\oplus M_\infty$ into its finite and properly infinite parts, omitting zero summands, and write $p_1,p_2$ for their units. Averaging by $p_1+e^{2\pi it}p_2$ removes the off-diagonal terms of every $x\in\CA$:
\begin{equation}\label{eq:central-block-average}
\int_0^1(p_1+e^{2\pi it}p_2)x(p_1+e^{2\pi it}p_2)^*\,\d t
=p_1xp_1+p_2 xp_2.
\end{equation}
The restrictions $E_i:p_i\CA p_i\to p_iM$ are normal faithful conditional expectations for $i=1,2$.

Take a normal faithful trace $\tau$ on $M_{\mathrm{fin}}$, a countable adjoint-closed family of unitaries $(u_k)$ generating $M_{\mathrm{fin}}$, and $c_k>0$ with $\sum_kc_k=1$. The state
$$\varphi_1=\sum_k c_k\langle\,\cdot\,u_k\xi_\tau,u_k\xi_\tau\rangle$$
is normal and regular. Both strong generation conditions follow from the weak density of $\C[(u_k)]$, since $(u_k)$ is adjoint-closed. Composing $E_1$ with $M_{\mathrm{fin}}\subset\CB_{\varphi_1}$, by Theorem~\ref{Thm ucp to PB}(ii), we have that $M_{\mathrm{fin}}\subset p_1\CA p_1$ has the weak Dixmier property.

For $M_\infty$, take an nrsg hyperstate $\varphi_2$ by Example~\ref{exists nrsg hyperstate}. Composing $E_2$ with $M_\infty\subset\CB_{\varphi_2}$ and applying Theorem~\ref{Thm ucp to PB}(i), we have the same conclusion for this corner. Taking product nets of the two unitary averaging nets after applying (\ref{eq:central-block-average}), we have an element of $M'\cap\CA$ in the closed convex unitary orbit of $x$.
\end{proof}

As a semifinite analogue of Corollary~\ref{cor expectation}, Isono \cite[Theorem~A]{Iso25} proved that, for a W$^*$-inclusion $M\subset\CA$ with a normal faithful semifinite operator-valued weight $T:\CA_+\to\widehat M_+$, every $x\in\CA_+$ with $T(x)$ bounded satisfies
$$\overline{\operatorname{conv}}^{w}\{uxu^*\mid u\in\mathcal U(M)\} \cap(M'\cap\CA)\ne\emptyset.$$
To generalize Isono's result, we first introduce a generalization of such weights.
\begin{definition}\label{def:M-bimodular-semifinite-map}
Let $M\subset\CA$ and $M\subset\CB$ be W$^*$-inclusions. A map $T:\CA_+\to\widehat{\CB}_+$ is called an \textbf{$M$-bimodular semifinite positive map} if it is additive, positively homogeneous and satisfies
$$T(a^*xa)=a^*T(x)a\qquad(x\in\CA_+,\ a\in M).$$
We require normality, meaning preservation of increasing bounded nets, faithfulness, and ultraweak density of
$$\mathfrak m_{T}^+=\{x\in\CA_+\mid T(x)\in \CB_+\}.$$
Here $\widehat{\CB}_+$ is the extended positive cone, and $T(x)<\infty$ will always mean that $T(x)$ is bounded.
\end{definition}

For every $\rho\in \CB_*^+$, the function $\rho\circ T$ is a normal weight. In particular, its sublevel sets are ultraweakly closed. It follows that, for $C<\infty$, the set
\begin{equation}\label{eq:weight-sublevel}
\{x\in\CA_+\mid T(x)\leq C1\}
\end{equation}
is ultraweakly closed. We use this scalar interpretation of lower semicontinuity throughout, see \cite[Sections~1--2]{Ha79OVW}.

The following lemma is a semifinite analogue of Lemma~\ref{lem:faithful-state-fixed-points}. When the invariant state is replaced by an invariant semifinite weight, Ces\`aro limits still satisfy the following multiplicativity.

\begin{lemma}\label{lemma weight preserving}
Let $\psi_\CA$ be a normal faithful semifinite weight on a von Neumann algebra $\CA$, and let $\CP:\CA\to\CA$ be a normal ucp map satisfying $\psi_\CA\circ \CP=\psi_\CA$. Let $\mathcal C_n=(1/n)\sum_{k=0}^{n-1}\CP^k$. Fix a free ultrafilter $\omega$ and define
$$E(x)=\operatorname*{wo-lim}_{n\to\omega}\mathcal C_n(x).$$
If $x\in\CA_+$ and $\psi_\CA(x)<\infty$, then $\CP(E(x))=E(x)$ and $\CP(E(x)^2)=E(x)^2$.
\end{lemma}
\begin{proof}
Let $y=E(x)$. Since $\|\CP(\mathcal C_n(x))-\mathcal C_n(x)\|\le2\|x\|/n$, by normality of $\CP$, we have $\CP(y)=y$. By lower semicontinuity and invariance of the weight, we have
$$\psi_\CA(y)\le\psi_\CA(x)<\infty, \qquad \psi_\CA(y^2)\le\|y\|\psi_\CA(y)<\infty.$$
By Kadison's inequality, $\CP(y^2)-y^2\geq0$. Since $\psi_\CA\circ \CP=\psi_\CA$ and $\psi_\CA(y^2)<\infty$,
$$\psi_\CA(\CP(y^2)-y^2) =\psi_\CA(\CP(y^2))-\psi_\CA(y^2)=0.$$
By faithfulness, we have $\CP(y^2)=y^2$.
\end{proof}

By this lemma, we have the following semifinite version of Proposition~\ref{double ergodicity}.

\begin{theorem}\label{stronger double ergodicity}
Let $M\subset\CA$ be a W$^*$-inclusion. Assume that $M$ has separable predual, and let $\varphi$ be an nrsg $M$-hyperstate. Assume that there exists a normal faithful semifinite $M$-bimodular positive map $T:\CA_+\to\widehat{\CB_\varphi}_+$ in the preceding sense. Define the normal $M$-inner ucp map $\CQ:\CA\to\CA$ by $\CQ(x)=\sum_k y_kxy_k^*$. For any $x\in\CA_+$ with $T(x)<\infty$, every ultraweak cluster point of the Ces\`aro means $\mathcal C_n(x)=(1/n)\sum_{k=0}^{n-1}\CQ^k(x)$ belongs to $M'\cap\CA$.
\end{theorem}
\begin{proof}
Let $\zeta$ be the canonical normal faithful state on $\CB_\varphi$ and let $\psi_\CA=\zeta\circ T$. This is a normal faithful weight. It is semifinite because its finite positive domain contains the ultraweakly dense cone $\mathfrak m_{T}^+$.

By normality and $M$-bimodularity, in the extended positive cone we have
$$T(\CQ(x))=\sum_k y_kT(x)y_k^*=\CQ_\varphi(T(x)).$$
By (\ref{eq:boundary-stationarity}), $\zeta\circ\CQ_\varphi=\zeta$ on $(\CB_\varphi)_+$. This identity extends to the extended positive cone by monotone approximation. Hence $\psi_\CA\circ \CQ=\psi_\CA$.

Let $a$ be an ultraweak cluster point of $(\mathcal C_n(x))$. By the proof of Lemma~\ref{lemma weight preserving}, applied along a convergent subnet, we have $\CQ(a)=a$ and $\CQ(a^2)=a^2$. Hence,
$$\sum_k(y_ka-ay_k)(y_ka-ay_k)^* =\CQ(a^2)-\CQ(a)a-a\CQ(a)+a^2=0.$$
Hence $[a,y_k]=0$ for every $k$. Since $\overline{\C[(y_k)]}^{\,w}=M$, $a\in M'\cap\CA$.
\end{proof}

As in Theorem~\ref{Thm ucp to PB}, we can convert these Ces\`aro limits into unitary averages.

\begin{theorem}\label{Thm semifinite map}
Let $M$ have separable predual and let $\varphi$ be an nrsg hyperstate. Assume that one of the following conditions holds:
\begin{itemize}
\item[(i)] $M$ is properly infinite;
\item[(ii)] $M$ is finite with a normal faithful tracial state $\tau$ and
\begin{equation}\label{eq:unitary-measure-hyperstate}
\varphi(a)=\int_{\mathcal U(M)}
\langle au\xi_\tau,u\xi_\tau\rangle\,\mathrm d\mu(u)
\end{equation}
for some $\mu\in\prob(\mathcal U(M))$.
\end{itemize}
Let $M\subset\CA$ be a W$^*$-inclusion and let $T:\CA_+\to\widehat{\CB_\varphi}_+$ be a normal faithful semifinite $M$-bimodular positive map. Then, for $x\in\CA_+$ with $T(x)<\infty$,
$$\overline{\operatorname{conv}}^{w} \{uxu^*\mid u\in\mathcal U(M)\}\cap(M'\cap\CA)\ne\emptyset.$$
\end{theorem}
\begin{proof}
Let $\CQ$ and its Ces\`aro means $\mathcal C_n$ be as in Theorem~\ref{stronger double ergodicity}. Let $C=\|T(x)\|$ and $K_x=\overline{\operatorname{conv}}^{w}\{uxu^*\mid u\in\mathcal U(M)\}$. For $u\in\mathcal U(M)$, $T(uxu^*)=uT(x)u^*\leq C1$. By the ultraweak closedness of (\ref{eq:weight-sublevel}),
$$T(b)\leq C1\qquad(b\in K_x).$$

Suppose first that $M$ is properly infinite. Take isometries $(v_i)\subset M$ with $v_i^*v_j=\delta_{ij}1$ and $\sum_i v_iv_i^*=1$, and use the isomorphism
$$\pi:\CA\xrightarrow{\sim}\CA\bar\otimes B(\ell^2(\mathbb N)), \qquad \pi(b)=\sum_{i,j}v_i^*bv_j\otimes e_{ij}.$$
It maps $M$ onto $M\bar\otimes B(\ell^2(\mathbb N))$. As in the proof of Theorem~\ref{Thm inner UCP eq to Dixmier}, by unitary averaging over the type I tensor factor, we have $a\otimes1\in\pi(K_x)$ for some $a\in\CA_+$. Let $b=\pi^{-1}(a\otimes1)$. Since $b\in K_x$ and $a=v_1^*bv_1$,
$$T(a)=v_1^*T(b)v_1\le C1.$$
By Theorem~\ref{stronger double ergodicity}, an ultraweak cluster point $c$ of $(\mathcal C_n(a))$ belongs to $M'\cap\CA$.

The matrix-dilation argument in the proof of Theorem~\ref{Thm inner UCP eq to Dixmier} applies to each $M$-inner ucp map. Hence $\mathcal C_n(a)\otimes1$ belongs to the weakly closed unitary convex hull of $a\otimes1$ for every $n$. That hull is contained in $\pi(K_x)$, so $c\otimes1\in\pi(K_x)$. Since $c\otimes1\in(M\bar\otimes B(\ell^2(\mathbb N)))'$,
$$\pi^{-1}(c\otimes1)\in K_x\cap(M'\cap\CA).$$

Suppose now that $M$ is finite and $\varphi$ satisfies (\ref{eq:unitary-measure-hyperstate}). With the present vector-state convention,
$$\CQ(b)=\int_{\mathcal U(M)}u^*bu\,\mathrm d\mu(u).$$
Its iterates and their Ces\`aro means $\mathcal C_n$ preserve $K_x$. The ultraweak cluster points of $(\mathcal C_n(x))$ belong to $M'\cap\CA$ by Theorem~\ref{stronger double ergodicity}, which proves the assertion.
\end{proof}

Composing an operator-valued weight onto $M$ with the boundary inclusion recovers Isono's result.

\begin{corollary}[\cite{Iso25}]\label{cor OVW}
Let $M\subset\CA$ be an inclusion of countably decomposable von Neumann algebras, and let $T:\CA_+\to\widehat M_+$ be a normal faithful semifinite operator-valued weight. For every $x\in\CA_+$ such that $T(x)$ is bounded,
$$\overline{\operatorname{conv}}^{w} \{uxu^*\mid u\in\mathcal U(M)\}\cap(M'\cap\CA)\ne\emptyset.$$
\end{corollary}
\begin{proof}
As in the beginning of the proof of Corollary~\ref{cor expectation}, we may assume that $M$ has separable predual. Write $M=M_{\mathrm{fin}}\oplus M_\infty$ for its finite and properly infinite parts, with units $p_1,p_2$. Averaging over the unitaries $p_1+e^{2\pi it}p_2$ sends $x$ to $p_1xp_1+p_2 xp_2$ within its unitary convex hull.

The restrictions $T_i:(p_i\CA p_i)_+\to\widehat{p_iM}_+$ are normal, faithful and semifinite. Indeed, compression of the ultraweakly dense finite positive domain of $T$ is ultraweakly dense in $(p_i\CA p_i)_+$, and
$$T_i(p_iap_i)=p_iT(a)p_i.$$
Moreover, $T_i(p_ixp_i)$ is bounded.

On $M_{\mathrm{fin}}$ take an atomic measure on a countable symmetric family of unitaries generating $M_{\mathrm{fin}}$, with every atom having positive mass. Then we have an nrsg hyperstate $\varphi_1$ with respect to a normal faithful tracial state. On $M_\infty$ take an nrsg hyperstate $\varphi_2$ by Example~\ref{exists nrsg hyperstate}. Compose each $T_i$ with the extension to extended positive cones of the canonical embedding $p_iM\subset\CB_{\varphi_i}$. These composites satisfy the hypotheses of Theorem~\ref{Thm semifinite map}. No normalization on $p_iM$ is required. Apply that theorem separately to the two positive diagonal blocks.

Finally, block-diagonal unitaries in $\mathcal U(M_{\mathrm{fin}})\oplus\mathcal U(M_\infty)$ belong to $\mathcal U(M)$. Taking product nets of the two averaging nets, we have a limit in $M'\cap\CA$, still in the weakly closed convex hull of the unitary orbit of $x$.
\end{proof}

\section{The relative bicentralizer conjecture}
\label{sec relative bicentralizer}

In this section, we prove Theorem~\ref{letterthm relative bicentralizer}. Fix an inclusion $M\subset P$ of countably decomposable von Neumann algebras with a normal faithful conditional expectation $E:P\to M$. We only need to prove that $M\subset c(P)$ has the weak Dixmier property.

By the same argument as in the beginning of the proof of Corollary~\ref{cor expectation}, we may assume that $M$ has separable predual, see also \cite[Theorems~D and~E(6)]{Ma25}. For the remainder of this section, we therefore assume that $M$ has separable predual, while $P$ need only be countably decomposable. Fix a normal faithful state $\phi$ on $M$, and let $\rho=\phi\circ E$. We use the compatible realizations
$$c_\phi(M)\subset c_\rho(P)=P\rtimes_{\sigma^\rho}\R.$$
We write $\phi^{it}$ for the canonical implementing unitaries in $c_\phi(M)$. Let $\CA=c_\rho(P)$, and let $\psi_\CA=\rho\circ T$ be the dual weight of $\rho$, where $T:\CA_+\to\widehat P_+$ is the canonical operator-valued weight. We use the notation $h=\frac{\d\psi_\CA}{\d\tau_\CA}$ and $H=\log h$ from Subsection~\ref{sec:core-weights}, where $\tau_\CA$ is the canonical trace on $\CA$. Here $\lambda_t=h^{it}=\phi^{it}$. The modular action $\sigma_t^{\psi_\CA}=\operatorname{Ad}(\lambda_t)$ restricts to $\sigma_t^\phi$ on $M$. Take the $(M,\phi)$-hyperstate $\varphi$ and the inner ucp coefficients $(y_k)$ from Example~\ref{exists nrsg hyperstate}, and define
\begin{equation}\label{eq:core-inner-map}
\CQ(x)=\sum_k y_kxy_k^*\qquad(x\in\CA).
\end{equation}

To prove that $M\subset\CA$ has the inner UCP-Dixmier property, we only need to prove that $d:=\CQ(aa^*)-aa^*=0$ for every $a\in\har(\CQ,\CA)$. As in the proof of Proposition~\ref{double ergodicity} and Lemma~\ref{lemma weight preserving}, followed by Ces\`aro averaging, we have
$$
\begin{aligned}
&\CQ(aa^*)=aa^*\quad(a\in\har(\CQ,\CA))\\
\Longrightarrow\quad &\har(\CQ,\CA)=M'\cap\CA\\
\Longrightarrow\quad &M\subset\CA\text{ has the inner UCP-Dixmier property}.
\end{aligned}
$$

Fix $a\in\har(\CQ,\CA)$. By Kadison's inequality, $d\geq0$. Lemma~\ref{lem:positive-localization} reduces $d=0$ to proving that every $x\in\CA_+$ with $0\leq x\leq d$ and $\psi_\CA(x)<\infty$ must satisfy $x=0$. Such an $x$ satisfies
$$\sum_{k=0}^{n-1}\CQ^k(x)\leq\sum_{k=0}^{n-1}\CQ^k(d) =\CQ^n(aa^*)-aa^*\leq C1\qquad(n\geq1),$$
with $C=\|a\|^2$.

To explain how this bound and $\psi_\CA(x)<\infty$ imply $x=0$, we first present a commutative model, which is not used in the proof of Theorem~\ref{letterthm relative bicentralizer}. This model illustrates the classical cutoff argument for recurrence of symmetric Markov operators, see \cite[Section~1.6, Theorem~1.6.3]{FOT11}. In this model, $L^\infty(\R,\d t)$, equipped with the weight induced by Lebesgue measure, plays the role of $(\CA,\psi_\CA)$, and a symmetric measure $\nu\in\prob(\R)$ plays the role of the symmetric nrsg hyperstate $\varphi$.

\begin{proposition}\label{prop:commutative-recurrence-model}
Let $\nu\in\prob(\R)$ be symmetric, and assume that $m_2=\int_\R s^2\,\d\nu(s)<\infty$. Define the Markov operator $Q$ on $L^\infty(\R,\d t)$ by
$$Q(b)(t)=\int_\R b(t+s)\,\d\nu(s)\qquad(b\in L^\infty(\R,\d t)).$$
Let $X\in L^\infty(\R,\d t)_+$ with $\int_\R X(t)\,\d t<\infty$. If there exists $0\leq C<\infty$ such that
\begin{equation}\label{eq:commutative-bounded-sums}
\sum_{k=0}^{n-1}Q^k(X)\leq C1\qquad(n\geq1),
\end{equation}
then $X=0$.
\end{proposition}

\begin{proof}
We first prove that $Q$ preserves the integral and is a self-adjoint contraction on $L^2(\R)$ (cf. Lemma~\ref{lem:pairing}). By Tonelli's theorem and translation invariance, we have $\int_\R Q(b)\,\d t=\int_\R b\,\d t$ for $b\in L^1(\R)_+$. By Kadison's inequality, we then have $\|Q(b)\|_2\leq\|b\|_2$ for $b\in L^2\cap L^\infty(\R)$, and hence on $L^2(\R)$ by density. Since $\nu$ is symmetric, $Q$ is self-adjoint, so $1-Q\geq0$.

We may assume that $\int_\R X\,\d t>0$. Take a real $f\in C_c^\infty(\R)$ with $0\leq f\leq1$, equal to $1$ near $0$, even and nonincreasing in $|t|$. Let $F_R(t)=f(t/R)$ for $R\geq1$. Then $F_R\nearrow1$ pointwise $(R\to\infty)$ and $\int_\R F_R\,\d t<\infty$.

We next prove the partial-sum estimate (cf. Lemma~\ref{lem:partial-sum-estimate})
\begin{equation}\label{eq:commutative-pairing-bound}
\int_\R X\,\d t\leq\liminf_{R\to\infty}C\|(1-Q)^{1/2}F_R\|_2^2.
\end{equation}
Let $a_n=\sum_{k=0}^{n-1}Q^k(X)$. Then $0\leq a_n\leq C1$, $\int_\R a_n\,\d t=n\int_\R X\,\d t$, and $a_n\in L^1\cap L^\infty(\R)\subset L^2(\R)$. Since $(1-Q)a_n=X-Q^n(X)$,
$$\|(1-Q)^{1/2}a_n\|_2^2=\int_\R a_n\bigl(X-Q^n(X)\bigr)\,\d t\leq C\int_\R X\,\d t.$$
For fixed $R$, by the Cauchy--Schwarz inequality we therefore have
\begin{equation}\label{eq:commutative-partial-cutoff}
\left|\int_\R F_R\bigl(X-Q^n(X)\bigr)\,\d t\right|^2\leq C\left(\int_\R X\,\d t\right)\|(1-Q)^{1/2}F_R\|_2^2.
\end{equation}
By (\ref{eq:commutative-bounded-sums}) and monotone convergence,
$$\sum_{n=0}^\infty\int_\R F_RQ^n(X)\,\d t\leq C\int_\R F_R\,\d t<\infty.$$
Hence $\int_\R F_RQ^n(X)\,\d t\to0$ $(n\to\infty)$, so by (\ref{eq:commutative-partial-cutoff}) we have
$$\left(\int_\R F_RX\,\d t\right)^2\leq C\left(\int_\R X\,\d t\right)\|(1-Q)^{1/2}F_R\|_2^2.$$
By dominated convergence as $R\to\infty$ and division by $\int_\R X\,\d t>0$, we have (\ref{eq:commutative-pairing-bound}).

We now prove the cutoff estimate (cf. Lemma~\ref{lem:cutoff-estimates}). Expanding the squared difference, for $b\in L^2(\R)$ we have
\begin{equation}\label{eq:commutative-difference}
\|(1-Q)^{1/2}b\|_2^2=\frac12\int_\R\|b(\,\cdot\,+s)-b\|_2^2\,\d\nu(s).
\end{equation}
By the fundamental theorem of calculus and the Cauchy--Schwarz inequality, we have
\begin{equation}\label{eq:commutative-translation}
\|F_R(\,\cdot\,+s)-F_R\|_2\leq|s|\|F_R'\|_2=\frac{|s|}{\sqrt R}\|f'\|_2\qquad(s\in\R).
\end{equation}
Recall that $m_2=\int_\R s^2\,\d\nu(s)<\infty$. By (\ref{eq:commutative-difference}) and (\ref{eq:commutative-translation}), we have
\begin{equation}\label{eq:commutative-cutoff-bound}
\|(1-Q)^{1/2}F_R\|_2^2\leq\frac{\|f'\|_2^2}{2R}\int_\R s^2\,\d\nu(s)=\frac{m_2\|f'\|_2^2}{2R}.
\end{equation}
Hence $\|(1-Q)^{1/2}F_R\|_2^2\to0$ $(R\to\infty)$, so by (\ref{eq:commutative-pairing-bound}) we have $\int_\R X\,\d t=0$ and $X=0$.
\end{proof}

\begin{remark}\label{rem:commutative-symmetry}
Symmetry of $Q$ is essential in Proposition~\ref{prop:commutative-recurrence-model}. Indeed, let $Q_0(f)(t)=f(t-1)$ be the nonsymmetric translation on $L^\infty(\R)$, and let $X=\mathds{1}_{[0,1)}$. Then $X\ne0$, but
$$\sum_{n=0}^\infty Q_0^n(X)=\mathds{1}_{[0,\infty)}\leq1.$$
This is why our proof of Theorem~\ref{letterthm relative bicentralizer} uses the symmetric nrsg hyperstate constructed in Example~\ref{exists nrsg hyperstate}, whereas the arguments in Section~\ref{sec:isono-theorem} and later in Section~\ref{sec:kadison} work with an arbitrary nrsg hyperstate.
\end{remark}

We now return to $\CQ:\CA\to\CA$ defined by (\ref{eq:core-inner-map}) and to $x\in\CA_+$ with $\psi_\CA(x)<\infty$ and $\sum_{k=0}^{n-1}\CQ^k(x)\leq C1$ $(n\geq1)$. The proof of $x=0$ follows the three steps in the proof of Proposition~\ref{prop:commutative-recurrence-model}. We first prove in (\ref{eq:weight-preserving}) that $\CQ$ preserves $\psi_\CA$. By the symmetry (\ref{eq symmetric coefficient involution}) of the coefficients $(y_k)$, we then construct in Lemma~\ref{lem:pairing} a self-adjoint contraction $S$ on $L^2(\CA,\tau_\CA)$, which plays the role of $Q$ on $L^2(\R)$. Lemmas~\ref{lem:partial-sum-estimate} and~\ref{lem:cutoff-estimates} are the analogues of (\ref{eq:commutative-pairing-bound}) and (\ref{eq:commutative-cutoff-bound}), respectively. Together they prove that $\psi_\CA(x)=0$. By faithfulness of $\psi_\CA$, we have $x=0$.

Hence $d=0$ for every $a\in\har(\CQ,\CA)$, and $M\subset\CA$ has the inner UCP-Dixmier property. By Theorem~\ref{Thm inner UCP eq to Dixmier}, we then have the weak Dixmier property when $M$ is properly infinite. For finite $M$, we take $\phi$ tracial so that the Ces\`aro means are already unitary averages. The general case follows by decomposing $M$ into its finite and properly infinite parts.

Since $\phi\circ \CQ=\phi$, by normality and $M$-bimodularity of $E$ and $T$, we have
\begin{equation}\label{eq:weight-preserving}
\psi_\CA\circ \CQ =\phi\circ E\circ T\circ \CQ =\phi\circ \CQ\circ E\circ T =\phi\circ E\circ T=\psi_\CA,
\end{equation}
where the compositions are understood on the extended positive cone.

Define
$$\mathfrak m_{\psi_\CA}
=\operatorname{span}\{a\in\CA_+\mid \psi_\CA(a)<\infty\},\qquad
\xi(a)=h^{1/4}ah^{1/4}\quad(a\in\mathfrak m_{\psi_\CA}).$$
On $L^2(\CA,\tau_\CA)$ we use the inner product $\langle X,Y\rangle_2=\tau_\CA(Y^*X)$, linear in the first variable.

Define $S:D\to D$, where $D=\{\xi(a)\mid a\in\mathfrak m_{\psi_\CA}\}$, by
\begin{equation}\label{eq:implementation}
S\xi(a)=\xi(\CQ(a))\qquad(a\in\mathfrak m_{\psi_\CA}).
\end{equation}
This is well-defined since $\xi$ is injective and $\CQ(\mathfrak m_{\psi_\CA})\subset\mathfrak m_{\psi_\CA}$ by (\ref{eq:weight-preserving}).

\begin{lemma}\label{lem:pairing}\label{lem:kms}
For every $a\in\mathfrak m_{\psi_\CA}$, we have $\xi(a)\in L^2(\CA,\tau_\CA)$. Let $x\in\CA_+$ satisfy $\psi_\CA(x)<\infty$. Then
$$\omega_x(a)=\tau_\CA(h^{1/2}xh^{1/2}a)\qquad(a\in\CA)$$
defines a finite normal positive functional on $\CA$, with $\omega_x(1)=\psi_\CA(x)$, and
\begin{equation}\label{eq:pairing}
\langle \xi(a),\xi(x)\rangle_2=\omega_x(a)
\qquad(a\in\mathfrak m_{\psi_\CA}).
\end{equation}
Moreover, $S$ extends to a self-adjoint contraction on $L^2(\CA,\tau_\CA)$.
\end{lemma}
\begin{proof}

For $a\in\mathfrak m_{\psi_\CA}$, let
$$\iota(a)=h^{1/2}ah^{1/2}\in L^1(\CA,\tau_\CA).$$
We first prove the assertions about $\xi$ and $\omega_x$. Both $ah^{1/2}$ and $h^{1/2}a$ belong to $L^2(\CA,\tau_\CA)$. By traciality and the Cauchy--Schwarz inequality, we have
\begin{equation}\label{eq symmetric embedding bound}
\begin{aligned}
\|\xi(a)\|_2^2
&=\tau_\CA(h^{1/4}a^*h^{1/2}ah^{1/4})\\
&=\tau_\CA\bigl((ah^{1/2})^*(h^{1/2}a)\bigr)\\
&\leq\|ah^{1/2}\|_2\|h^{1/2}a\|_2\\
&=\tau_\CA(h^{1/2}a^*ah^{1/2})^{1/2}
  \tau_\CA(h^{1/2}aa^*h^{1/2})^{1/2}\\
&=\psi_\CA(a^*a)^{1/2}\psi_\CA(aa^*)^{1/2}.
\end{aligned}
\end{equation}
We first make this calculation for $a_n=e_nae_n$, where $e_n=\mathds{1}_{[1/n,n]}(h)$. The convergences $a_nh^{1/2}\to ah^{1/2}$ and $h^{1/2}a_n\to h^{1/2}a$ in $L^2$ justify passage to the limit.

Since $\iota(x)\in L^1(\CA,\tau_\CA)_+$ and $\tau_\CA(\iota(x))=\psi_\CA(x)$, $\omega_x$ is a finite normal positive functional with $\omega_x(1)=\psi_\CA(x)$. To prove (\ref{eq:pairing}), let
$$e_n=\mathds{1}_{[\frac1n,n]}(h),\qquad x_n=e_nxe_n,\qquad a_n=e_nae_n.$$
By (\ref{eq:tail}), the projections $e_n$ have finite trace. By traciality in $e_n\CA e_n$, we have
$$\langle \xi(a_n),\xi(x_n)\rangle_2 =\tau_\CA(h^{1/2}x_nh^{1/2}a_n).$$
Since $e_nh=he_n$ and $e_n\nearrow1$ strongly,
$$
\iota(x_n)=e_n\iota(x)e_n\longrightarrow \iota(x)\quad\text{in }L^1,\qquad \xi(x_n)=e_n\xi(x)e_n\longrightarrow\xi(x)\quad\text{in }L^2,\qquad \xi(a_n)\longrightarrow\xi(a)\quad\text{in }L^2.
$$
Moreover,
$$\sup_n\|a_n\|\leq\|a\|,\qquad a_n\longrightarrow a\quad\text{$\ast$-strongly}.$$
Passing to the limit proves (\ref{eq:pairing}).

We next prove that $S$ extends to a self-adjoint contraction. Identifying $\CA_*$ with $L^1(\CA,\tau_\CA)$, we first show that the preadjoint of $\CQ$ satisfies
\begin{equation}\label{eq:kms}
\CQ_*\iota(a)=\iota(\CQ(a))\qquad(a\in\mathfrak m_{\psi_\CA}).
\end{equation}
Let $v_k=\sigma_{i/2}^\phi(y_k^*)$. By Example~\ref{exists nrsg hyperstate}, $(v_k)$ is a permutation of $(y_k)$. We have
\begin{equation}\label{eq:half-strip-products}
y_k^*h^{1/2}=h^{1/2}v_k,\qquad
h^{1/2}y_k=v_k^*h^{1/2}.
\end{equation}
Formula (\ref{eq:half-strip-products}) follows first on compact spectral vectors for $H$, and then for the closed measurable products. For $a\in\CA_+$ with $\psi_\CA(a)<\infty$, we have
$$y_k^*\iota(a)y_k=\iota(v_kav_k^*).$$
By (\ref{eq:half-strip-products}), the positive partial sums $\sum_{k\leq n}y_k^*\iota(a)y_k$ converge in $L^1$, since
$$\tau_\CA\left(\sum_{k\leq n}y_k^*\iota(a)y_k\right) =\tau_\CA\left(\iota(a)\sum_{k\leq n}y_ky_k^*\right) \nearrow\psi_\CA(a).$$
By monotone convergence of the corresponding quadratic forms, we have
$$\CQ_*\iota(a)=\sum_k y_k^*\iota(a)y_k =\iota\left(\sum_kv_kav_k^*\right)=\iota(\CQ(a)).$$
This proves (\ref{eq:kms}) for positive $a$, and the general case follows by linearity.

By (\ref{eq:kms}) and (\ref{eq:pairing}), for $a,b\in\mathfrak m_{\psi_\CA}\cap\CA_+$,
$$
\langle S\xi(a),\xi(b)\rangle_2 =\tau_\CA(\iota(b)\CQ(a))=\tau_\CA(\CQ_*\iota(b)a) =\tau_\CA(\iota(\CQ(b))a)=\langle\xi(a),S\xi(b)\rangle_2.
$$
By linearity, we have $\langle Sv,w\rangle_2=\langle v,Sw\rangle_2$ for $v,w\in D$. Also $S(D)\subset D$ by (\ref{eq:weight-preserving}). For $a\in\mathfrak m_{\psi_\CA}$, by (\ref{eq:implementation}), (\ref{eq symmetric embedding bound}), Kadison's inequality for $\CQ^n$, and $\psi_\CA\circ \CQ^n=\psi_\CA$, we have
\begin{align*}
\|S^n \xi(a)\|_2^2
&=\|\xi(\CQ^n(a))\|_2^2\\
&=\|h^{1/4}\CQ^n(a)h^{1/4}\|_2^2\\
&=\tau_\CA(h^{1/2}\CQ^n(a)^*h^{1/2}\CQ^n(a))\\
&\leq\psi_\CA((\CQ^n(a))^*\CQ^n(a))^{1/2}
       \psi_\CA(\CQ^n(a)(\CQ^n(a))^*)^{1/2}\\
&\leq\psi_\CA(\CQ^n(a^*a))^{1/2}
       \psi_\CA(\CQ^n(aa^*))^{1/2}\\
&=\psi_\CA(a^*a)^{1/2}\psi_\CA(aa^*)^{1/2}.
\end{align*}
Hence
$$C_v:=\sup_{m\geq0}\|S^m v\|_2<\infty\qquad(v\in D).$$
For $v\in D$, by symmetry, we have
\begin{equation}\label{eq:S-square-bound}
\|Sv\|_2^2=\langle S^2v,v\rangle_2
\leq\|S^2v\|_2\|v\|_2.
\end{equation}
Iterating (\ref{eq:S-square-bound}), we have
$$\|Sv\|_2 \leq\|S^{2^n}v\|_2^{1/2^n}\|v\|_2^{1-1/2^n} \leq C_v^{1/2^n}\|v\|_2^{1-1/2^n} \longrightarrow\|v\|_2.$$

It remains to check that $D$ is dense. With $e_n=\mathds{1}_{[1/n,n]}(h)$, both $h$ and $h^{-1}$ are bounded on $e_n$, and $\tau_\CA(e_n)<\infty$. Every bounded $b\in e_n\CA e_n$ has the form
$$b=\xi(h^{-1/4}bh^{-1/4}),\qquad h^{-1/4}bh^{-1/4}\in e_n\CA e_n\subset\mathfrak m_{\psi_\CA}.$$
These elements are dense in $L^2(\CA,\tau_\CA)$. Therefore $S$ extends to a self-adjoint contraction.
\end{proof}

Let $f\in C_c^\infty(\R)$ be real with $0\leq f\leq1$, equal to $1$ near $0$, even and nonincreasing in $|s|$. Define $H_R=f(H/R)$, $R\geq1$. Then $H_R\nearrow1$ strongly as $R\to\infty$, and by (\ref{eq:core-normalization}), we have
$$\psi_\CA(H_R)=R\int_\R f(s)\,\d s<\infty.$$

\begin{lemma}\label{lem:partial-sum-estimate}
Let $x\in\CA_+$ satisfy $\psi_\CA(x)<\infty$. Assume that there is a constant $0\leq C<\infty$ such that
\begin{equation}\label{eq:bounded-positive-sums}
\sum_{k=0}^{n-1}\CQ^k(x)\leq C1\qquad(n\geq1).
\end{equation}
Then
\begin{equation}\label{eq:finite-weight-pairing-bound}
\psi_\CA(x)
\leq\liminf_{R\to\infty}
C\|(1-S)^{1/2}\xi(H_R)\|_2^2.
\end{equation}
\end{lemma}
\begin{proof}
The case $\psi_\CA(x)=0$ is immediate. Assume that $\psi_\CA(x)>0$. Let $a_n=\sum_{k=0}^{n-1}\CQ^k(x)\leq C1$. By (\ref{eq:weight-preserving}) and (\ref{eq:implementation}),
$$\psi_\CA(a_n)=n\psi_\CA(x)<\infty,\qquad (1-S)\xi(a_n)=\xi(x-\CQ^n(x)).$$
By (\ref{eq:pairing}) and positivity,
\begin{equation}\label{eq:partial-sum-bound}
\|(1-S)^{1/2}\xi(a_n)\|_2^2 =\langle\xi(x-\CQ^n(x)),\xi(a_n)\rangle_2 =\omega_x(a_n)-\omega_{a_n}(\CQ^n(x)) \leq\omega_x(a_n) \leq C\omega_x(1)=C\psi_\CA(x).
\end{equation}
Fix $R\geq1$. By (\ref{eq:pairing}) and positivity, $\omega_x(H_R)=\overline{\omega_{H_R}(x)}=\omega_{H_R}(x)$. Since $a_n-\CQ(a_n)=x-\CQ^n(x)$ and $1-S\geq0$, by (\ref{eq:implementation}) and (\ref{eq:partial-sum-bound}), we have
\begin{equation}\label{eq:partial-sum-cutoff-bound}
\begin{aligned}
&|\omega_x(H_R)-\omega_{H_R}(\CQ^n(x))|^2\\
=&\bigl|\omega_{H_R}(x-\CQ^n(x))\bigr|^2\\
=&\bigl|\langle\xi(x-\CQ^n(x)),\xi(H_R)\rangle_2\bigr|^2\\
=&\bigl|\langle(1-S)\xi(a_n),\xi(H_R)\rangle_2\bigr|^2\\
=&\bigl|\langle(1-S)^{1/2}\xi(a_n),
                  (1-S)^{1/2}\xi(H_R)\rangle_2\bigr|^2\\
\leq&\|(1-S)^{1/2}\xi(a_n)\|_2^2
          \|(1-S)^{1/2}\xi(H_R)\|_2^2\\
\leq&C\psi_\CA(x)\|(1-S)^{1/2}\xi(H_R)\|_2^2.
\end{aligned}
\end{equation}
By (\ref{eq:bounded-positive-sums}), for every $N\geq1$,
$$\sum_{n=0}^{N-1}\omega_{H_R}(\CQ^n(x)) =\omega_{H_R}\left(\sum_{n=0}^{N-1}\CQ^n(x)\right) \leq C\omega_{H_R}(1)=C\psi_\CA(H_R)<\infty.$$
Letting $N\to\infty$, we have
$$\sum_{n=0}^{\infty}\omega_{H_R}(\CQ^n(x)) \leq C\psi_\CA(H_R)<\infty,$$
so $\omega_{H_R}(\CQ^n(x))\to0$ as $n\to\infty$. Letting $n\to\infty$ in (\ref{eq:partial-sum-cutoff-bound}), we have
\begin{equation}\label{eq:cutoff-pairing-bound}
|\omega_x(H_R)|^2
\leq C\psi_\CA(x)\|(1-S)^{1/2}\xi(H_R)\|_2^2
\qquad(R\geq1).
\end{equation}
Since $H_R\nearrow1$ strongly as $R\to\infty$, by normality of $\omega_x$ and (\ref{eq:cutoff-pairing-bound}), we have
$$\psi_\CA(x)^2=|\omega_x(1)|^2 =\lim_{R\to\infty}|\omega_x(H_R)|^2 \leq\liminf_{R\to\infty} C\psi_\CA(x)\|(1-S)^{1/2}\xi(H_R)\|_2^2.$$
Dividing by $\psi_\CA(x)$ proves (\ref{eq:finite-weight-pairing-bound}).
\end{proof}

Let $\nu\in\prob(\R)$ be the spectral distribution of $\log\Delta_\phi$ under $\varphi$, i.e., $\nu(E)=\varphi(\mathds{1}_E(\log\Delta_\phi))$ for every Borel set $E\subset\R$. Define
$$m_2:=\int_\R s^2\,\d\nu(s) =\varphi((\log\Delta_\phi)^2).$$

\begin{lemma}\label{lem:cutoff-estimates}
We have $m_2<\infty$ and, for every $R\geq1$,
\begin{equation}\label{eq:cutoff-quadratic-form}
\|(1-S)^{1/2}\xi(H_R)\|_2^2\leq\frac{m_2\|f'\|_2^2}{2R}.
\end{equation}
Moreover, if $x\in\CA_+$ satisfies $\psi_\CA(x)<\infty$ and $\sum_{k=0}^{n-1}\CQ^k(x)\leq C1$ for some $0\leq C<\infty$ and every $n\geq1$, then $x=0$.
\end{lemma}
\begin{proof}
We first prove that $m_2<\infty$. By (\ref{eq:nrsg-vector-representations}), the Tomita identity, and Example~\ref{exists nrsg hyperstate}(ii), we have
\begin{equation}\label{eq symmetric hyperstate moments}
\varphi(\Delta_\phi)
=\sum_k\langle S_\phi y_k^*\xi_\phi,S_\phi y_k^*\xi_\phi\rangle
=\sum_k\|y_k\xi_\phi\|^2
\leq\sum_k\|y_k\|_\infty^2<\infty.
\end{equation}
Since $J\Delta_\phi J=\Delta_\phi^{-1}$, by the symmetry in (\ref{eq:nrsg-J-symmetry}), we have $\varphi(\Delta_\phi^{-1})=\varphi(\Delta_\phi)<\infty$. By (\ref{eq symmetric hyperstate moments}) and $(\log t)^2\leq t+t^{-1}-2$ $(t>0)$, we have $m_2=\varphi((\log\Delta_\phi)^2)<\infty$.

Let $g\in C_c^\infty(\R)$ be real. By (\ref{eq:core-normalization}), $g(H)\in\mathfrak m_{\psi_\CA}$, and $g(H)\in\CA_{\psi_\CA}$. Expanding the commutators, we have
\begin{equation}\label{eq:cutoff-commutator-expansion}
\sum_k[y_k,g(H)][y_k,g(H)]^*
=\CQ(g(H)^2)-\CQ(g(H))g(H)-g(H)\CQ(g(H))+g(H)^2.
\end{equation}
By Kadison's inequality, (\ref{eq:weight-preserving}), and (\ref{eq:core-normalization}),
$$\psi_\CA(\CQ(g(H))^2)\leq\psi_\CA(\CQ(g(H)^2)) =\psi_\CA(g(H)^2)=\int_\R |g(s)|^2\,\d s<\infty.$$
Hence all four terms on the right of (\ref{eq:cutoff-commutator-expansion}) belong to $\mathfrak m_{\psi_\CA}$. Applying $\psi_\CA$ and using the $\psi_\CA$-centralizer property of $g(H)$, the linear extension of (\ref{eq:pairing}) to self-adjoint elements, and (\ref{eq:implementation}), we have
\begin{equation}\label{eq:commutator-quadratic-form}
\begin{aligned}
\sum_k\psi_\CA([y_k,g(H)][y_k,g(H)]^*)
&=2\bigl(\psi_\CA(g(H)^2)-\psi_\CA(g(H)\CQ(g(H)))\bigr)\\
&=2\bigl(\langle\xi(g(H)),\xi(g(H))\rangle_2
       -\langle S\xi(g(H)),\xi(g(H))\rangle_2\bigr)\\
&=2\|(1-S)^{1/2}\xi(g(H))\|_2^2.
\end{aligned}
\end{equation}
We compute the left-hand side by proving
\begin{equation}\label{eq:commutator-integral}
\sum_k\psi_\CA([y_k,g(H)][y_k,g(H)]^*)
=\int_\R\int_\R|g(r+s)-g(r)|^2\,\d r\,\d\nu(s).
\end{equation}
We use the unitary Fourier transform
$$\widehat g(t)=(2\pi)^{-1/2}\int_\R e^{-itr}g(r)\,\d r.$$
By the crossed-product relation, we have
\begin{equation}\label{eq:commutator-fourier}
[y_k,g(H)]=(2\pi)^{-1/2}\int_\R
\widehat g(t)(y_k-\sigma_t^\phi(y_k))\lambda_t\,\d t.
\end{equation}
With the normalization in (\ref{eq:core-normalization}), for $X=(2\pi)^{-1/2}\int_\R a_t\lambda_t\,\d t=(2\pi)^{-1/2}\int_\R a_t e^{itH}\,\d t$, we have
\begin{equation}\label{eq:dual-weight-plancherel}
\psi_\CA(X^*X)=\int_\R\rho(a_t^*a_t)\,\d t,\qquad
\psi_\CA(XX^*)=\int_\R\rho(a_ta_t^*)\,\d t,
\end{equation}
The identities in (\ref{eq:dual-weight-plancherel}) hold first on the crossed-product Hilbert algebra and extend by its $L^2$ isometry, see \cite{Ha78Dual}. In the first identity in (\ref{eq:dual-weight-plancherel}), the zeroth coefficient is $\int_\R\sigma_{-t}^\rho(a_t^*a_t)\,\d t$, and the modular invariance of $\rho$ removes this twist.

By $\varphi=\sum_k\langle\,\cdot\,y_k^*\xi_\phi,y_k^*\xi_\phi\rangle$ and the definition of $\nu$ as the spectral distribution of $\log\Delta_\phi$ under $\varphi$, we have
\begin{align*}
\sum_k\phi\bigl((y_k-\sigma_t^\phi(y_k))
                  (y_k-\sigma_t^\phi(y_k))^*\bigr)
&=\sum_k\|(1-\Delta_\phi^{it})y_k^*\xi_\phi\|^2\\
&=\sum_k\big\langle(1-\Delta_\phi^{-it})(1-\Delta_\phi^{it})
       y_k^*\xi_\phi,y_k^*\xi_\phi\big\rangle\\
&=\varphi\bigl((1-\Delta_\phi^{-it})(1-\Delta_\phi^{it})\bigr)\\
&=\int_\R(1-e^{-its})(1-e^{its})\,\d\nu(s)\\
&=\int_\R|1-e^{its}|^2\,\d\nu(s).
\end{align*}
By (\ref{eq:commutator-fourier}), (\ref{eq:dual-weight-plancherel}), and scalar Plancherel,
$$
\begin{aligned}
\sum_k\psi_\CA([y_k,g(H)][y_k,g(H)]^*)
&=\int_\R|\widehat g(t)|^2
       \int_\R|1-e^{its}|^2\,\d\nu(s)\,\d t\\
&=\int_\R\int_\R|\widehat g(t)|^2|1-e^{its}|^2\,\d t\,\d\nu(s)\\
&=\int_\R\int_\R|g(r+s)-g(r)|^2\,\d r\,\d\nu(s).
\end{aligned}
$$
This proves (\ref{eq:commutator-integral}). Now take $g(r)=f(r/R)$, so $g(H)=H_R$. Since $\|g(\,\cdot\,+s)-g\|_2 \leq|s|\|g'\|_2 =\frac{|s|}{\sqrt R}\|f'\|_2,$ by (\ref{eq:commutator-integral}), we have
\begin{equation}\label{eq:cutoff-commutator-bound}
\sum_k\psi_\CA([y_k,H_R][y_k,H_R]^*)
\leq\frac{m_2\|f'\|_2^2}{R}.
\end{equation}

Combining (\ref{eq:commutator-quadratic-form}) with $g(H)=H_R$ and (\ref{eq:cutoff-commutator-bound}) proves (\ref{eq:cutoff-quadratic-form}).

Now let $x$ satisfy the hypotheses of Lemma~\ref{lem:partial-sum-estimate}. By (\ref{eq:finite-weight-pairing-bound}) and (\ref{eq:cutoff-quadratic-form}),
$$0\leq\psi_\CA(x) \leq\liminf_{R\to\infty}\frac{Cm_2\|f'\|_2^2}{2R}=0.$$
Hence $\psi_\CA(x)=0$, and by faithfulness of $\psi_\CA$, we have $x=0$.
\end{proof}

We can now identify $\har(\CQ,\CA)$.
\begin{proposition}\label{prop:fixed}
The fixed-point space of $\CQ$ is
$$\har(\CQ,\CA)=M'\cap\CA.$$
\end{proposition}
\begin{proof}
We only need to prove that $\har(\CQ,\CA)\subset M'\cap\CA$. Let $a\in\har(\CQ,\CA)$. By Kadison's inequality,
$$d:=\CQ(aa^*)-aa^*=\CQ(aa^*)-\CQ(a)\CQ(a)^*\geq0.$$
We have
$$\sum_{k=0}^{n-1}\CQ^k(d)=\CQ^n(aa^*)-aa^* \leq\|a\|^2 1\qquad(n\geq1).$$
Suppose that $d\ne0$. By Lemma~\ref{lem:positive-localization}, take $0\ne x\leq d$ with $\psi_\CA(x)<\infty$. Then
$$\sum_{k=0}^{n-1}\CQ^k(x)\leq\sum_{k=0}^{n-1}\CQ^k(d) \leq\|a\|^2 1\qquad(n\geq1).$$
By Lemma~\ref{lem:cutoff-estimates}, we have $x=0$, a contradiction. Hence $\CQ(aa^*)=aa^*$. Applying the same argument to $a^*\in\har(\CQ,\CA)$, we have $\CQ(a^*a)=a^*a$, so $a$ belongs to the multiplicative domain of $\CQ$. Expanding as in (\ref{yT-Ty}), we have
$$\sum_k(y_ka-ay_k)(y_ka-ay_k)^*=0.$$
Hence $[a,y_k]=0$ for every $k$. Since $\overline{\C[(y_k)]}^{\,w}=M$, $a\in M'\cap\CA$.
\end{proof}

Marrakchi proved the relative bicentralizer conjecture, in particular, when $M$ is amenable, has no type $\mathrm{III}_1$ direct summand, is a tensor product of two type $\mathrm{III}_1$ von Neumann algebras, or admits a Cartan subalgebra \cite[Theorem~E(1)--(4)]{Ma25}. For completeness, we give a new proof covering these cases as well. We now finish the proof by applying Theorem~\ref{Thm inner UCP eq to Dixmier}.
\begin{proof}[Proof of Theorem~\ref{letterthm relative bicentralizer}]
As in the beginning of the proof of Corollary~\ref{cor expectation}, we may assume that $M$ has separable predual.

Suppose first that $M$ is either finite or properly infinite. Fix a normal faithful state $\phi$ on $M$, taking $\phi$ tracial when $M$ is finite. Let $\CQ$ be the map in (\ref{eq:core-inner-map}) on $\CA=c_{\phi\circ E}(P)$, with coefficients from Example~\ref{exists nrsg hyperstate}. For $x\in\CA$, define
$$\mathcal C_n(x)=\frac1n\sum_{k=0}^{n-1}\CQ^k(x),\qquad K_x=\overline{\operatorname{conv}}^{\,w} \{uxu^*\mid u\in\mathcal U(M)\}.$$
Since
$$\|\CQ(\mathcal C_n(x))-\mathcal C_n(x)\|\leq\frac{2\|x\|}{n},$$
normality of $\CQ$ and Proposition~\ref{prop:fixed} imply that every ultraweak cluster point of $(\mathcal C_n(x))$ belongs to $M'\cap\CA$.

If $M$ is properly infinite, each $\mathcal C_n$ is an $M$-inner ucp map. Hence $M\subset\CA$ has the inner UCP-Dixmier property, and by Theorem~\ref{Thm inner UCP eq to Dixmier}, we have the weak Dixmier property.

If $M$ is finite, by the tracial choice in Example~\ref{exists nrsg hyperstate}, we have
$$y_k=\lambda_ku_k,\qquad \CQ(x)=\sum_k\lambda_k^2u_kxu_k^*,\qquad \sum_k\lambda_k^2=1.$$
Hence $\CQ(K_x)\subset K_x$ and $\mathcal C_n(x)\in K_x$ for every $n$. Every ultraweak cluster point belongs to $K_x\cap(M'\cap\CA)$, which proves the finite case.

For general $M$, write
$$M=M_{\mathrm{fin}}\oplus M_\infty$$
for its finite and properly infinite parts. We may assume that both summands are nonzero. Let $p_1,p_2$ be their units and let $P_i=p_iPp_i$. There is a normal faithful conditional expectation
$$F:P\to P_1\oplus P_2,\qquad F(a)=p_1ap_1+p_2ap_2.$$
The restrictions
$$E_1=E|_{P_1}:P_1\to M_{\mathrm{fin}},\qquad E_2=E|_{P_2}:P_2\to M_\infty$$
are normal faithful conditional expectations.

Take a normal faithful state $\phi$ on $M$ whose restriction to $M_{\mathrm{fin}}$ is tracial, and let $\rho=\phi\circ E$, $\CA=c_\rho(P)$, and $\rho_i=\rho|_{P_i}/\rho(p_i)$. Since $p_i\in M_\phi$, by bimodularity, we have $p_i\in P_\rho$. Hence $F$ preserves $\rho$ and extends to a normal faithful conditional expectation
$$\widetilde F:\CA\to p_1\CA p_1\oplus p_2\CA p_2,\qquad p_i\CA p_i\cong c_{\rho_i}(P_i).$$
For $x\in\CA$, since $p_1-p_2\in\mathcal U(M)$,
$$\widetilde F(x)=p_1xp_1+p_2xp_2 =\frac12\bigl(x+(p_1-p_2)x(p_1-p_2)\bigr)\in K_x.$$
By the finite and properly infinite cases, we have the weak Dixmier property for $M_{\mathrm{fin}}\subset p_1\CA p_1$ and $M_\infty\subset p_2\CA p_2$, respectively. As in the proof of Corollary~\ref{cor expectation}, taking product nets of the two unitary averaging nets, we have
$$K_x\cap(M'\cap\CA)\ne\emptyset.$$
This proves that $M\subset c(P)$ has the weak Dixmier property. Taking $P=M$, we have the last assertion.
\end{proof}

\begin{proof}[Proof of Corollary~\ref{lettercor Connes bicentralizer}]
By Theorem~\ref{letterthm relative bicentralizer}, $M\subset c(M)$ has the weak Dixmier property. Since $M$ is of type $\mathrm{III}_1$, by \cite[Theorem~3.1]{Haagerup1987Uniqueness}, we have $\mathrm B(M,\phi)=\C1$ for every normal faithful state $\phi$ on $M$.
\end{proof}

\begin{proof}[Proof of Corollary~\ref{lettercor relative bicentralizer}]
By Theorem~\ref{letterthm relative bicentralizer}, we have the weak Dixmier property for $M\subset c(P)$. Hence by \cite[Theorem~D]{Ma25}, we have (\ref{eq:target}) for every normal faithful state $\phi$ on $M$.
\end{proof}

\section{Kadison's problem on MASA}\label{sec:kadison}

In this section, we use the generating inner ucp map associated with an nrsg hyperstate to give a new proof of the solution to Kadison's problem (Theorem~\ref{letterthm kadison}).

Fix an inclusion $M\subset P$ of von Neumann algebras with separable preduals and a normal faithful conditional expectation $E:P\to M$ as in Theorem~\ref{letterthm kadison}.

Following the well-established route of Popa \cite{Popa1981Kadison}, also used in \cite{Ma25,Ma26}, we only need to construct an increasing sequence of finite-dimensional abelian subalgebras $D_1\subset D_2\subset\cdots\subset M$ such that $\|\theta|_{D_n'\cap P}\|\to0$ as $n\to\infty$ for every $\theta\in M_\perp=\{\theta\in P_*\mid\theta|_M=0\}$. Then $A_0=\bigvee_nD_n$ satisfies $\|\theta|_{A_0'\cap P}\|\leq\|\theta|_{D_n'\cap P}\|\to0$ as $n\to\infty$ for every $\theta\in M_\perp$. Hence, by the Hahn--Banach theorem, $A_0'\cap P\subset M$. Any MASA $A\subset M$ containing $A_0$ then satisfies $A'\cap P=A$, which finishes the proof.

Given $\varepsilon>0$ and $\theta_1,\ldots,\theta_s\in M_\perp$, we first construct a single finite-dimensional abelian subalgebra $D\subset M$ such that $\|\theta_j|_{D'\cap P}\|<\varepsilon$ for $1\leq j\leq s$. We will then strengthen this in Proposition~\ref{lem:kadison-refinement} to construct an increasing sequence $D_1\subset D_2\subset\cdots\subset M$ with the required convergence.

As mentioned in Subsection~\ref{sec:intro-kadison}, our main contribution is to use Poisson boundary theory to substantially simplify the construction of $D$ in the continuous properly infinite case. This replaces the passage through a hyperfinite subalgebra used by Marrakchi \cite[Theorem~B]{Ma26}, which in the type $\mathrm{III}_1$ case relies on the ergodicity of the relative bicentralizer flow \cite[Theorem~C]{Ma26}.

Following Popa \cite[Definition~2.2]{Popa2021Ergodic}, we call an inclusion $M\subset P$ \textbf{MASA-ergodic} if there is an abelian subalgebra $A\subset M$ with $A'\cap P=A$. The following lemma reduces the proof of Theorem~\ref{letterthm kadison} and the construction of $D$ with $\|\theta_j|_{D'\cap P}\|<\varepsilon$ to the type $\mathrm I$, finite, and properly infinite cases.

\begin{lemma}\label{lem:corners}
Let $M\subset P$ be as in Theorem~\ref{letterthm kadison}. Then the following hold.
\begin{itemize}
\item[(i)] For every nonzero projection $e\in M$, the restriction $E|_{ePe}:ePe\to eMe$ is a normal faithful conditional expectation, and
\begin{equation}\label{eq:kadison-corner-relative-commutant}
(eMe)'\cap ePe=e(M'\cap P)=e\mathcal Z(M)=\mathcal Z(eMe).
\end{equation}
Moreover, let $(e_k)\subset M$ be a countable family of mutually orthogonal nonzero projections with $\sum_ke_k=1$. If each inclusion $e_kMe_k\subset e_kPe_k$ is MASA-ergodic, then $M\subset P$ is MASA-ergodic.
\item[(ii)] For every finite-dimensional abelian subalgebra $D\subset M$,
$$(D'\cap M)'\cap P=D\vee\mathcal Z(M)\subset M.$$
\item[(iii)] Let $e_1,\ldots,e_m\in M$ be mutually orthogonal nonzero projections with $\sum_ie_i=1$, and let $D_i\subset e_iMe_i$ be finite-dimensional abelian subalgebras. Then $D=\bigoplus_{i=1}^mD_i$ is a finite-dimensional abelian subalgebra of $M$, and
$$\|\theta|_{D'\cap P}\|\leq\sum_{i=1}^m \|\theta|_{D_i'\cap e_iPe_i}\|\qquad(\theta\in P_*).$$
\end{itemize}
\end{lemma}

\begin{proof}
(i) Since $E$ is $M$-bimodular, $E(exe)=eE(x)e$ for $x\in P$. Hence $E|_{ePe}$ is a normal faithful conditional expectation onto $eMe$. By \cite[Propositions~4.2.1 and~4.2.2]{AP17} (see also \cite{Jo83}) and $M'\cap P=\mathcal Z(M)$, we have (\ref{eq:kadison-corner-relative-commutant}). For each $k$, take $A_k\subset e_kMe_k$ maximal abelian in $e_kPe_k$, and let $A=\bigoplus_kA_k\subset M$. Since $e_k\in A$, for every $x\in A'\cap P$ we have $x=\sum_ke_kxe_k$ and $e_kxe_k\in A_k$ for every $k$. Hence $A'\cap P=A$.

(ii) Let $e_1,\ldots,e_r$ be the minimal projections of $D$, and let $N=D'\cap M$, so that
$$D=\bigoplus_{i=1}^r\C e_i,\qquad N=\bigoplus_{i=1}^r e_iMe_i.$$
If $x\in N'\cap P$, then $x$ commutes with every $e_i\in N$, and by (\ref{eq:kadison-corner-relative-commutant}) we have
$$x=\sum_{i=1}^r e_i x e_i,\qquad e_i x e_i\in(e_iMe_i)'\cap e_iPe_i=e_i\mathcal Z(M).$$
Hence $N'\cap P\subset\bigoplus_i e_i\mathcal Z(M)=D\vee\mathcal Z(M)$. The reverse inclusion follows from $D\vee\mathcal Z(M)\subset\mathcal Z(N)$.

(iii) Since $e_i\in D$ for every $i$, we have
$$D'\cap P=\bigoplus_{i=1}^m(D_i'\cap e_iPe_i).$$
The estimate follows by the triangle inequality.

\end{proof}

By Lemma~\ref{lem:corners}(i), we only need to prove Theorem~\ref{letterthm kadison} in the type $\mathrm I$ and continuous cases. We treat the type $\mathrm I$ case in the final proof. For continuous $M$, Lemma~\ref{lem:corners}(iii) reduces the construction of $D$ to the finite and properly infinite cases. We now construct, for any $\theta_1,\ldots,\theta_s\in M_\perp$ and $\varepsilon>0$, a finite-dimensional abelian subalgebra $D\subset M$ such that $\|\theta_j|_{D'\cap P}\|<\varepsilon$ for $1\leq j\leq s$.

\subsection{The continuous properly infinite case.}\label{sec:kadison-continuous-infinite}

We begin with the case where $M$ is continuous and properly infinite, i.e., $M=M_{\mathrm{II}_\infty}\oplus M_{\mathrm{III}}$. We proceed in two steps.
\begin{enumerate}
\item Using our Poisson boundary tools, we first construct a full matrix algebra $F\subset M$ such that $\|\theta_j|_{F'\cap P}\|<\varepsilon/2$ for $1\leq j\leq s$ (Lemma~\ref{lem:kadison-finite-averages} and Proposition~\ref{lem:kadison-matrix}).
\item We then prove the existence of a finite-dimensional abelian subalgebra $D\subset M$ such that
$$\|\theta_j|_{D'\cap P}\|<\|\theta_j|_{F'\cap P}\|+\varepsilon/2<\varepsilon\qquad(1\leq j\leq s)$$
(Lemma~\ref{lem:kadison-haar} and Proposition~\ref{prop:kadison-abelian}).
\end{enumerate}

We use our Poisson boundary theory to construct $M$-inner ucp maps $R_n$ with finitely many coefficients such that $\|\theta\circ R_n\|\to0$ for every $\theta\in \mathcal Z(M)_\perp$. We will use these maps to construct a finite-dimensional matrix algebra $F\subset M$ with $\|\theta|_{F'\cap P}\|<\varepsilon$.

\begin{lemma}\label{lem:kadison-finite-averages}
Let $M\subset P$ be as in Theorem~\ref{letterthm kadison}. Then there is a sequence of normal $M$-inner ucp maps $R_n:P\to P,\ a\mapsto\sum_{k=1}^{r(n)}c_{n,k}ac_{n,k}^*$, each with finitely many coefficients $c_{n,k}\in M$, such that
$$\|\theta\circ R_n\|\longrightarrow0\qquad(n\to\infty)$$
for every $\theta\in \mathcal Z(M)_\perp=\{\theta\in P_*\mid\theta|_{\mathcal Z(M)}=0\}$.
\end{lemma}

\begin{proof}
Fix a normal faithful state $\phi$ on $M$. Take an nrsg $(M,\phi)$-hyperstate $\varphi$, which always exists by Example~\ref{exists nrsg hyperstate}, with inner ucp coefficients $(y_k)\subset M$. Define the normal $M$-inner ucp map
\begin{equation}\label{eq:kadison-inner-map}
\CQ(x)=\sum_k y_kxy_k^*\qquad(x\in P).
\end{equation}
By Proposition~\ref{double ergodicity}, $\har(\CQ,P)=M'\cap P=\mathcal Z(M)$.

Since $P_*$ is separable, take a norm-dense sequence $(\theta_j)_{j\geq1}$ in $\mathcal Z(M)_\perp$. Let $\mathcal C_n=(1/n)\sum_{\ell=0}^{n-1}\CQ^\ell$. Let $T=\CQ_*$ be the predual map on $P_*$, so that $T\theta=\theta\circ \CQ$. The annihilator of $\operatorname{Ran}(\id-T)$ in $P=(P_*)^*$ is $\har(\CQ,P)=\mathcal Z(M)$. Hence by the Hahn--Banach theorem, we have
$$\overline{\operatorname{Ran}(\id-T)}^{\|\,\cdot\,\|} =\{\theta\in P_*\mid\theta|_{\mathcal Z(M)}=0\}.$$
For $\theta=(\id-T)\eta$ with $\eta\in P_*$, we have
$$\|\theta\circ \mathcal C_n\| =\frac{1}{n}\|\eta-T^n\eta\|\leq\frac{2\|\eta\|}{n}.$$
Since $\|\mathcal C_n\|=1$, this convergence extends to the norm closure of $\operatorname{Ran}(\id-T)$. In particular, $\|\theta_j\circ \mathcal C_n\|\to0$ as $n\to\infty$ for every $j\geq1$.

Each $\mathcal C_n$ is a normal $M$-inner ucp map. Write $\mathcal C_n(x)=\sum_{i\geq1}a_{n,i}xa_{n,i}^*$ $(x\in P)$ with $a_{n,i}\in M$ and $\sum_i a_{n,i}a_{n,i}^*=1$. For $m\geq1$, define
$$b_{n,m}=\left(1-\sum_{i=1}^m a_{n,i}a_{n,i}^*\right)^{1/2},\qquad R_{n,m}(x)=\sum_{i=1}^m a_{n,i}xa_{n,i}^*+b_{n,m}xb_{n,m}.$$
Then $R_{n,m}$ is an $M$-inner ucp map with finitely many coefficients. For $\theta\in P_*^+$, the functional $\theta\circ R_{n,m}-\theta\circ \mathcal C_n$ is the difference of the positive functionals $\theta(b_{n,m}\,\cdot\,b_{n,m})$ and $\sum_{i>m}\theta(a_{n,i}\,\cdot\,a_{n,i}^*)$, both of norm $\theta(b_{n,m}^2)$. Hence, for fixed $n$, we have
$$\|\theta\circ R_{n,m}-\theta\circ \mathcal C_n\| \leq2\theta(b_{n,m}^2)\longrightarrow0\qquad(m\to\infty).$$
Every normal functional is a linear combination of positive normal functionals, so the convergence holds for every $\theta\in P_*$. For each $n\geq1$, take $m(n)$ sufficiently large and let $R_n=R_{n,m(n)}$ so that
$$\|\theta_j\circ R_n-\theta_j\circ \mathcal C_n\|<2^{-n} \qquad(1\leq j\leq n).$$
For every fixed $j$ and $n\geq j$, we have
$$\|\theta_j\circ R_n\|\leq2^{-n}+\|\theta_j\circ \mathcal C_n\| \longrightarrow0\qquad(n\to\infty).$$
Since $\|R_n\|=1$, norm density gives the same convergence for every $\theta\in \mathcal Z(M)_\perp$.
\end{proof}

We use the proper infiniteness of $M$ and the maps $R_n$ constructed above to construct a full matrix algebra $F\subset M$ with $\|\theta|_{F'\cap P}\|<\varepsilon$.

\begin{proposition}\label{lem:kadison-matrix}
Let $M\subset P$ be as in Theorem~\ref{letterthm kadison}. Assume that $M$ is properly infinite. Let $\theta_1,\ldots,\theta_s\in \mathcal Z(M)_\perp$, and let $\varepsilon>0$. Then there is a full matrix algebra $F\subset M$ such that
$$\|\theta_j\big|_{F'\cap P}\|<\varepsilon \qquad(1\leq j\leq s).$$
\end{proposition}

\begin{proof}
By Lemma~\ref{lem:kadison-finite-averages}, take a normal $M$-inner ucp map $R(x)=\sum_{i=1}^r a_ixa_i^*$ $(x\in P)$ with finitely many coefficients $a_i\in M$ such that $\|\theta_j\circ R\|<\varepsilon/2$ for $1\leq j\leq s$. We use the same matrix construction as in the proof of Theorem~\ref{Thm inner UCP eq to Dixmier}. Let $d=r+1$, and fix a projection $e\in M$ with $e\sim1-e\sim1$. Take matrix units $(e_{ij})_{1\leq i,j\leq d}$ in $\mathrm M_d(\C)\otimes1$. Take an isometry $v_e\in M$ with $v_ev_e^*=e$ and an isomorphism
$$\pi_e:M\xrightarrow{\sim}\mathrm M_d(\C)\bar\otimes M, \qquad \pi_e(e)=e_{11},$$
which extends to an isomorphism $P\xrightarrow{\sim}\mathrm M_d(\C)\bar\otimes P$, still denoted by $\pi_e$, such that
$$
\pi_e^{-1}\begin{pmatrix}
x&0_{1\times r}\\
0_{r\times1}&0_{r\times r}
\end{pmatrix}=v_exv_e^*
\qquad(x\in P).
$$
In $\mathrm M_d(\C)\bar\otimes M$, let
$$
Z=\begin{pmatrix}
0&a_1&\cdots&a_r\\
0&0&\cdots&0\\
\vdots&\vdots&&\vdots\\
0&0&\cdots&0
\end{pmatrix}.
$$
Since $\sum_{i=1}^ra_ia_i^*=1$ and the first column of $Z$ is $0$,
$$ZZ^*=e_{11},\qquad 1-Z^*Z\geq e_{11}\sim1,\qquad 1-e_{11}\sim1.$$
Take a partial isometry $W\in\mathrm M_d(\C)\bar\otimes M$ with
$$W^*W=1-Z^*Z,\qquad WW^*=1-e_{11}.$$
Then $U=Z+W$ is unitary and its first row is $(0,a_1,\ldots,a_r)$. For $x\in P$, by matrix multiplication, we have
$$
U\begin{pmatrix}
x&0&\cdots&0\\
0&x&\cdots&0\\
\vdots&\vdots&\ddots&\vdots\\
0&0&\cdots&x
\end{pmatrix}U^*
=\begin{pmatrix}R(x)&\ast\\\ast&\ast\end{pmatrix}.
$$
Define
$$\alpha_e(x)=\pi_e^{-1}\bigl(U(1_d\otimes x)U^*\bigr),\qquad F_e=\pi_e^{-1}\bigl(U(\mathrm M_d(\C)\otimes1)U^*\bigr).$$
Then $\alpha_e$ is a normal unital injective endomorphism, $F_e\subset M$ is a full matrix algebra, and
$$\alpha_e(P)=F_e'\cap P.$$
Since $\pi_e(e)=e_{11}$, we have
\begin{equation}\label{eq:kadison-endomorphism-compression}
e\alpha_e(x)e
=\pi_e^{-1}\begin{pmatrix}R(x)&0_{1\times r}\\0_{r\times1}&0_{r\times r}\end{pmatrix}
=v_eR(x)v_e^*
\qquad(x\in P).
\end{equation}

Finally, take projections $e_n\nearrow1$ strongly in $M$, with $e_0=0$ and $e_n-e_{n-1}\sim1-e_n\sim1$ for $n\geq1$. Take partial isometries $v_n\in M$ with $v_n^*v_n=1-e_{n-1}$ and $v_nv_n^*=e_n-e_{n-1}$, and use $v_{e_n}=e_{n-1}+v_n$ in the construction above. Since $v_{e_n}e_{n-1}=v_{e_n}^*e_{n-1}=e_{n-1}$, we have $v_{e_n}\to1$ $\ast$-strongly. Let $\alpha_n=\alpha_{e_n}$ and $F_n=F_{e_n}$. For $\theta\in P_*$, by (\ref{eq:kadison-endomorphism-compression}), we have
$$
\|\theta\circ\alpha_n-\theta\circ R\|
\leq\|\theta-\theta(e_n\,\cdot\,e_n)\|
    +\|\theta(v_{e_n}\,\cdot\,v_{e_n}^*)-\theta\|
\longrightarrow0\qquad(n\to\infty).
$$
The convergence follows first for normal states by the Cauchy--Schwarz inequality and then for every $\theta\in P_*$ by linearity. Take $n$ sufficiently large that $\|\theta_j\circ\alpha_n-\theta_j\circ R\|<\varepsilon/2$ for $1\leq j\leq s$. Since $\alpha_n$ is a $\ast$-isomorphism of $P$ onto $F_n'\cap P$, we have
$$\|\theta_j|_{F_n'\cap P}\|=\|\theta_j\circ\alpha_n\|
\leq\|\theta_j\circ\alpha_n-\theta_j\circ R\|
+\|\theta_j\circ R\|<\varepsilon\qquad(1\leq j\leq s).$$
Taking $F=F_n$ finishes the proof.
\end{proof}

For a finite-dimensional abelian subalgebra $D\subset P$, let $\operatorname{At}(D)$ denote the set of minimal projections of $D$, and define
$$\Phi_D(x)=\sum_{p\in\operatorname{At}(D)}pxp\qquad(x\in P).$$
For a finite-dimensional subalgebra $F\subset P$, define
$$\Psi_F(x)=\int_{\mathcal U(F)}uxu^*\,\d u\qquad(x\in P),$$
where $\d u$ is the normalized Haar measure on $\mathcal U(F)$. Then $\Phi_D$ and $\Psi_F$ are the normal conditional expectations onto $D'\cap P$ and $F'\cap P$, respectively. In particular,
$$\|\theta|_{D'\cap P}\|=\|\theta\circ\Phi_D\|,\qquad \|\theta|_{F'\cap P}\|=\|\theta\circ\Psi_F\|\qquad(\theta\in P_*).$$

The following simple lemma will help us pass from the matrix algebra $F$ to an abelian subalgebra $D$ with $\|\theta|_{D'\cap P}\|<\varepsilon$.

\begin{lemma}\label{lem:kadison-haar}
Let $M\subset P$ be as in Theorem~\ref{letterthm kadison}. Let $F\subset P$ be a full matrix algebra with $F\cong\mathrm M_d(\C)$, and let $D\subset F$ be a maximal abelian subalgebra. Then
\begin{equation}\label{eq:kadison-haar}
\int_{\mathcal U(F)}\Phi_{uDu^*}\,\d u=\frac{\id+d\Psi_F}{d+1}.
\end{equation}
If $\rho$ is a normal faithful state on $P$ with $F\subset P_\rho$, then
\begin{equation}\label{eq:kadison-haar-L2}
\int_{\mathcal U(F)}\|\Phi_{uDu^*}(b)\|_{2,\rho}^2\,\d u
 =\frac{\|b\|_{2,\rho}^2}{d+1}
\end{equation}
for every $b\in P$ with $\Psi_F(b)=0$.
\end{lemma}

\begin{proof}
Under the identification $P=F\bar\otimes(F'\cap P)$, both sides of (\ref{eq:kadison-haar}) act only on the first tensor factor. Hence we only need to prove the identity on $F\cong\mathrm M_d(\C)$. The case $d=1$ is immediate. Let $\Phi=\int_{\mathcal U(F)}\Phi_{uDu^*}\,\d u$. Let $p\in F$ be a rank-one projection. For any $u\in\mathcal U(\{p\}'\cap F)$, invariance of Haar measure gives
$$u\Phi(p)u^*=\Phi(upu^*)=\Phi(p).$$
Hence $\Phi(p)\in(\{p\}'\cap F)'\cap F=\{p\}''\cap F=\C p+\C(1-p)$. Since all rank-one projections are unitarily conjugate and span $F$, we have $\Phi=\lambda\id+\mu\operatorname{tr}_d(\,\cdot\,)1$ for some scalars $\lambda,\mu$, where $\operatorname{tr}_d$ is the normalized trace on $F$. Each $\Phi_{uDu^*}$ is a rank-$d$ orthogonal projection in $B(L^2(F,\operatorname{tr}_d))$. Taking the operator trace $\operatorname{Tr}$ on $B(L^2(F,\operatorname{tr}_d))$ and using $\Phi(1)=1$, we have
$$d^2\lambda+\mu=d,\qquad \lambda+\mu=1.$$
Hence $\lambda=1/(d+1)$ and $\mu=d/(d+1)$. Tensoring with the identity on $F'\cap P$ proves (\ref{eq:kadison-haar}).

Assume now that $F\subset P_\rho$. Then each $\Phi_{uDu^*}$ is the $\rho$-preserving conditional expectation onto $(uDu^*)'\cap P$, and hence acts as an orthogonal projection on $L^2(P,\rho)$. For $b\in P$ with $\Psi_F(b)=0$, by (\ref{eq:kadison-haar}) we have
$$
\int_{\mathcal U(F)}\|\Phi_{uDu^*}(b)\|_{2,\rho}^2\,\d u =\int_{\mathcal U(F)} \langle \Phi_{uDu^*}(b),b\rangle_{2,\rho}\,\d u =\frac{\langle b+d\Psi_F(b),b\rangle_{2,\rho}}{d+1} =\frac{\|b\|_{2,\rho}^2}{d+1}.
$$
This proves (\ref{eq:kadison-haar-L2}).
\end{proof}

We now construct the required finite-dimensional abelian subalgebra when $M$ is continuous and properly infinite.

\begin{proposition}\label{prop:kadison-abelian}
Let $M\subset P$ be as in Theorem~\ref{letterthm kadison}. Assume that $M$ is continuous and properly infinite. Let $\theta_1,\ldots,\theta_s\in M_\perp$, and let $\varepsilon>0$. Then there is a finite-dimensional abelian subalgebra $D\subset M$ such that
$$\|\theta_j|_{D'\cap P}\|=\|\theta_j\circ \Phi_D\|<\varepsilon\qquad(1\leq j\leq s).$$
\end{proposition}

\begin{proof}
By Proposition~\ref{lem:kadison-matrix}, take a full matrix algebra $F\subset M$ such that $\|\theta_j|_{F'\cap P}\|<\varepsilon/2$ for $1\leq j\leq s$. Let $F\cong\mathrm M_r(\C)$ and $N=F'\cap M$. Then $M=F\bar\otimes N$, and $N$ is a continuous von Neumann algebra. By \cite[Theorem~11.1]{HS90}, take a normal faithful state $\chi$ on $N$ such that $N_\chi$ is of type $\mathrm{II}_1$ and hence contains a copy of $\mathrm M_k(\C)$ for every $k\geq1$. Indeed, there exists a projection $p\in N_\chi$ with $\chi(p)=1/k$ such that $N_\chi\cong\mathrm M_k(\C)\bar\otimes pN_\chi p$. Let $\operatorname{tr}_r$ be the normalized trace on $\mathrm M_r(\C)$, and let
$$\psi=\operatorname{tr}_r\otimes\chi\in M_*,\qquad\rho=\psi\circ E\in P_*.$$

Since $\rho$ is faithful, the only $x\in P$ with $\rho(ax)=0$ for all $a\in P$ is $x=0$. By the Hahn--Banach theorem,
\begin{equation}\label{eq:faithful-density}
\overline{\{\rho(a\,\cdot\,)\mid a\in P\}}^{\|\,\cdot\,\|}=P_*.
\end{equation}

By (\ref{eq:faithful-density}), we can take $a_1,\ldots,a_s\in P$ and $k\geq1$ such that
$$
\|\theta_j-\rho(a_j\,\cdot\,)\|<\frac{\varepsilon}{8},\qquad
\left(\frac{4\sum_{j=1}^s\|a_j\|^2}{rk+1}\right)^{1/2}
<\frac{\varepsilon}{4}.
$$
Let $K\subset N_\chi$ be a subalgebra with $K\cong\mathrm M_k(\C)$, and let $F_1=F\vee K$. Then
$$F_1\cong\mathrm M_{rk}(\C),\qquad F_1\subset M_\psi\subset P_\rho.$$

Let $c_j=a_j-\Psi_F(a_j)$ for $1\leq j\leq s$. Since $\mathcal U(F)\subset\mathcal U(F_1)$, by the invariance of Haar measure we have $\Psi_{F_1}\Psi_F=\Psi_{F_1}$, and hence $\Psi_{F_1}(c_j^*)=0$. Fix a maximal abelian subalgebra $A_1\subset F_1$. Since $\rho(1)=1$ and $\|\Psi_F\|=1$, we have
$$\|c_j^*\|_{2,\rho}\leq\|c_j\| \leq\|a_j\|+\|\Psi_F(a_j)\|\leq2\|a_j\|.$$
Using $F_1\cong\mathrm M_{rk}(\C)$, $F_1\subset P_\rho$ and $\Psi_{F_1}(c_j^*)=0$, by (\ref{eq:kadison-haar-L2}) we have
$$
\int_{\mathcal U(F_1)} \sum_{j=1}^s\|\Phi_{uA_1u^*}(c_j^*)\|_{2,\rho}^2\,\d u =\sum_{j=1}^s\int_{\mathcal U(F_1)} \|\Phi_{uA_1u^*}(c_j^*)\|_{2,\rho}^2\,\d u =\frac{\sum_{j=1}^s\|c_j^*\|_{2,\rho}^2}{rk+1} \leq\frac{4\sum_{j=1}^s\|a_j\|^2}{rk+1} <\frac{\varepsilon^2}{16}.
$$
Since $\d u$ has total mass $1$, there is $u\in\mathcal U(F_1)$ for which the sum is no greater than its average. Hence the finite-dimensional abelian subalgebra $D=uA_1u^*\subset M$ satisfies
$$\sum_{j=1}^s\|\Phi_D(c_j^*)\|_{2,\rho}^2 \leq\frac{\sum_{j=1}^s\|c_j^*\|_{2,\rho}^2}{rk+1} <\frac{\varepsilon^2}{16}.$$

Since $F,D\subset P_\rho$, the conditional expectations $\Psi_F$ and $\Phi_D$ preserve $\rho$, and by bimodularity we have
$$\rho(c_j\,\cdot\,)=\rho(a_j\,\cdot\,)\circ(\id-\Psi_F),\qquad \rho(c_j\,\cdot\,)\circ \Phi_D=\rho(\Phi_D(c_j)\,\cdot\,).$$
Hence
\begin{equation}\label{eq:kadison-total-error}
\theta_j\circ \Phi_D =\theta_j\circ \Psi_F\circ \Phi_D+\rho(\Phi_D(c_j)\,\cdot\,)+(\theta_j-\rho(a_j\,\cdot\,))\circ(\id-\Psi_F)\circ \Phi_D.
\end{equation}
By (\ref{eq:kadison-total-error}), the triangle inequality, and the Cauchy--Schwarz inequality,
$$
\|\theta_j|_{D'\cap P}\|=\|\theta_j\circ \Phi_D\| \leq\|\theta_j\circ \Psi_F\|+\|\Phi_D(c_j^*)\|_{2,\rho} +2\|\theta_j-\rho(a_j\,\cdot\,)\| <\|\theta_j\big|_{F'\cap P}\|+\frac{\varepsilon}{2}<\varepsilon \qquad(1\leq j\leq s).
$$
\end{proof}

\subsection{The finite case.}\label{sec:kadison-finite-case}

For finite $M$, the proof of the following proposition is essentially Popa's averaging and refinement argument \cite{Popa1981Kadison}. We include it for completeness, using the tracial case of Example~\ref{exists nrsg hyperstate} to express the unitary averaging through an inner ucp map.

\begin{proposition}\label{lem:kadison-finite}
Let $M\subset P$ be as in Theorem~\ref{letterthm kadison}. Assume that $M$ is finite. Let $\theta_1,\ldots,\theta_s\in M_\perp$, and let $\varepsilon>0$. Then there is a finite-dimensional abelian subalgebra $D\subset M$ such that
$$\|\theta_j|_{D'\cap P}\|=\|\theta_j\circ \Phi_D\|<\varepsilon\qquad(1\leq j\leq s).$$
\end{proposition}

\begin{proof}
Let $\tau$ be a normal faithful tracial state on $M$, and let $\rho=\tau\circ E$. Since $M_\tau=M$, we have $M\subset P_\rho$.

We first reduce the proposition to an $L^2$ estimate. Since $\rho$ is faithful, we still have (\ref{eq:faithful-density}). Hence we can take $a_1,\ldots,a_s\in P$ such that
$$\|\theta_j-\rho(a_j^*\,\cdot\,)\|<\frac{\varepsilon}{8}.$$
Let $b_j=a_j-E(a_j)$ and $\eta_j=\rho(b_j^*\,\cdot\,)$. Since $\rho\circ E=\rho$ and $E$ is $M$-bimodular,
$$\rho(a_j^*E(x))=\rho(E(a_j)^*x)\qquad(x\in P).$$
Using $\theta_j\circ E=0$ and $\|\id-E\|\leq2$, we have
$$E(b_j)=0,\qquad \theta_j-\eta_j=(\theta_j-\rho(a_j^*\,\cdot\,))\circ(\id-E),\qquad \|\theta_j-\eta_j\|<\frac{\varepsilon}{4}.$$
Let $D\subset M$ be a finite-dimensional abelian subalgebra. Since $D\subset M\subset P_\rho$, we have
$$\eta_j(\Phi_D(x)) =\sum_{p\in\operatorname{At}(D)}\rho(b_j^*pxp) =\rho(\Phi_D(b_j)^*x)\qquad(x\in P).$$
By the Cauchy--Schwarz inequality, we have
$$\|\rho(\Phi_D(b_j)^*\,\cdot\,)\| \leq\|\Phi_D(b_j)\|_{2,\rho}.$$
Hence, by the triangle inequality and $\|\Phi_D\|=1$,
$$\|\theta_j\circ \Phi_D\| \leq\|\theta_j-\eta_j\|+\|\rho(\Phi_D(b_j)^*\,\cdot\,)\| <\frac{\varepsilon}{4}+\|\Phi_D(b_j)\|_{2,\rho}.$$
Let $\mathcal D$ be the family of finite-dimensional abelian subalgebras of $M$, and let
$$t=\inf_{D\in\mathcal D}\sum_{j=1}^s\|\Phi_D(b_j)\|_{2,\rho}^2.$$
We only need to prove that $t=0$, since then some $D\in\mathcal D$ satisfies $\|\Phi_D(b_j)\|_{2,\rho}<\varepsilon/2$ for $1\leq j\leq s$.

Fix $\delta>0$, and take $D\in\mathcal D$ such that
$$\sum_{j=1}^s\|\Phi_D(b_j)\|_{2,\rho}^2<t+\delta.$$
Let $c_j=\Phi_D(b_j)$ and $N=D'\cap M$. By Lemma~\ref{lem:corners}(ii), we have $N'\cap P=D\vee\mathcal Z(M)=\mathcal Z(N)\subset M$. Moreover, by the $M$-bimodularity of $E$, we have
$$E(c_j)=\Phi_D(E(b_j))=0.$$

The map $\Phi_D\circ E:P\to N$ is a normal faithful conditional expectation. Using Example~\ref{exists nrsg hyperstate}(iii) for $(N,\tau|_N)$, take a normal $N$-inner ucp map $\CQ:P\to P$ whose coefficients are scalar multiples of unitaries. By the first part of the proof of Lemma~\ref{lem:kadison-finite-averages}, its Ces\`aro means $\mathcal C_n=(1/n)\sum_{\ell=0}^{n-1}\CQ^\ell$ satisfy $\|\theta\circ \mathcal C_n\|\to0$ as $n\to\infty$ for every $\theta\in N_\perp=\{\theta\in P_*\mid\theta|_N=0\}$. Each $\mathcal C_n$ is a countable convex combination of conjugations by unitaries of $N$. Since $E(c_j)=0$, the functionals $\rho(c_j^*\,\cdot\,)$ vanish on $N$. Hence we have
$$\sum_{j=1}^s\operatorname{Re}\rho(c_j^*\mathcal C_n(c_j)) \longrightarrow0\qquad(n\to\infty).$$
For sufficiently large $n$, one of the unitaries $u\in\mathcal U(N)$ in this convex combination satisfies
$$\sum_{j=1}^s\operatorname{Re}\rho(c_j^*uc_ju^*)<\delta.$$
Approximating $u$ in norm by a finite-spectrum unitary $v\in W^*(u)\subset N$, we may also arrange
$$\sum_{j=1}^s\operatorname{Re}\rho(c_j^*vc_jv^*)<\delta.$$
Let $B=D\vee W^*(v)$. Since $D$ is finite-dimensional and $v\in D'\cap M$ has finite spectrum, $B$ is a finite-dimensional abelian subalgebra of $M$ containing $D$, so $B\in\mathcal D$. Since $M\subset P_\rho$, the map $\Phi_B$ is a contraction on $L^2(P,\rho)$ and $\|vc_jv^*\|_{2,\rho}=\|c_j\|_{2,\rho}$. Using $\Phi_B\Phi_D=\Phi_B$ and $\Phi_B\circ\AD(v)=\Phi_B$, we have
$$
t\leq\sum_j\|\Phi_B(b_j)\|_{2,\rho}^2 =\sum_j\left\|\Phi_B\left(\frac{c_j+vc_jv^*}{2}\right)\right\|_{2,\rho}^2 \leq\frac{1}{2}\sum_j\|c_j\|_{2,\rho}^2 +\frac{1}{2}\sum_j\operatorname{Re}\rho(c_j^*vc_jv^*) <\frac{t}{2}+\delta.
$$
Letting $\delta\to0$, we have $t=0$, which proves the proposition.
\end{proof}

\subsection{Proof of Theorem~\ref{letterthm kadison}.}\label{sec:kadison-proof}

We are ready to put everything together.

\begin{proposition}\label{lem:kadison-refinement}
Let $M\subset P$ be as in Theorem~\ref{letterthm kadison}. Assume that $M$ is continuous, i.e., $M=M_{\mathrm{II}}\oplus M_{\mathrm{III}}$. Let $D_0\subset M$ be a finite-dimensional abelian subalgebra, let $\theta_1,\ldots,\theta_s\in M_\perp$, and let $\varepsilon>0$. Then there is a finite-dimensional abelian subalgebra $D_1\subset M$ such that
$$D_0\subset D_1\qquad\text{and}\qquad \|\theta_j|_{D_1'\cap P}\|<\varepsilon\quad(1\leq j\leq s).$$
\end{proposition}

\begin{proof}
Let $e_1,\ldots,e_m$ be the minimal projections of $D_0$. For each $k$, take a central projection $z_k\in\mathcal Z(e_kMe_k)$ such that
$$e_kMe_k=z_k(e_kMe_k)\oplus(1-z_k)(e_kMe_k),$$
with $z_k(e_kMe_k)$ finite and $(1-z_k)(e_kMe_k)$ properly infinite. Replacing $D_0$ with $D_0\vee\{z_k\}''$ and relabelling, we may assume that each $e_kMe_k$ is either finite or properly infinite. By Lemma~\ref{lem:corners}(i), each inclusion $e_kMe_k\subset e_kPe_k$ satisfies the hypotheses of Theorem~\ref{letterthm kadison}, and the restrictions $\theta_j|_{e_kPe_k}$ vanish on $e_kMe_k$. For each $k$, apply Proposition~\ref{lem:kadison-finite} if $e_kMe_k$ is finite and Proposition~\ref{prop:kadison-abelian} if it is properly infinite, with error $\varepsilon/m$. We have a finite-dimensional abelian subalgebra $D_{1,k}\subset e_kMe_k$ such that
$$\|\theta_j|_{D_{1,k}'\cap e_kPe_k}\| <\frac{\varepsilon}{m}\qquad(1\leq j\leq s).$$
Let $D_1=\bigoplus_{k=1}^m D_{1,k}$. Then $D_0\subset D_1\subset M$. By Lemma~\ref{lem:corners}(iii), we have
$$\|\theta_j|_{D_1'\cap P}\| \leq\sum_{k=1}^m\|\theta_j|_{D_{1,k}'\cap e_kPe_k}\| <\varepsilon\qquad(1\leq j\leq s).$$
\end{proof}

\begin{proof}[Proof of Theorem~\ref{letterthm kadison}]
Write $M=(\bigoplus_{1\leq n\leq\infty}M_{\mathrm I_n})\oplus M_{\mathrm c}$, where $M_{\mathrm I_n}$ is homogeneous of type $\mathrm I_n$ and $M_{\mathrm c}$ is continuous. By Lemma~\ref{lem:corners}(i), we only need to prove the homogeneous type $\mathrm I$ case and the continuous case. In the former case, write $M=B(\mathcal{H})\bar\otimes\mathcal Z(M)$. Let $P_0=\mathcal Z(M)'\cap P$. Since $B(\mathcal{H})'\cap P_0=M'\cap P=\mathcal Z(M)$, we have
$$P_0=B(\mathcal{H})\bar\otimes(B(\mathcal{H})'\cap P_0) =B(\mathcal{H})\bar\otimes\mathcal Z(M)=M.$$
For any maximal abelian subalgebra $A\subset M$, we then have $A'\cap P\subset\mathcal Z(M)'\cap P=M$, and hence $A'\cap P=A$.

Now assume that $M$ is continuous. Since $P_*$ is separable, take a norm-dense sequence $(\theta_j)_{j\geq1}$ in $M_\perp$. Starting with $D_0=\C1$ and applying Proposition~\ref{lem:kadison-refinement} to $\theta_1,\ldots,\theta_n$ with $\varepsilon=2^{-n}$ at the $n$th step, we have finite-dimensional abelian subalgebras $D_n\subset M$ such that
$$D_{n-1}\subset D_n,\qquad \|\theta_j|_{D_n'\cap P}\|<2^{-n}\qquad(1\leq j\leq n).$$
Let $A_0=\bigvee_nD_n$. For every $j\geq1$, we have
$$\|\theta_j|_{A_0'\cap P}\|\leq\|\theta_j|_{D_n'\cap P}\|\longrightarrow0\qquad(n\to\infty).$$
By norm density and the Hahn--Banach theorem, $A_0'\cap P\subset M$. Let $A$ be a maximal abelian subalgebra of $M$ containing $A_0$. Then
$$A'\cap P\subset A_0'\cap P\subset M,$$
and hence $A'\cap P=A'\cap M=A$, which finishes the proof.
\end{proof}

\subsection{\texorpdfstring{$\varphi$}{varphi}-harmonic ergodicity and AFD ergodicity.}\label{sec:harmonic-afd}

Let $M\subset P$ be an inclusion of von Neumann algebras with separable preduals such that $M'\cap P=\mathcal Z(M)$. Following the terminology for ergodic embeddings of Popa \cite[Definition~2.2]{Popa2021Ergodic} and the notion of AFD-ergodicity used by Marrakchi \cite{Ma26}, we use the following terminology, with $\mathcal Z(M)$ in place of $\C1$ when $M$ is not a factor. We say that $M\subset P$ is \textbf{AFD-ergodic} if there is a hyperfinite von Neumann subalgebra $N\subset M$ such that $N'\cap P=\mathcal Z(M)$. We say that $M\subset P$ is \textbf{MV-ergodic} if it has the weak Dixmier property. Clearly, AFD-ergodicity implies MV-ergodicity by averaging over the unitary groups of an increasing sequence of finite-dimensional subalgebras generating $N$.

We also use Poisson boundaries to define the following notion of ergodicity.

\begin{definition}\label{def:harmonic-ergodicity}
Let $M\subset P$ be an inclusion of von Neumann algebras with separable preduals, and let $\varphi$ be an nrsg $M$-hyperstate with Poisson boundary $\CB_\varphi$. We say that $M\subset P$ has \textbf{$\varphi$-harmonic ergodicity} if there is a normal faithful ucp map $\Phi:P\to\CB_\varphi$ with $\Phi|_M=\id_M$.
\end{definition}

By Theorem~\ref{Thm ucp to PB}(i), we have already proved that $\varphi$-harmonic ergodicity implies MV-ergodicity when $M$ is properly infinite.

In the case $M'\cap P=\mathcal Z(M)$, Popa \cite[Theorem~3.2 and Remark~3.5(2)]{Popa1981Kadison} proved that, when $M$ is a semifinite factor, the existence of a normal faithful conditional expectation $E:P\to M$ implies AFD-ergodicity. Marrakchi \cite[Theorem~B]{Ma26} proved this implication for general $M$. As explained in Subsection~\ref{sec:intro-kadison}, this is their method for constructing finite-dimensional subalgebras $F_n\subset M$ with $\|\theta|_{F_n'\cap P}\|\to0$ for every $\theta\in\mathcal Z(M)_\perp$.

We strengthen these results when $M$ is properly infinite by proving that $\varphi$-harmonic ergodicity directly implies AFD-ergodicity. We first strengthen Proposition~\ref{lem:kadison-matrix} by passing from a single full matrix algebra $F\subset M$ with $\|\theta_j|_{F'\cap P}\|<\varepsilon$ to an increasing sequence.

\begin{proposition}\label{lem:harmonic-matrix-refinement}
Let $M\subset P$ be an inclusion of von Neumann algebras with separable preduals, and assume that $M$ is properly infinite. Let $\varphi$ be an nrsg $M$-hyperstate, and assume that $M\subset P$ has $\varphi$-harmonic ergodicity. Write $(M'\cap P)_\perp=\{\theta\in P_*\mid\theta|_{M'\cap P}=0\}$. For every full matrix algebra $F_0\subset M$, every $\theta_1,\ldots,\theta_s\in(M'\cap P)_\perp$ and every $\varepsilon>0$, there is a full matrix algebra $F\subset M$ such that
$$F_0\subset F\qquad\text{and}\qquad\|\theta_j|_{F'\cap P}\|<\varepsilon\qquad(1\leq j\leq s).$$
\end{proposition}

\begin{proof}
Let $\CQ:P\to P$ be the normal $M$-inner ucp map $\CQ(x)=\sum_k y_kxy_k^*$ associated with $\varphi$, as in (\ref{eq:kadison-inner-map}). By Proposition~\ref{double ergodicity}, we have $\har(\CQ,P)=M'\cap P$. The proofs of Lemma~\ref{lem:kadison-finite-averages} and Proposition~\ref{lem:kadison-matrix} use the conditional expectation only to establish this identity. Hence the same arguments apply to functionals $\theta$ with $\theta|_{M'\cap P}=0$. For each $n\geq1$, we apply the matrix construction of Proposition~\ref{lem:kadison-matrix} with error $2^{-n}$ to the first $n$ elements of a norm-dense sequence in $(M'\cap P)_\perp$. By norm density, the resulting full matrix algebras $\mathcal F_n\subset M$ satisfy
$$\|\theta|_{\mathcal F_n'\cap P}\|\longrightarrow0\qquad(n\to\infty)$$
for every $\theta\in(M'\cap P)_\perp$.

Let $M_0=F_0'\cap M$ and $P_0=F_0'\cap P$, so that $M=F_0\bar\otimes M_0$ and $P=F_0\bar\otimes P_0$. Take matrix units $(e_{ij})_{1\leq i,j\leq d}$ for $F_0$. Since $M$ is properly infinite, $e_{11}\sim1$. Take an isometry $v\in M$ with $vv^*=e_{11}$. Then
$$\pi:P\xrightarrow{\sim}P_0\qquad\text{with } \pi(x)=\sum_{i=1}^d e_{i1}vxv^*e_{1i}$$
is a normal $*$-isomorphism with $\pi(M)=M_0$. Since $v,e_{ij}\in M$, we have $\pi(z)=z$ for every $z\in M'\cap P$. Hence $\theta_j\circ\pi\in(M'\cap P)_\perp$. Take $n$ sufficiently large that $\|(\theta_j\circ\pi)|_{\mathcal F_n'\cap P}\|<\varepsilon$ for $1\leq j\leq s$, and let $F=F_0\vee\pi(\mathcal F_n)\cong F_0\bar\otimes\pi(\mathcal F_n)$. Then $F$ is a full matrix algebra containing $F_0$, and $F'\cap P=\pi(\mathcal F_n'\cap P)$. Therefore,
$$\|\theta_j|_{F'\cap P}\| =\|(\theta_j\circ\pi)|_{\mathcal F_n'\cap P}\| <\varepsilon\qquad(1\leq j\leq s),$$
which finishes the proof.
\end{proof}

\begin{theorem}[Theorem~\ref{letterthm harmonic afd}]\label{thm:harmonic-implies-afd}
Let $M\subset P$ be an inclusion of von Neumann algebras with separable preduals, and assume that $M$ is properly infinite. Let $\varphi$ be an nrsg $M$-hyperstate. If there is a normal faithful ucp map $\Phi:P\to\CB_\varphi$ with $\Phi|_M=\id_M$, then there is a hyperfinite von Neumann subalgebra $N\subset M$ such that $N'\cap P=M'\cap P$.
\end{theorem}

\begin{proof}
Take a norm-dense sequence $(\theta_j)_{j\geq1}$ in $(M'\cap P)_\perp$. Starting with $F_0=\C1$, apply Proposition~\ref{lem:harmonic-matrix-refinement} recursively. We have full matrix algebras $F_n\subset M$ such that
$$F_{n-1}\subset F_n,\qquad \|\theta_j|_{F_n'\cap P}\|<2^{-n}\qquad(1\leq j\leq n).$$
Let $N=\bigvee_nF_n$. Then $N$ is hyperfinite. For every $j\geq1$, we have
$$\|\theta_j|_{N'\cap P}\|\leq\|\theta_j|_{F_n'\cap P}\|\longrightarrow0\qquad(n\to\infty).$$
By norm density and the Hahn--Banach theorem, $N'\cap P\subset M'\cap P$. Since $N\subset M$, the reverse inclusion is immediate. Hence $N'\cap P=M'\cap P$.
\end{proof}

\bibliographystyle{alphaurl}
\bibliography{refv9}
\end{document}